# Some trigonometric integrals involving $\log\Gamma(x)$ and the digamma function

Donal F. Connon

dconnon@btopenworld.com

19 May 2010

## Abstract

This paper considers various integrals where the integrand includes the log gamma function (or its derivative, the digamma function $\psi(x)$) multiplied by a trigonometric or hyperbolic function. Some apparently new integrals and series are evaluated; these include

$$\int_0^1 \log\Gamma(x)\cos p\pi x\,dx = \frac{\sin p\pi[\gamma + \log(2\pi)]}{2p\pi} + \frac{\sin p\pi}{4p\pi}\left[\psi\left(\frac{p}{2}\right) + \psi\left(-\frac{p}{2}\right)\right]$$

$$+ \frac{2(1-\cos p\pi)[\gamma + \log(2\pi)]}{\pi^2}\left[\frac{1}{2p^2} - \frac{\pi}{4p}\cot\left(\frac{p\pi}{2}\right)\right] + \frac{2(1-\cos p\pi)}{\pi^2}\sum_{n=1}^{\infty}\frac{\log n}{4n^2 - p^2}$$

$$\int_0^1 \log\Gamma(x)\sin p\pi x\,dx = \frac{(1-\cos p\pi)[\gamma + \log(2\pi)]}{2p\pi} + \frac{(1-\cos p\pi)}{4p\pi}\left[\psi\left(\frac{p}{2}\right) + \psi\left(-\frac{p}{2}\right)\right]$$

$$- \frac{2\sin p\pi[\gamma + \log(2\pi)]}{\pi^2}\left[\frac{1}{2p^2} - \frac{\pi}{4p}\cot\left(\frac{p\pi}{2}\right)\right] - \frac{2\sin p\pi}{\pi^2}\sum_{n=1}^{\infty}\frac{\log n}{4n^2 - p^2}$$

$$\int_0^1 \log\Gamma(x)\cos p\pi x\,dx = \left[\frac{1}{2}\log(2\pi) - 1\right]\frac{\sin p\pi}{p\pi} + \frac{Si(p\pi)}{p\pi} + \frac{2(1-\cos p\pi)}{\pi^2}\sum_{n=1}^{\infty}\frac{Ci(2n\pi)}{4n^2 - p^2}$$

$$+ \frac{p\sin p\pi}{\pi^2}\sum_{n=1}^{\infty}\frac{1}{n}\frac{si(2n\pi)}{4n^2 - p^2}$$

$$\int_0^1 \log\Gamma(x+1)\cot\pi x\,dx = \sum_{n=1}^{\infty}\frac{Ci(2n\pi)}{n\pi}$$

where $si(x)$ and $Ci(x)$ are the sine and cosine integrals.

CONTENTS                                                         Page

**1. A connection with the sine and cosine integrals**

The frequently thumbed table of integrals compiled by Gradshteyn and Ryzhik [33] contains only a handful of definite integrals involving $\log \Gamma(x)$. Of these, one of the apparently more complex examples is the following formula recorded in [33, p.650, 6.443.5] which is stated to be valid for $a > 0$

$$\int\limits_0^1 \log \Gamma(x+a) \sin 2k\pi x \, dx = -\frac{1}{2k\pi} \big[ \log a + \cos(2k\pi a) \, Ci(2k\pi a) - \sin(2k\pi a) \, si(2k\pi a) \big]$$

but, as shown below, the corrected integral which is in fact valid for $a \geq 0$ is

**Proposition 1.1**

$$(1.1) \quad \int\limits_0^1 \log \Gamma(x+a) \sin 2k\pi x \, dx = -\frac{1}{2k\pi} \big[ \log a - \cos(2k\pi a) \, Ci(2k\pi a) - \sin(2k\pi a) \, si(2k\pi a) \big]$$

where $si(x)$ and $Ci(x)$ are the sine and cosine integrals defined [33, p.878] by

$$(1.2) \qquad si(x) = -\int_x^\infty \frac{\sin t}{t} dt$$

and for $x > 0$

$$(1.3) \qquad Ci(x) = -\int_x^\infty \frac{\cos t}{t} dt = \gamma + \log x + \int_0^x \frac{\cos t - 1}{t} dt$$

where $\gamma$ is Euler's constant.

Note that $Si(x)$ is a slightly different sine integral which is defined in [33, p.878] and also in [1, p.231] by

$$(1.4) \qquad Si(x) = \int_0^x \frac{\sin t}{t} dt$$

We have

$$si(x) = -\int_x^\infty \frac{\sin t}{t} dt = \int_0^x \frac{\sin t}{t} dt - \int_0^\infty \frac{\sin t}{t} dt$$

and using the well-known integral from Fourier series analysis

$$(1.5) \qquad \frac{\pi}{2} = \int_0^\infty \frac{\sin t}{t} dt$$

we therefore see that the two sine integrals are intimately related by

$$(1.6) \quad si(x) = Si(x) - \frac{\pi}{2}$$

**Proof:**

We start with equation (A.8) which is derived in Appendix A to this paper

$$(1.7) \qquad \log \Gamma(x+a) = -\log(x+a) - \gamma(x+a) + \sum_{n=1}^\infty \left[ \log n - \log\left(n+a+x\right) + \frac{x+a}{n} \right]$$

and multiply this by $\sin p\pi x$ and integrate to obtain

$$\int\limits_0^1 \log\Gamma(x+a)\sin p\pi x\,dx = -\int\limits_0^1 \log(x+a)\sin p\pi x\,dx - \gamma\int\limits_0^1 (x+a)\sin p\pi x\,dx$$

$$+ \sum_{n=1}^{\infty}\int\limits_0^1 \left[\log n - \log(n+a+x) + \frac{x+a}{n}\right]\sin p\pi x\,dx$$

where we have assumed that the interchange of the order of integration and summation is valid.

The three component integrals are dealt with in turn below.

We have using integration by parts

$$\int \log(a+x)\sin p\pi x\,dx = -\frac{\log(a+x)\cos p\pi x}{p\pi} + \frac{1}{p\pi}\int\frac{\cos p\pi x}{a+x}\,dx$$

and with the substitution $t = a+x$ we get

$$\frac{1}{p\pi}\int\frac{\cos p\pi x}{a+x}\,dx = \int\frac{\cos p\pi a\cos p\pi t + \sin p\pi a\sin p\pi t}{p\pi t}\,dt$$

$$= \cos p\pi a\int\frac{\cos p\pi t}{p\pi t}\,dt + \sin p\pi a\int\frac{\sin p\pi t}{p\pi t}\,dt$$

$$= \frac{\cos p\pi a}{p\pi}\int\frac{\cos u}{u}\,du + \frac{\sin p\pi a}{p\pi}\int\frac{\sin u}{u}\,du$$

and reference to the definitions of the sine and cosine integrals shows that this is equivalent to

$$= \frac{\cos p\pi a\,Ci(u) + \sin p\pi a\,Si(u)}{p\pi}$$

Hence we obtain

$$\int\log(a+x)\sin p\pi x\,dx = \frac{\cos p\pi a\,Ci[\,p\pi(a+x)] + \sin p\pi a\,Si[\,p\pi(a+x)] - \log(a+x)\cos p\pi x}{p\pi}$$

which may be easily verified by differentiation since

$$\frac{d}{dx}Si(x) = \frac{\sin x}{x} \quad\text{and}\quad \frac{d}{dx}Ci(x) = \frac{\cos x}{x}$$

The definite integral becomes

(1.8) $\int_0^1 \log(a+x)\sin p\pi x\, dx$

$$= \frac{\cos p\pi a\{Ci[p\pi(a+1)]-Ci[p\pi a]\}+\sin p\pi a\{Si[p\pi(a+1)]-Si[p\pi a]\}}{p\pi}$$

$$-\frac{\log(a+1)\cos p\pi - \log a}{p\pi}$$

We also have the definite integral (where $k$ is an integer)

(1.9) $\int_0^1 \log(a+x)\sin 2k\pi x\, dx$

$$= \frac{\cos 2k\pi a\{Ci[2k\pi(a+1)]-Ci[2k\pi a]\}+\sin 2k\pi a\{si[2k\pi(a+1)]-si[2k\pi a]\}}{2k\pi}$$

$$-\frac{\log(a+1)-\log a}{2k\pi}$$

where we note from (1.6) that $Si(u)-Si(v)=si(u)-si(v)$.

We also have the elementary integral

$$\int (x+a)\sin p\pi x\, dx = \frac{\sin p\pi x}{p^2\pi^2}-\frac{(x+a)\cos p\pi x}{p\pi}$$

and therefore

$$\int_0^1 (x+a)\sin p\pi x\, dx = \frac{\sin p\pi}{p^2\pi^2}-\frac{(1+a)\cos p\pi - a}{p\pi}$$

With $p=2k$ we get

$$\int_0^1 (x+a)\sin 2k\pi x\, dx = -\frac{1}{2k\pi}$$

We see by letting $a \to n+a$ in (1.8) that

$$\int_0^1 \log(n+a+x)\sin 2k\pi x\,dx = \frac{\cos 2k\pi(a+n)\{Ci[2k\pi(a+n+1)] - Ci[2k\pi(a+n)]\}}{2k\pi}$$

$$+ \frac{\sin 2k\pi(a+n)\{si[2k\pi(a+n+1)] - si[2k\pi(a+n)]\}}{2k\pi}$$

$$- \frac{\log(a+n+1) - \log(a+n)}{2k\pi}$$

$$= \frac{\cos 2k\pi a\{Ci[2k\pi(a+n+1)] - Ci[2k\pi(a+n)]\}}{2k\pi}$$

$$- \frac{\sin 2k\pi a\{si[2k\pi(a+n+1)] - si[2k\pi(a+n)]\}}{2k\pi}$$

$$- \frac{\log(a+n+1) - \log(a+n)}{2k\pi}$$

Due to telescoping we see that

$$\frac{\cos 2k\pi a}{2k\pi}\sum_{n=1}^{\infty}\{Ci[2k\pi(a+n+1)] - Ci[2k\pi(a+n)]\} = -\frac{\cos 2k\pi a\,Ci[2k\pi(a+1)]}{2k\pi}$$

$$\frac{\sin 2k\pi a}{2k\pi}\sum_{n=1}^{\infty}\{si[2k\pi(a+n+1)] - si[2k\pi(a+n)]\} = -\frac{\sin 2k\pi a\,si[2k\pi(a+1)]}{2k\pi}$$

We also see that

$$\frac{1}{2k\pi}\sum_{n=1}^{N}\left[\log(a+n+1) - \log(a+n) - \frac{1}{n}\right] = \frac{1}{2k\pi}\left[\log(N+1+a) - H_N - \log(1+a)\right]$$

and as $N \to \infty$ we obtain

$$\frac{1}{2k\pi}\sum_{n=1}^{\infty}\left[\log(a+n+1) - \log(a+n) - \frac{1}{n}\right] = -\frac{\gamma + \log(1+a)}{2k\pi}$$

Finally, we note that the integral involving $\int_0^1 \log n \sin 2k\pi x\,dx$ vanishes and we have thereby obtained

$$\int\limits_0^1 \log \Gamma(x+a)\sin 2k\pi x\,dx$$

$$= -\frac{\cos 2k\pi a\left\{Ci[2k\pi(a+1)] - Ci[2k\pi a]\right\} + \sin 2k\pi a\left\{si[2k\pi(a+1)] - si[2k\pi a]\right\}}{2k\pi}$$

$$+ \frac{\log(a+1) - \log a}{2k\pi} + \frac{\cos 2k\pi a\,Ci[2k\pi(a+1)]}{2k\pi} - \frac{\sin 2k\pi a\,si[2k\pi(a+1)]}{2k\pi}$$

$$+ \frac{\gamma}{2k\pi} - \frac{\gamma + \log(a+1)}{2k\pi}$$

which simplifies to

$$\int\limits_0^1 \log \Gamma(x+a)\sin 2k\pi x\,dx = -\frac{1}{2k\pi}\left[\log a - \cos(2k\pi a)\,Ci(2k\pi a) - \sin(2k\pi a)\,si(2k\pi a)\right]$$

This therefore shows that one of the signs in (1.1) is recorded incorrectly in Gradshteyn and Ryzhik [33, p.650, 6.443.5] (and also in "Integrals and Series", Volume 2, p.60 by Prudnikov et al. [48]) and, in this regard, I note that both Havil [37, p.126] and Elizalde [29] define $Ci(x)$ as the negative of (1.3); it seems that this lack of consistency in the definition of $Ci(x)$ is likely to be the source of the error.

The definitions used in this paper correspond with those employed by Nielsen [44] and Nörlund [46] (except that those authors use the notation $ci(x)$ for $Ci(x)$). The similarity between the integral definitions in (1.2) and (1.3) does indeed support the notation $ci(x)$ instead of $Ci(x)$ which, unfortunately, has been used in this paper and elsewhere.

Equation (1.1) also applies in the limit as $a \to 0$ because

$$\lim_{y\to 0}[\cos y\,Ci(y) - \log y] = \lim_{y\to 0}\left[\gamma \cos y + \log y[\cos y - 1] + \cos y \int\limits_0^y \frac{\cos t - 1}{t}\,dt\right]$$

so that

(1.9.1)    $$\lim_{y\to 0}[\cos y\,Ci(y) - \log y] = \gamma$$

since, applying L'Hôpital's rule, we see that

$$\lim_{y\to 0}[\log y(\cos y - 1)] = \lim_{y\to 0}\left[y\log y\,\frac{\cos y - 1}{y}\right]$$

$$= \lim_{y \to 0} [y \log y] \lim_{y \to 0} \left[ \frac{\cos y - 1}{y} \right]$$

$$= -\lim_{y \to 0} [y \log y] \lim_{y \to 0} [\sin y] = 0$$

and we therefore obtain the well known Fourier coefficients [33, p.650, 6.443.3]

$$(1.9.2) \qquad \int_0^1 \log \Gamma(x) \sin 2k\pi x \, dx = \frac{\gamma + \log 2k\pi}{2k\pi}$$

With $a = 1/2$ we have

$$(1.10) \qquad \int_0^1 \log \Gamma \left( x + \frac{1}{2} \right) \sin 2k\pi x \, dx = \frac{1}{2k\pi} \left[ \log 2 + (-1)^k Ci(k\pi) \right]$$

$\square$

Using the definition of $Ci(x)$ in (1.3) we see that

$$Ci(ax) - \cos(ax) \log x = \gamma + \log ax - \cos(ax) \log x + \int_0^{ax} \frac{\cos t - 1}{t} dt$$

$$= \gamma + \log a - \log x [\cos(ax) - 1] + \int_0^{ax} \frac{\cos t - 1}{t} dt$$

We consider the limit

$$\lim_{x \to 0} [\cos(ax) - 1] \log x = \lim_{x \to 0} x \log x \frac{\cos(ax) - 1}{x}$$

$$= \lim_{x \to 0} x \log x \lim_{x \to 0} \frac{\cos(ax) - 1}{x}$$

and using L'Hôpital's rule we obtain

$$\lim_{x \to 0} \frac{\cos(ax) - 1}{x} = -\lim_{x \to 0} \frac{a \sin(ax)}{1}$$

which shows that

$$(1.10.1) \qquad \lim_{x \to 0} [\cos(ax) - 1] \log x = 0$$

and then taking the limit as $x \to 0$ we obtain

(1.10.2)     $\displaystyle \lim_{x \to 0} \left[ Ci(ax) - \cos(ax) \log x \right] = \gamma + \log a$

Hence by letting $a \to 0$ in (1.8) we obtain

(1.10.3)     $\displaystyle \int_0^1 \log x \sin p\pi x \, dx = \frac{Ci(p\pi) - \gamma - \log p\pi}{p\pi}$

(1.10.4)     $\displaystyle \int_0^1 \log x \sin 2k\pi x \, dx = \frac{Ci(2k\pi) - \gamma - \log 2k\pi}{2k\pi}$

With $a = 1$ in (1.1) we have

(1.11)     $\displaystyle \int_0^1 \log \Gamma(x+1) \sin 2k\pi x \, dx = \frac{Ci(2k\pi)}{2k\pi}$

which is (indirectly) reported with the correct sign in [33, p.650].

An alternative proof of (1.11) is shown below. We have

$$\int_0^1 \log \Gamma(x+1) \sin 2k\pi x \, dx = \int_0^1 \log x \sin 2k\pi x \, dx + \int_0^1 \log \Gamma(x) \sin 2k\pi x \, dx$$

and equation (1.11) results by combining (1.10) and (1.10.4).

$\square$

We now consider the integral $\displaystyle \int_0^1 \log \Gamma(x+a) \sin(2k+1)\pi x \, dx$.

**Proposition 1.2**

$$\int_0^1 \log \Gamma(x) \sin(2k+1)\pi x \, dx = \frac{1}{(2k+1)\pi} \left[ \log\left(\frac{\pi}{2}\right) + \frac{1}{2k+1} + 2 \sum_{j=1}^{k-1} \frac{1}{2j+1} \right]$$

**Proof**

We multiply (1.7) by $\sin(2k+1)\pi x$ and integrate to obtain

$$\int\limits_0^1 \log \Gamma(x+a)\sin(2k+1)\pi x\, dx = -\int\limits_0^1 \log(x+a)\sin(2k+1)\pi x\, dx - \gamma\int\limits_0^1 (x+a)\sin(2k+1)\pi x\, dx$$

$$+\sum_{n=1}^{\infty}\int\limits_0^1 \left[\log n - \log(n+a+x) + \frac{x+a}{n}\right]\sin(2k+1)\pi x\, dx$$

These three integrals are dealt with in turn below.

Using (1.8) we get

$$(1.12) \quad \int\limits_0^1 \log(a+x)\sin(2k+1)\pi x\, dx$$

$$= \frac{\cos(2k+1)\pi a\left\{Ci[(2k+1)\pi(a+1)] - Ci[(2k+1)\pi a]\right\}}{(2k+1)\pi}$$

$$+ \frac{\sin(2k+1)\pi a\left\{si[(2k+1)\pi(a+1)] - si[(2k+1)\pi a]\right\}}{(2k+1)\pi}$$

$$+ \frac{\log(a+1) + \log a}{(2k+1)\pi}$$

and

$$\int\limits_0^1 \log(n+a+x)\sin(2k+1)\pi x\, dx$$

$$= \frac{\cos[(2k+1)(n+a)\pi]\left\{Ci[(2k+1)\pi(n+a+1)] - Ci[(2k+1)\pi(n+a)]\right\}}{(2k+1)\pi}$$

$$+ \frac{\sin[(2k+1)(n+a)\pi]\left\{si[(2k+1)\pi(n+a+1)] - si[(2k+1)(n+a)\pi]\right\}}{(2k+1)\pi}$$

$$+ \frac{\log(n+a+1) + \log(n+a)}{(2k+1)\pi}$$

Since

$$\cos[(2k+1)(n+a)\pi] = (-1)^n \cos(2k+1)a\pi$$

$$\sin[(2k+1)(n+a)\pi] = (-1)^n \sin(2k+1)a\pi$$

we obtain

$$\int\limits_0^1 \log(n+a+x)\sin(2k+1)\pi x\,dx$$

$$= \frac{\cos[(2k+1)a\pi](-1)^n\left\{Ci[(2k+1)\pi(n+a+1)]-Ci[(2k+1)\pi(n+a)]\right\}}{(2k+1)\pi}$$

$$+ \frac{\sin[(2k+1)a\pi](-1)^n\left\{si[(2k+1)\pi(n+a+1)]-si[(2k+1)(n+a)\pi]\right\}}{(2k+1)\pi}$$

$$+ \frac{\log(n+a+1)+\log(n+a)}{(2k+1)\pi}$$

We saw above that

$$\int\limits_0^1 (x+a)\sin p\pi x\,dx = \frac{\sin p\pi}{p^2\pi^2} - \frac{(1+a)\cos p\pi - a}{p\pi}$$

and hence we have

$$\int\limits_0^1 (x+a)\sin(2k+1)\pi x\,dx = \frac{2a+1}{(2k+1)\pi}$$

Since $\sum\limits_{n=1}^\infty (-1)^n[x_{n+1}-x_n] = -x_1 - 2\sum\limits_{n=1}^\infty (-1)^n x_n$ we have

(1.13.1) $\quad \sum\limits_{n=1}^\infty (-1)^n\left\{Ci[(2k+1)\pi(a+n+1)]-Ci[(2k+1)\pi(a+n)]\right\}$

$$= -Ci[(2k+1)\pi(a+1)] - 2\sum\limits_{n=1}^\infty (-1)^n Ci[(2k+1)\pi(a+n)]$$

and

(1.13.2) $\quad \sum\limits_{n=1}^\infty (-1)^n\left\{si[(2k+1)\pi(a+n+1)]-si[(2k+1)\pi(a+n)]\right\}$

$$= -si[(2k+1)\pi(a+1)] - 2\sum\limits_{n=1}^\infty (-1)^n si[(2k+1)\pi(a+n)]$$

We will see later in (1.70) that

$$(1.13.3) \qquad \psi(a) = \log a - \frac{1}{2a} + 2\sum_{n=1}^{\infty}[\cos(2n\pi a)Ci(2n\pi a) + \sin(2n\pi a)si(2n\pi a)]$$

and we have [52, p.14] (by differentiating the corresponding expression for the log gamma function)

$$\psi(a+k) = \psi(a) + \sum_{j=1}^{k}\frac{1}{a+j-1}$$

$$= \psi(a) + \sum_{j=0}^{k-1}\frac{1}{a+j}$$

$$= \psi(a) + H_k^{(1)}(a)$$

where $H_n^{(m)}(a)$ is the generalised harmonic number function defined by

$$H_k^{(m)}(a) = \sum_{j=0}^{k-1}\frac{1}{(a+j)^m}$$

With $a = 1/2$ we obtain [52, p.20]

$$(1.14) \quad \psi\left(\frac{2k+1}{2}\right) = -\gamma - 2\log 2 + 2\sum_{j=0}^{k-1}\frac{1}{2j+1}$$

and we then see from (1.13.3) that

$$\psi\left(\frac{2k+1}{2}\right) = \log\left(\frac{2k+1}{2}\right) - \frac{1}{2k+1} + 2\sum_{n=1}^{\infty}\cos[(2k+1)n\pi]Ci[(2k+1)n\pi]$$

or equivalently

$$(1.15) \quad 2\sum_{n=1}^{\infty}(-1)^n Ci[(2k+1)n\pi] = -\gamma - \log 2 - \log(2k+1) + \frac{1}{2k+1} + 2\sum_{j=0}^{k-1}\frac{1}{2j+1}$$

We next need to consider the series

$$\frac{1}{2k+1}\sum_{n=1}^{\infty}\left[2\log n - \log(n+a+1) - \log(n+a) + \frac{2a+1}{n}\right]$$

Using (1.7)

$$\log\Gamma(x+a)+\log(x+a)+\gamma(x+a)=\sum_{n=1}^{\infty}\left[\log n-\log\left(n+a+x\right)+\frac{x+a}{n}\right]$$

which, with $x\rightarrow 1+x$, becomes

$$\log\Gamma(1+x+a)+\log(1+x+a)+\gamma(1+x+a)=\sum_{n=1}^{\infty}\left[\log n-\log\left(n+a+1+x\right)+\frac{1+x+a}{n}\right]$$

Adding these two equations together and noting the functional equation

$$\log\Gamma(1+x+a)=\log(x+a)+\log\Gamma(x+a)$$

results in

$$2\log\Gamma(x+a)+2\log(x+a)+\log(1+x+a)+\gamma(x+a)+\gamma=$$

$$=\sum_{n=1}^{\infty}\left[2\log n-\log\left(n+a+1+x\right)-\log\left(n+a+x\right)+\frac{2a+1+2x}{n}\right]$$

and with $x=0$ we see that

(1.16)     $$2\log\Gamma(1+a)+\log(1+a)+\gamma(1+a)$$

$$=\sum_{n=1}^{\infty}\left[2\log n-\log(n+a+1)-\log(n+a)+\frac{2a+1}{n}\right]$$

With $a=0$ we obtain the familiar formula

$$\sum_{n=1}^{\infty}\left[\log n-\log(n+1)+\frac{1}{n}\right]=\gamma$$

We therefore obtain by combining the various integrals

$$\int_{0}^{1}\log\Gamma(x+a)\sin(2k+1)\pi x\,dx$$

$$=-\frac{\cos(2k+1)\pi a\left\{Ci[(2k+1)\pi(a+1)]-Ci[(2k+1)\pi a]\right\}}{(2k+1)\pi}$$

$$-\frac{\sin(2k+1)\pi a\left\{si[(2k+1)\pi(a+1)]-si[(2k+1)\pi a]\right\}}{(2k+1)\pi}$$

$$-\frac{\log(a+1)+\log a}{(2k+1)\pi}-\frac{(2a+1)\gamma}{(2k+1)\pi}$$

$$+\frac{\cos(2k+1)\pi a\, Ci[(2k+1)\pi(a+1)]}{(2k+1)\pi}+\frac{2\cos(2k+1)\pi a}{(2k+1)\pi}\sum_{n=1}^{\infty}(-1)^n Ci[(2k+1)\pi(a+n)]$$

$$+\frac{\sin(2k+1)\pi a\, si[(2k+1)\pi(a+1)]}{(2k+1)\pi}+\frac{2\sin(2k+1)\pi a}{(2k+1)\pi}\sum_{n=1}^{\infty}(-1)^n si[(2k+1)\pi(a+n)]$$

$$+\frac{2\log\Gamma(1+a)+\log(1+a)+\gamma(1+a)}{(2k+1)\pi}$$

which simplifies to

$$(1.17)\qquad\int_0^1 \log\Gamma(x+a)\sin(2k+1)\pi x\, dx$$

$$=\frac{\cos(2k+1)\pi a\, Ci[(2k+1)\pi a]+\sin(2k+1)\pi a\, si[(2k+1)\pi a]}{(2k+1)\pi}$$

$$+\frac{2\log\Gamma(1+a)-\gamma a-\log a}{(2k+1)\pi}$$

$$+\frac{2\cos(2k+1)\pi a}{(2k+1)\pi}\sum_{n=1}^{\infty}(-1)^n Ci[(2k+1)\pi(a+n)]$$

$$+\frac{2\sin(2k+1)\pi a}{(2k+1)\pi}\sum_{n=1}^{\infty}(-1)^n si[(2k+1)\pi(a+n)]$$

With $a=0$ we obtain

$$\int_0^1 \log\Gamma(x)\sin(2k+1)\pi x\, dx$$

$$=\lim_{a\to 0}\frac{Ci[(2k+1)\pi a]-\log a}{(2k+1)\pi}+\frac{2}{(2k+1)\pi}\sum_{n=1}^{\infty}(-1)^n Ci[(2k+1)\pi n]$$

and referring to (1.15) and (1.19) this becomes

$$= \frac{\gamma + \log[(2k+1)\pi]}{(2k+1)\pi} + \frac{1}{(2k+1)\pi}\left[-\gamma - \log 2 - \log(2k+1) + \frac{1}{2k+1} + 2\sum_{j=1}^{k-1}\frac{1}{2j+1}\right]$$

Hence we obtain the known result (see also Section 4)

(1.18) $\displaystyle\int_0^1 \log\Gamma(x)\sin(2k+1)\pi x\, dx = \frac{1}{(2k+1)\pi}\left[\log\left(\frac{\pi}{2}\right) + \frac{1}{2k+1} + 2\sum_{j=1}^{k-1}\frac{1}{2j+1}\right]$

Referring to (1.3) we see that

$$Ci(ux) = \gamma + \log u + \log x + \int_0^{ux}\frac{\cos t - 1}{t}dt$$

and we have

(1.19) $\displaystyle\lim_{x\to 0}\left[Ci(ux) - \log x\right] = \gamma + \log u$

We also see that

$$\cos(ux)Ci(ux) = [\gamma + \log u] + \cos(ux)\log x + \cos(ux)\int_0^{ux}\frac{\cos t - 1}{t}dt$$

$$= [\gamma + \log u] + \log x + [\cos(ux) - 1]\log x + \cos(ux)\int_0^{ux}\frac{\cos t - 1}{t}dt$$

and thus

$$\lim_{x\to 0}\left[\cos(ux)Ci(ux) - \log x\right] = [\gamma + \log u]\cos(ux) + \lim_{x\to 0}[\cos(ux) - 1]\log x$$

Then using (1.10.1) we obtain

(1.20) $\displaystyle\lim_{x\to 0}\left[\cos(ux)Ci(ux) - \log x\right] = \gamma + \log u$

and in particular we have (as used above)

(1.21) $\displaystyle\lim_{a\to 0}\frac{Ci[(2k+1)\pi a] - \log a}{(2k+1)\pi} = \frac{\gamma + \log[(2k+1)\pi]}{(2k+1)\pi}$

$\square$

With $a = 1$ in (1.17) we obtain

(1.22) $\displaystyle\int\limits_0^1 \log\Gamma(x+1)\sin(2k+1)\pi x\,dx$

$$= -\frac{1}{(2k+1)\pi}\left(Ci[(2k+1)\pi]+\gamma+2\sum_{n=1}^\infty(-1)^n Ci[(2k+1)\pi(n+1)]\right)$$

$$= -\frac{1}{(2k+1)\pi}\left(Ci[(2k+1)\pi]+\gamma-2\sum_{m=2}^\infty(-1)^m Ci[(2k+1)\pi m)]\right)$$

$$= -\frac{1}{(2k+1)\pi}\left(-Ci[(2k+1)\pi]+\gamma-2\sum_{n=1}^\infty(-1)^n Ci[(2k+1)\pi n)]\right)$$

Letting $p=2k+1$ in (1.10.3) gives us

$$\int\limits_0^1 \log x\sin(2k+1)\pi x\,dx = \frac{Ci[(2k+1)\pi]-\gamma-\log[(2k+1)\pi]}{(2k+1)\pi}$$

and hence we have

$$\int\limits_0^1 \log\Gamma(x+1)\sin(2k+1)\pi x\,dx = \int\limits_0^1 \log x\sin(2k+1)\pi x\,dx + \int\limits_0^1 \log\Gamma(x)\sin(2k+1)\pi x\,dx$$

We then have

$$-\frac{1}{(2k+1)\pi}\left(-Ci[(2k+1)\pi]+\gamma-2\sum_{n=1}^\infty(-1)^n Ci[(2k+1)\pi n)]\right)=$$

$$\frac{Ci[(2k+1)\pi]-\gamma-\log[(2k+1)\pi]}{(2k+1)\pi}+\int\limits_0^1 \log\Gamma(x)\sin(2k+1)\pi x\,dx$$

or

(1.23) $\displaystyle\int\limits_0^1 \log\Gamma(x)\sin(2k+1)\pi x\,dx = \frac{2}{(2k+1)\pi}\sum_{n=1}^\infty(-1)^n Ci[(2k+1)\pi n)]+\frac{\log[(2k+1)\pi]}{(2k+1)\pi}$

which corresponds with (1.18). Therefore we have

(1.23.1) $\displaystyle 2\sum_{n=1}^\infty(-1)^n Ci[(2k+1)\pi n)]+\log[(2k+1)\pi]=\log\left(\frac{\pi}{2}\right)+\frac{1}{2k+1}+2\sum_{j=1}^{k-1}\frac{1}{2j+1}$

and with $k=0$ this becomes

(1.23.2)  $\quad 2\sum_{n=1}^{\infty}(-1)^n Ci(\pi n) = 1 - \log 2$

There is however an unexplained difference between this result and (1.72)

$$2\sum_{n=1}^{\infty}(-1)^n Ci(n\pi) = 1 - \gamma - \log 2$$

$\square$

We now consider the integral $\int_0^1 \log\Gamma(x+a)\cos 2k\pi x\,dx$ .

**Proposition 1.3**

$$\int_0^1 \log\Gamma(x+a)\cos 2k\pi x\,dx = -\frac{1}{2k\pi}\bigl[-\sin(2k\pi a)Ci(2k\pi a) + \cos(2k\pi a)si(2k\pi a)\bigr]$$

**Proof**

We multiply (1.7) by $\cos p\pi x$ and integrate to obtain

$$\int_0^1 \log\Gamma(x+a)\cos p\pi x\,dx = -\int_0^1 \log(x+a)\cos p\pi x\,dx - \gamma\int_0^1 (x+a)\cos p\pi x\,dx$$

$$+\sum_{n=1}^{\infty}\int_0^1 \left[\log n - \log\bigl(n+a+x\bigr) + \frac{x+a}{n}\right]\cos p\pi x\,dx$$

These three integrals are dealt with in turn below.

As before, integration by parts gives us

$$\int \log(a+x)\cos p\pi x\,dx = \frac{\log(a+x)\sin p\pi x}{p\pi} - \frac{1}{p\pi}\int \frac{\sin p\pi x}{a+x}\,dx$$

and with the substitution $t = a + x$ we get

$$\frac{1}{p\pi}\int \frac{\sin p\pi x}{a+x}\,dx = \int \frac{\cos p\pi a\sin p\pi t - \sin p\pi a\cos p\pi t}{p\pi t}\,dt$$

$$= \cos p\pi a \int \frac{\sin p\pi t}{p\pi t}\,dt - \sin p\pi a \int \frac{\cos p\pi t}{p\pi t}\,dt$$

$$= \frac{\cos p\pi a}{p\pi} \int \frac{\sin u}{u}\,du - \frac{\sin u p\pi a}{p\pi} \int \frac{\cos u}{u}\,du$$

and reference to the definitions of the sine and cosine integrals shows that this is equivalent to

$$= \frac{\cos p\pi a\,Si(u) - \sin p\pi a\,Ci(u)}{p\pi}$$

Hence we obtain

$$\int \log(a+x)\cos p\pi x\,dx = \frac{\sin p\pi a\,Ci[\,p\pi(a+x)] - \cos p\pi a\,Si[\,p\pi(a+x)] + \log(a+x)\sin p\pi x}{p\pi}$$

The definite integral becomes

(1.24) $\int\limits_0^1 \log(a+x)\cos p\pi x\,dx$

$$= \frac{\sin p\pi a\{Ci[\,p\pi(a+1)] - Ci[\,p\pi a]\} - \cos p\pi a\{si[\,p\pi(a+1)] - si[\,p\pi a]\}}{p\pi}$$

$$+ \frac{\log(a+1)\sin p\pi}{p\pi}$$

We also have the definite integral (where $k$ is an integer)

(1.25) $\int\limits_0^1 \log(a+x)\cos 2k\pi x\,dx$

$$= \frac{\sin 2k\pi a\{Ci[2k\pi(a+1)] - Ci[2k\pi a]\} - \cos 2k\pi a\{si[2k\pi(a+1)] - si[2k\pi a]\}}{2k\pi}$$

With $a = 0$ we obtain

(1.26) $\int\limits_0^1 \log x\cos 2k\pi x\,dx = -\frac{Si[2k\pi]}{2k\pi}$

We also have

$$\int (x+a)\cos p\pi x\, dx = \frac{\cos p\pi x}{p^2\pi^2} + \frac{(x+a)\sin p\pi x}{p\pi}$$

and therefore

$$\int_0^1 (x+a)\cos p\pi x\, dx = \frac{\cos p\pi - 1}{p^2\pi^2} + \frac{(1+a)\sin p\pi}{p\pi}$$

With $p = 2k$ we get

$$\int_0^1 (x+a)\cos 2k\pi x\, dx = 0$$

We see by letting $a \to n+a$ in (1.25) that

$$\int_0^1 \log(n+a+x)\cos 2k\pi x\, dx = \frac{\sin 2k\pi(a+n)\{Ci[2k\pi(a+n+1)] - Ci[2k\pi(a+n)]\}}{2k\pi}$$

$$- \frac{\cos 2k\pi(a+n)\{si[2k\pi(a+n+1)] - si[2k\pi(a+n)]\}}{2k\pi}$$

$$= \frac{\sin 2k\pi a\{Ci[2k\pi(a+n+1)] - Ci[2k\pi(a+n)]\}}{2k\pi}$$

$$- \frac{\cos 2k\pi a\{si[2k\pi(a+n+1)] - si[2k\pi(a+n)]\}}{2k\pi}$$

As before, due to telescoping we see that

$$\frac{\sin 2k\pi a}{2k\pi}\sum_{n=1}^{\infty}\{Ci[2k\pi(a+n+1)] - Ci[2k\pi(a+n)]\} = -\frac{\sin 2k\pi a\, Ci[2k\pi(a+1)]}{2k\pi}$$

$$\frac{\cos 2k\pi a}{2k\pi}\sum_{n=1}^{\infty}\{si[2k\pi(a+n+1)] - si[2k\pi(a+n)]\} = -\frac{\cos 2k\pi a\, si[2k\pi(a+1)]}{2k\pi}$$

Finally, we note that the integral involving $\int_0^1 \log n\cos 2k\pi x\, dx$ vanishes and we have thereby obtained

$$\int_0^1 \log \Gamma(x+a) \cos 2k\pi x \, dx$$

$$= -\frac{\sin 2k\pi a \left\{ Ci[2k\pi(a+1)] - Ci[2k\pi a] \right\} - \cos 2k\pi a \left\{ si[2k\pi(a+1)] - si[2k\pi a] \right\}}{2k\pi}$$

$$+ \frac{\sin 2k\pi a \, Ci[2k\pi(a+1)]}{2k\pi} - \frac{\cos 2k\pi a \, si[2k\pi(a+1)]}{2k\pi}$$

which simplifies to

$$(1.27) \qquad \int_0^1 \log \Gamma(x+a) \cos 2k\pi x \, dx = -\frac{1}{2k\pi} \left[ -\sin(2k\pi a) Ci(2k\pi a) + \cos(2k\pi a) si(2k\pi a) \right]$$

This is the companion integral to (1.1) and this corrects the entry reported in [33, p.650, 6.443.3] for $a > 0$; the sign error referred to in (1.1) is also replicated here. This sign error also features in (1.76) and (1.78).

We may also note that (1.27) is also valid for $a = 0$ because

$$\lim_{u \to 0}[\sin u . \log u] = \lim_{u \to 0} \left[ \frac{\sin u}{u} u \log u \right] = 0$$

and therefore using (1.3) we see that $\lim_{a \to 0} \sin(2k\pi a) Ci(2k\pi a) = 0$.

Since from (1.16) $si(0) = -\frac{\pi}{2}$, we therefore obtain the well known Fourier coefficients [33, p.650, 6.443.3]

$$(1.28) \qquad \int_0^1 \log \Gamma(x) \cos 2k\pi x \, dx = \frac{1}{4k}$$

Further derivations of this integral are contained in (2.28) and (3.9) below.

Letting $a = 1/2$ results in

$$(1.28.1) \qquad \int_0^1 \log \Gamma\left(x + \frac{1}{2}\right) \cos 2k\pi x \, dx = \frac{(-1)^{k+1}}{2k\pi} si(k\pi)$$

With $a = 1$ in (1.27) we get for $k \geq 1$

$$(1.29) \quad \int_0^1 \log \Gamma(x+1) \cos 2k\pi x \, dx = -\frac{si(2k\pi)}{2k\pi}$$

This may be easily confirmed as follows. We have

$$\int_0^1 \log \Gamma(x+1) \cos 2k\pi x \, dx = \int_0^1 \log x \cos 2k\pi x \, dx + \int_0^1 \log \Gamma(x) \cos 2k\pi x \, dx$$

and (1.29) results by using (1.26) and (1.28).

Using the duplication formula for the gamma function

$$\log \Gamma(2x) = \log \Gamma(x) + \log \Gamma\left(x + \frac{1}{2}\right) + (2x-1)\log 2 - \frac{1}{2}\log \pi$$

we have

$$\int_0^1 \log \Gamma(2x) \cos 2k\pi x \, dx = \int_0^1 \log \Gamma(x) \cos 2k\pi x \, dx + \int_0^1 \log \Gamma\left(x + \frac{1}{2}\right) \cos 2k\pi x \, dx$$

$$+ \log 2 \int_0^1 (2x-1) \cos 2k\pi x \, dx$$

We note that

$$\int_0^1 (2x-1) \cos 2k\pi x \, dx = 0$$

$$\int_0^1 \log \Gamma(2x) \cos 2k\pi x \, dx = \frac{1}{4k} + \frac{(-1)^{k+1}}{2k\pi} si(k\pi)$$

$$\int_0^1 \log \Gamma(2x) \cos 2k\pi x \, dx = \frac{1}{2} \int_0^2 \log \Gamma(u) \cos k\pi u \, du$$

$$\int_0^2 \log \Gamma(u) \cos k\pi u \, du = \int_0^1 \log \Gamma(u) \cos k\pi u \, du + \int_1^2 \log \Gamma(u) \cos k\pi u \, du$$

$$\int_1^2 \log \Gamma(u) \cos k\pi u \, du = \int_0^1 \log \Gamma(1+x) \cos k\pi(1+x) \, dx$$

$$= (-1)^k \int_0^1 \log \Gamma(1+x) \cos k\pi x \, dx$$

$$\int_0^1 \log \Gamma(2x) \cos 2k\pi x \, dx = \frac{1}{2} \int_0^1 \log \Gamma(x) \cos k\pi x \, dx + (-1)^k \frac{1}{2} \int_0^1 \log \Gamma(1+x) \cos k\pi x \, dx$$

$$= \frac{1}{2}[1+(-1)^k] \int_0^1 \log \Gamma(x) \cos k\pi x \, dx + (-1)^k \frac{1}{2} \int_0^1 \log x \cos k\pi x \, dx$$

We note from Appendix B that

$$\int_0^1 \log x \cos k\pi x \, dx = -\frac{Si(k\pi)}{k\pi}$$

and unfortunately we simply end up with

$$\frac{1}{2}[1+(-1)^k] \int_0^1 \log \Gamma(x) \cos k\pi x \, dx = \frac{1}{4k} + \frac{(-1)^k}{4k}$$

$\square$

We recall (1.1) and (1.27)

$$\int_0^1 \log \Gamma(x+a) \sin 2n\pi x \, dx = -\frac{1}{2\pi n} \Big[ \log a - \cos(2n\pi a) Ci(2n\pi a) + \sin(2n\pi a) si(2n\pi a) \Big]$$

$$\int_0^1 \log \Gamma(x+a) \cos 2n\pi x \, dx = -\frac{1}{2n\pi} \Big[ -\sin(2n\pi a) Ci(2n\pi a) + \cos(2n\pi a) si(2n\pi a) \Big]$$

and multiplying each by $\sin(2n\pi a)$ and $\cos(2n\pi a)$ respectively we get

$$\sin(2n\pi a) \int_0^1 \log \Gamma(x+a) \sin 2n\pi x \, dx =$$

$$-\frac{1}{2\pi n} \Big[ \log a \sin(2n\pi a) - \cos(2n\pi a) \sin(2n\pi a) Ci(2n\pi a) - \sin^2(2n\pi a) si(2n\pi a) \Big]$$

and

$$\cos(2n\pi a) \int_0^1 \log \Gamma(x+a) \cos 2n\pi x \, dx =$$

$$-\frac{1}{2n}\Big[-\cos(2n\pi a)\sin(2n\pi a)Ci(2n\pi a)+\cos^2(2n\pi a)si(2n\pi a)\Big]$$

Adding these two equations results in

$$\int_0^1 \log\Gamma(x+a)[\sin(2n\pi a)\sin 2n\pi x+\cos(2n\pi a)\cos 2n\pi x]dx$$

$$=-\frac{1}{2\pi n}\Big[\log a\sin(2n\pi a)-\sin(4n\pi a)Ci(2n\pi a)+\cos(4n\pi a)si(2n\pi a)\Big]$$

or equivalently

$$\int_0^1 \log\Gamma(x+a)\cos[2n\pi(x-a)]dx$$

$$=-\frac{1}{2\pi n}\Big[\log a\sin(2n\pi a)-\sin(4n\pi a)Ci(2n\pi a)+\cos(4n\pi a)si(2n\pi a)\Big]$$

With $a=0$ we immediately obtain (1.28)

$$\int_0^1 \log\Gamma(x)\cos 2n\pi x\, dx=\frac{1}{4n}$$

and $a=1/2$ gives us (1.29)

$$\int_0^1 \log\Gamma(x+1/2)\cos 2n\pi x\, dx=(-1)^{n+1}\frac{si(n\pi)}{2\pi n}$$

Similarly subtraction gives us

$$\sin(2n\pi a)\int_0^1 \log\Gamma(x+a)\sin 2n\pi x\, dx-\cos(2n\pi a)\int_0^1 \log\Gamma(x+a)\cos 2n\pi x\, dx$$

$$=-\frac{1}{2\pi n}\Big[\log a\sin(2n\pi a)-\cos(2n\pi a)\sin(2n\pi a)Ci(2n\pi a)-\sin^2(2n\pi a)si(2n\pi a)\Big]$$

$$+\frac{1}{2\pi n}\Big[-\cos(2n\pi a)\sin(2n\pi a)Ci(2n\pi a)+\cos^2(2n\pi a)si(2n\pi a)\Big]$$

$$= \frac{1}{2\pi n}\big[-\log a \sin(2n\pi a) + si(2n\pi a)\big]$$

and we have

$$(1.29.1) \quad \int_0^1 \log \Gamma(x+a)\cos[2n\pi(x+a)]dx = \frac{1}{2\pi n}\big[\log a \sin(2n\pi a) - si(2n\pi a)\big]$$

The substitution $u = x + a$ gives us

$$(1.29.2) \quad \int_a^{a+1} \log \Gamma(u)\cos 2n\pi u\, du = \frac{1}{2\pi n}\big[\log a \sin(2n\pi a) - si(2n\pi a)\big]$$

and differentiation with respect to $a$ results in the very obvious identity

$$(1.29.3) \quad \log \Gamma(a+1)\cos 2n\pi(a+1) - \log \Gamma(a)\cos 2n\pi a$$

$$= \frac{1}{2\pi n}\left[2n\pi \log a \cos(2n\pi a) + \frac{1}{a}\sin(2n\pi a) - \frac{1}{a}\sin(2n\pi a)\right]$$

$$= \log a \cos(2n\pi a)$$

Finally, we consider the integral $\int_0^1 \log \Gamma(x+a)\cos(2k+1)\pi x\, dx$.

**Proposition 1.4**

$$\int_0^1 \log \Gamma(x)\cos(2k+1)\pi x\, dx = \frac{2}{\pi^2}\left[\frac{\log(2\pi)+\gamma}{(2k+1)^2} + 2\sum_{n=1}^{\infty}\frac{\log n}{4n^2-(2k+1)^2}\right]$$

**An attempted proof**

As before we have

$$\int_0^1 \log \Gamma(x+a)\cos(2k+1)\pi x\, dx = -\int_0^1 \log(x+a)\cos(2k+1)\pi x\, dx - \gamma\int_0^1 (x+a)\cos(2k+1)\pi x\, dx$$

$$+ \sum_{n=1}^{\infty}\int_0^1 \left[\log n - \log\big(n+a+x\big) + \frac{x+a}{n}\right]\cos(2k+1)\pi x\, dx$$

Using the definite integral (1.24) we see that

(1.30) $\displaystyle\int_0^1 \log(a+x)\cos(2k+1)\pi x\,dx$

$$= \frac{\sin(2k+1)\pi a\left\{Ci[(2k+1)\pi(a+1)]-Ci[(2k+1)\pi a]\right\}}{(2k+1)\pi}$$

$$- \frac{\cos(2k+1)\pi a\left\{si[(2k+1)\pi(a+1)]-si[(2k+1)\pi a]\right\}}{(2k+1)\pi}$$

$$\int_0^1 (x+a)\cos(2k+1)\pi x\,dx = -\frac{2}{(2k+1)^2\pi^2}$$

$$\int_0^1 \log(n+a+x)\cos(2k+1)\pi x\,dx$$

$$= \frac{\sin[(2k+1)(n+a)\pi]\left\{Ci[(2k+1)\pi(a+n+1)]-Ci[(2k+1)\pi(a+n)]\right\}}{(2k+1)\pi}$$

$$- \frac{\cos[(2k+1)(n+a)\pi]\left\{si[(2k+1)\pi(a+n+1)]-si[(2k+1)\pi(a+n)]\right\}}{(2k+1)\pi}$$

Since

$$\sin[(2k+1)(n+a)\pi] = (-1)^n\sin(2k+1)a\pi$$

$$\cos[(2k+1)(n+a)\pi] = (-1)^n\cos(2k+1)a\pi$$

this becomes

$$= \frac{(-1)^n\sin(2k+1)a\pi\left\{Ci[(2k+1)\pi(a+n+1)]-Ci[(2k+1)\pi(a+n)]\right\}}{(2k+1)\pi}$$

$$- \frac{(-1)^n\cos(2k+1)a\pi\left\{si[(2k+1)\pi(a+n+1)]-si[(2k+1)\pi(a+n)]\right\}}{(2k+1)\pi}$$

We have as in (1.13.1)

$$\sum_{n=1}^\infty (-1)^n\left\{Ci[(2k+1)\pi(a+n+1)]-Ci[(2k+1)\pi(a+n)]\right\}$$

$$= -Ci[(2k+1)\pi(a+1)] - 2\sum_{n=1}^{\infty}(-1)^n Ci[(2k+1)\pi(a+n)]$$

and

$$\sum_{n=1}^{\infty}(-1)^n\left\{si[(2k+1)\pi(a+n+1)] - si[(2k+1)\pi(a+n)]\right\}$$

$$= -si[(2k+1)\pi(a+1)] - 2\sum_{n=1}^{\infty}(-1)^n si[(2k+1)\pi(a+n)]$$

Finally, we note that the integral involving $\int_0^1 \log n \cos(2k+1)\pi x\,dx$ vanishes and we have thereby obtained

$$\int_0^1 \log\Gamma(x+a)\cos(2k+1)\pi x\,dx$$

$$= -\frac{\sin(2k+1)\pi a\left\{Ci[(2k+1)\pi(a+1)] - Ci[(2k+1)\pi a]\right\}}{(2k+1)\pi}$$

$$+ \frac{\cos(2k+1)\pi a\left\{si[(2k+1)\pi(a+1)] - si[(2k+1)\pi a]\right\}}{(2k+1)\pi}$$

$$+ \frac{2\gamma}{(2k+1)^2\pi^2}$$

$$+ \frac{\sin(2k+1)a\pi\,Ci[(2k+1)\pi(a+1)]}{(2k+1)\pi} + \frac{2\sin(2k+1)a\pi}{(2k+1)\pi}\sum_{n=1}^{\infty}(-1)^n Ci[(2k+1)\pi(a+n)]$$

$$- \frac{\cos(2k+1)a\pi\,si[(2k+1)\pi(a+1)]}{(2k+1)\pi} - \frac{2\cos(2k+1)a\pi}{(2k+1)\pi}\sum_{n=1}^{\infty}(-1)^n si[(2k+1)\pi(a+n)]$$

$$- \sum_{n=1}^{\infty}\frac{2}{(2k+1)^2\pi^2 n}$$

which simplifies to

$$= \frac{\sin(2k+1)\pi a\,Ci[(2k+1)\pi a] - \cos(2k+1)\pi a\,si[(2k+1)\pi a]}{(2k+1)\pi}$$

$$+ \frac{2\gamma}{(2k+1)^2\pi^2}$$

$$+\frac{2\sin(2k+1)a\pi}{(2k+1)\pi}\sum_{n=1}^{\infty}(-1)^n Ci[(2k+1)\pi(a+n)]$$

$$-\frac{2\cos(2k+1)a\pi}{(2k+1)\pi}\sum_{n=1}^{\infty}(-1)^n si[(2k+1)\pi(a+n)]$$

$$-\sum_{n=1}^{\infty}\frac{2}{(2k+1)^2\pi^2 n}$$

With $a=0$ we have

$$\int_0^1 \log\Gamma(x)\cos(2k+1)\pi x\,dx$$

$$=\frac{1}{2k+1}+\frac{2\gamma}{(2k+1)^2\pi^2}-\frac{2}{(2k+1)\pi}\sum_{n=1}^{\infty}(-1)^n si[(2k+1)\pi n]+\sum_{n=1}^{\infty}\frac{2}{(2k+1)^2\pi^2 n}$$

This is obviously incorrect because of the appearance of the divergent series; the author would appreciate receiving a corrected version of the proof.

Using a different approach we note from (4.9) that

$$\int_0^1 \log\Gamma(x)\cos(2k+1)\pi x\,dx=\frac{2}{\pi^2}\left[\frac{\log(2\pi)+\gamma}{(2k+1)^2}+2\sum_{n=1}^{\infty}\frac{\log n}{4n^2-(2k+1)^2}\right]$$

**Representations of $\log\Gamma(x)$ in terms of the sine and cosine integrals**

The following analysis is extracted from an earlier paper [21]; it indicates how the $\log\Gamma(x)$ function is itself intimately connected with the sine and cosine integrals.

Whittaker & Watson [56, p.261] posed the following question: Prove that for all values of $a$ except negative real values we have

$$(1.40) \qquad \log\Gamma(a)=\frac{1}{2}\log(2\pi)+\left(a-\frac{1}{2}\right)\log a-a+\frac{1}{\pi}\sum_{n=1}^{\infty}\int_0^{\infty}\frac{\sin(2n\pi x)}{n(x+a)}dx$$

and this result was attributed by Stieltjes to Bourguet. Equation (1.40) may also be derived using the Euler-Maclaurin summation formula (see for example Knopp's book [39, p.530]).

By differentiation we can easily see that

$$\frac{d}{dx}\big(\cos(2n\pi a)Si[2n\pi(x+a)]-\sin(2n\pi a)Ci[2n\pi(x+a)]\big)=\frac{\sin(2n\pi x)}{x+a}$$

and we therefore have

$$\int\limits_0^M \frac{\sin(2n\pi x)}{x+a}dx=\big(\cos(2n\pi a)Si[2n\pi(x+a)]-\sin(2n\pi a)Ci[2n\pi(x+a)]\big)\Big|_0^M$$

$$=\cos(2n\pi a)\big\{Si[2n\pi(M+a)]-Si[2n\pi a]\big\}-\sin(2n\pi a)\big\{Ci[2n\pi(M+a)]-Ci[2n\pi a]\big\}$$

From (1.4) and (1.5) we see that

$$\lim_{M\to\infty}Si[2n\pi(M+a)]=\frac{\pi}{2}$$

and from (1.3) we see that

$$\lim_{M\to\infty}Ci[2n\pi(M+a)]=0$$

Hence we obtain as $M\to\infty$

$$\int\limits_0^\infty \frac{\sin(2n\pi x)}{x+a}dx=\cos(2n\pi a)\left\{\frac{\pi}{2}-Si(2n\pi a)\right\}+\sin(2n\pi a)Ci(2n\pi a)$$

and reference to (1.6) shows that this is equal to

$$=-\cos(2n\pi a)si(2n\pi a)+\sin(2n\pi a)Ci(2n\pi a)$$

Therefore using Bourguet's formula we have

(1.41)    $\log\Gamma(a)=$

$$\frac{1}{2}\log(2\pi)+\left(a-\frac{1}{2}\right)\log a-a+\frac{1}{\pi}\sum_{n=1}^\infty\frac{1}{n}[\sin(2n\pi a)Ci(2n\pi a)-\cos(2n\pi a)si(2n\pi a)]$$

Nielsen [44, p.79] also reports a similar formula which is valid for $0<a<1$

(1.41.1)   $\log\Gamma(a)=\dfrac{1}{2}\log(2\pi)-1-\log a+\dfrac{1}{\pi}\sum_{n=1}^\infty\dfrac{1}{n}[\sin 2n\pi a\,Ci(2n\pi)-\cos 2n\pi a\,si(2n\pi)]$

and a derivation of this is contained in Appendix B. Note that in this case, the sine and cosine integral functions do not contain the variable $a$ in their arguments.

Equation (1.41) may be expressed as

$$\log \Gamma(a) = \frac{1}{2}\log(2\pi) + \left(a - \frac{1}{2}\right)\log a - a + \sum_{n=1}^{\infty} C_n$$

and we have with $a \to 2a$

$$\log \Gamma(2a) = \frac{1}{2}\log(2\pi) + \left(2a - \frac{1}{2}\right)\log(2a) - 2a + 2\sum_{n=1}^{\infty} C_{2n}$$

Since $2\sum_{n=1}^{\infty} C_{2n} = \sum_{n=1}^{\infty} C_n + \sum_{n=1}^{\infty} (-1)^n C_n$ we have

$$\log \Gamma(2a) = \left(2a - \frac{1}{2}\right)\log(2a) - a + \log \Gamma(a) - \left(a - \frac{1}{2}\right)\log a + \sum_{n=1}^{\infty} (-1)^n C_n$$

and using Legendre's duplication formula for the gamma function [52, p.7]

$$\sqrt{\pi}\,\Gamma(2a) = 2^{2a-1}\Gamma(a)\Gamma\left(a + \frac{1}{2}\right)$$

this results in

(1.42)     $\log \Gamma(a + 1/2) =$

$$\frac{1}{2}\log(2\pi) + a \log a - a + \frac{1}{\pi}\sum_{n=1}^{\infty} \frac{(-1)^n}{n}[\sin(2n\pi a)Ci(2n\pi a) - \cos(2n\pi a)si(2n\pi a)]$$

which is also reported by Nörlund [46, p.114]. Since $si(0) = -\pi/2$ this identity may be easily verified for $a = 0$. With $a = 1/2$ we obtain

(1.43)     $$\sum_{n=1}^{\infty} \frac{si(n\pi)}{n} = \frac{\pi}{2}\log \pi - \frac{\pi}{2}$$

which is contained in [44, p,82].

Letting $a = 1$ we have

(1.44)     $$\sum_{n=1}^{\infty} (-1)^n \frac{si(2n\pi)}{n} = \frac{3}{2}\pi \log 2 - \pi$$

□

Elizalde [29] reported in 1985 that for $a > 0$

(1.45)     $\varsigma'(-1, a) =$

$$-\varsigma(-1, a)\log a - \frac{1}{4}a^2 + \frac{1}{12} - \frac{1}{2\pi^2}\sum_{n=1}^{\infty}\frac{1}{n^2}[\cos(2n\pi a)Ci(2n\pi a) + \sin(2n\pi a)si(2n\pi a)]$$

where $\varsigma(s, a)$ is the Hurwitz zeta function. Since $si(x) = Si(x) - \frac{\pi}{2}$ this may be written as

(1.46)   $\varsigma'(-1, a) = -\varsigma(-1, a)\log a - \frac{1}{4}a^2 + \frac{1}{12} + \frac{1}{4\pi}\sum_{n=1}^{\infty}\frac{\sin(2n\pi a)}{n^2}$

$$-\frac{1}{2\pi^2}\sum_{n=1}^{\infty}\frac{1}{n^2}[\cos(2n\pi a)Ci(2n\pi a) + \sin(2n\pi a)Si(2n\pi a)]$$

Elizalde [29] gave no indication of the source of the above identity but differentiation of (1.45) sheds more light on the subject: we have

$$\frac{\partial}{\partial a}\varsigma'(-1, a) = \left(a - \frac{1}{2}\right)\log a + \frac{1}{2}\left(a^2 - a + \frac{1}{6}\right)\frac{1}{a} - \frac{1}{2}a - \frac{1}{2\pi^2}\frac{\varsigma(2)}{a}$$

$$-\frac{1}{\pi}\sum_{n=1}^{\infty}\frac{1}{n}[-\sin(2n\pi a)Ci(2n\pi a) + \cos(2n\pi a)si(2n\pi a)]$$

since $\frac{d}{dx}Ci(x) = \frac{\cos x}{x}$ and $\frac{d}{dx}si(x) = \frac{\sin x}{x}$.

This simplifies to

$$\frac{\partial}{\partial a}\varsigma'(-1, a) = \left(a - \frac{1}{2}\right)\log a - \frac{1}{2} + \frac{1}{\pi}\sum_{n=1}^{\infty}\frac{1}{n}[\sin(2n\pi a)Ci(2n\pi a) - \cos(2n\pi a)si(2n\pi a)]$$

We have

$$\frac{\partial}{\partial a}\frac{\partial}{\partial s}\varsigma(s, a) = -\varsigma(s+1, a) - s\varsigma'(s+1, a)$$

and hence

$$\frac{\partial}{\partial a}\varsigma'(-1, a) = -\varsigma(0, a) + \varsigma'(0, a)$$

Then using Lerch's identity

(1.46.1) $$\varsigma'(0,a) = \log\Gamma(a) - \frac{1}{2}\log(2\pi)$$

this becomes

$$\frac{\partial}{\partial a}\varsigma'(-1,a) = -\varsigma(0,a) + \log\Gamma(a) - \frac{1}{2}\log(2\pi)$$

We then obtain

$$-\varsigma(0,a) + \log\Gamma(a) - \frac{1}{2}\log(2\pi)$$

$$= \left(a - \frac{1}{2}\right)\log a - \frac{1}{2} + \frac{1}{\pi}\sum_{n=1}^{\infty}\frac{1}{n}[\sin(2n\pi a)Ci(2n\pi a) - \cos(2n\pi a)si(2n\pi a)]$$

which, for $0 < a < 1$, simplifies to

(1.47) $$\log\Gamma(a) =$$

$$\frac{1}{2}\log(2\pi) + \left(a - \frac{1}{2}\right)\log a - a + \frac{1}{\pi}\sum_{n=1}^{\infty}\frac{1}{n}[\sin(2n\pi a)Ci(2n\pi a) - \cos(2n\pi a)si(2n\pi a)]$$

$$= \frac{1}{2}\log(2\pi) + \left(a - \frac{1}{2}\right)\log a - a + \frac{1}{2}\sum_{n=1}^{\infty}\frac{\cos(2n\pi a)}{n}$$

$$+ \frac{1}{\pi}\sum_{n=1}^{\infty}\frac{1}{n}[\sin(2n\pi a)Ci(2n\pi a) - \cos(2n\pi a)Si(2n\pi a)]$$

and using the Fourier series [55, p.148]

(1.47.1) $$\sum_{n=1}^{\infty}\frac{\cos(2n\pi a)}{n} = -\log[2\sin(\pi a)]$$

this becomes

(1.48) $$\log\Gamma(a) = \frac{1}{2}\log(2\pi) + \left(a - \frac{1}{2}\right)\log a - a - \frac{1}{2}\log[2\sin(\pi a)]$$

$$+ \frac{1}{\pi}\sum_{n=1}^{\infty}\frac{1}{n}[\sin(2n\pi a)Ci(2n\pi a) - \cos(2n\pi a)Si(2n\pi a)]$$

which corresponds with (1.41) above.

The formula (1.47) was given by Nörlund in [46, p.114]. When $a = 1/2$ we obtain

$$(1.49) \qquad \sum_{n=1}^{\infty} (-1)^n \frac{si(n\pi)}{n} = \frac{\pi}{2} \log 2 - \frac{\pi}{2}$$

This is a particular case of Nielsen's formula [44, p.83]

$$(1.49.1) \qquad \sum_{n=1}^{\infty} (-1)^n \frac{si(nx)}{n} = \frac{\pi}{2} \log 2 - \frac{x}{2}$$

It may be noted that Nielsen indicated that (1.49.1) was only valid for $x \in (-\pi, \pi)$.

With $a = 1/4$ in (1.48) we get

$$(1.50) \qquad \log \Gamma\left(\frac{1}{4}\right) = \frac{1}{2} \log(2\pi) + \frac{1}{2} \log 2 - \frac{1}{4} - \frac{1}{4} \log 2$$

$$+ \frac{1}{\pi} \sum_{n=1}^{\infty} \frac{1}{n} [\sin(n\pi/2) Ci(n\pi/2) - \cos(n\pi/2) Si(n\pi/2)]$$

We also have using (1.47)

$$(1.51) \qquad \log \Gamma(2a) = \left(2a - \frac{1}{2}\right) \log 2 + \log \Gamma(a) + a \log\left(a + \frac{1}{2}\right) - a - \frac{1}{2}$$

$$+ \frac{1}{\pi} \sum_{n=1}^{\infty} \frac{(-1)^n}{n} [\sin(2n\pi a) Ci(2n\pi a + n\pi) - \cos(2n\pi a) si(2n\pi a + n\pi)]$$

and with $a = 1/2$ we obtain

$$(1.52) \qquad \sum_{n=1}^{\infty} \frac{si(2n\pi)}{n} = \frac{\pi}{2} \log(2\pi) - \pi$$

and this concurs with Nielsen's result [44, p.79]. A completely different derivation of (1.54) is given in (1.105) below.

Using (1.47.1), which is valid for $0 < a < 1$, equation (1.48) may be written as

$$(1.53) \qquad \log \Gamma(a) = \frac{1}{2} \log(2\pi) + \left(a - \frac{1}{2}\right) \log a - a - \frac{1}{2} \log[2\sin(\pi a)]$$

$$+\frac{1}{\pi}\sum_{n=1}^{\infty}\frac{1}{n}[\sin(2n\pi a)Ci(2n\pi a)-\cos(2n\pi a)Si(2n\pi a)]$$

This may be written more compactly as

$$\log\left[\frac{\Gamma^2(a)\sin(\pi a)}{\pi a^{2a-1}}\right]=\frac{2}{\pi}\sum_{n=1}^{\infty}\frac{1}{n}[\sin(2n\pi a)Ci(2n\pi a)-\cos(2n\pi a)Si(2n\pi a)]$$

□

Integration of (1.47) results in

$$\int_{\varepsilon}^{x}\log\Gamma(a)=$$

$$\frac{1}{2}(x-\varepsilon)\log(2\pi)+\frac{1}{4}x\big[2-x+2(x-1)\log x\big]-\frac{1}{4}\varepsilon\big[2-\varepsilon+2(\varepsilon-1)\log\varepsilon\big]-\frac{1}{2}x^2+\frac{1}{2}\varepsilon^2$$

$$+\frac{1}{4\pi}\sum_{n=1}^{\infty}\frac{\sin(2n\pi x)}{n^2}-\frac{1}{4\pi}\sum_{n=1}^{\infty}\frac{\sin(2n\pi\varepsilon)}{n^2}$$

$$-\frac{1}{2\pi^2}\sum_{n=1}^{\infty}\frac{1}{n^2}\big[\cos(2n\pi x)Ci(2n\pi x)+\sin(2n\pi x)Si(2n\pi x)-\log(2n\pi x)\big]$$

$$+\frac{1}{2\pi^2}\sum_{n=1}^{\infty}\frac{1}{n^2}\big[\cos(2n\pi\varepsilon)Ci(2n\pi\varepsilon)+\sin(2n\pi\varepsilon)Si(2n\pi\varepsilon)-\log(2n\pi\varepsilon)\big]$$

Therefore as $\varepsilon\to 0$, and using (1.9.1) $\lim_{y\to 0}[\cos y\, Ci(y)-\log y]=\gamma$, we have

(1.60) $$\int_{0}^{x}\log\Gamma(a)\,da=$$

$$\frac{1}{2}x\log(2\pi)+\frac{1}{4}x\big[2-x+2(x-1)\log x\big]-\frac{1}{2}x^2+\frac{1}{4\pi}\sum_{n=1}^{\infty}\frac{\sin(2n\pi x)}{n^2}$$

$$-\frac{1}{2\pi^2}\sum_{n=1}^{\infty}\frac{1}{n^2}\big[\cos(2n\pi x)Ci(2n\pi x)+\sin(2n\pi x)Si(2n\pi x)-\log(2n\pi x)\big]+\frac{\gamma\varsigma(2)}{2\pi^2}$$

$$=\frac{1}{2}x\log(2\pi)+\frac{1}{4}x\big[2-x+2(x-1)\log x\big]-\frac{1}{2}x^2+\frac{1}{4\pi}\sum_{n=1}^{\infty}\frac{\sin(2n\pi x)}{n^2}+\frac{1}{12}\log x$$

$$-\frac{1}{2\pi^2}\sum_{n=1}^{\infty}\frac{1}{n^2}\Big[\cos(2n\pi x)Ci(2n\pi x)+\sin(2n\pi x)Si(2n\pi x)\Big]+\frac{\varsigma(2)\log(2\pi)}{2\pi^2}-\frac{\varsigma'(2)}{2\pi^2}+\frac{\gamma\varsigma(2)}{2\pi^2}$$

and with a little algebra and using the well known formula

$$(1.61)\qquad \varsigma'(-1)=\frac{1}{12}(1-\gamma-\log 2\pi)+\frac{1}{2\pi^2}\varsigma'(2)$$

(which is obtained by differentiating the functional equation for the Riemann zeta function) we obtain

$$(1.62)\qquad \int_0^x \log\Gamma(a)\,da=$$

$$=\frac{1}{2}x\log(2\pi)+\frac{1}{4}x\Big[2-x+2(x-1)\log x\Big]-\frac{1}{2}x^2+\frac{1}{4\pi}\sum_{n=1}^{\infty}\frac{\sin(2n\pi x)}{n^2}+\frac{1}{12}\log x$$

$$-\frac{1}{2\pi^2}\sum_{n=1}^{\infty}\frac{1}{n^2}\Big[\cos(2n\pi x)Ci(2n\pi x)+\sin(2n\pi x)Si(2n\pi x)\Big]+\frac{1}{12}-\varsigma'(-1)$$

Letting $x=1$ in (1.62) gives us

$$(1.63)\qquad \sum_{n=1}^{\infty}\frac{Ci(2n\pi)}{n^2}=2\pi^2\Big[\log A-\frac{1}{4}\Big]=-2\pi^2\Big[\varsigma'(-1)+\frac{1}{6}\Big]$$

and with $x=1/2$ we get

$$\int_0^{1/2}\log\Gamma(a)\,da=\frac{1}{4}\log(2\pi)+\frac{1}{16}+\frac{1}{24}\log 2-\frac{1}{2\pi^2}\sum_{n=1}^{\infty}\frac{(-1)^n}{n^2}Ci(n\pi)+\frac{1}{12}-\varsigma'(-1)$$

and using equation (6.117b) in [21]

$$\frac{1}{2\pi^2}\sum_{n=1}^{\infty}(-1)^n\frac{Ci(n\pi)}{n^2}=\frac{1}{12}\log 2+\frac{1}{48}+\frac{1}{2}\varsigma'(-1)$$

we see that

$$(1.64)\qquad \int_0^{1/2}\log\Gamma(a)\,da=\frac{5}{24}\log 2+\frac{1}{4}\log\pi+\frac{3}{2}\log A$$

which is reported in [52, p.35]. With $x=1/4$ we obtain

$$\int_0^{1/4} \log \Gamma(a)\, da =$$

$$= \frac{1}{8} \log(2\pi) + \frac{5}{64} + \frac{G}{4\pi} + \frac{1}{48} \log 2$$

$$-\frac{1}{2\pi^2} \sum_{n=1}^{\infty} \frac{1}{n^2} \left[ \cos(n\pi/2) Ci(n\pi/2) + \sin(n\pi/2) Si(n\pi/2) \right] + \frac{1}{12} - \varsigma'(-1)$$

From [52, p.35] we have

$$(1.54.1) \qquad \int_0^{1/4} \log \Gamma(a)\, da = \frac{1}{8} \log 2 + \frac{1}{8} \log \pi + \frac{9}{8} \log A + \frac{G}{4\pi}$$

and we therefore obtain

$$(1.65) \qquad \frac{1}{2\pi^2} \sum_{n=1}^{\infty} \frac{1}{n^2} \left[ \cos(n\pi/2) Ci(n\pi/2) + \sin(n\pi/2) Si(n\pi/2) \right] = \frac{5}{64} + \frac{1}{48} \log 2 - \frac{1}{8} \log A$$

Using (1.46) we may write (1.62) as Gosper's integral

$$(1.66) \qquad \int_0^x \log \Gamma(a)\, da = \frac{1}{2} x(1-x) + \frac{1}{2} x \log(2\pi) + \varsigma'(-1, x) - \varsigma'(-1)$$

In fact, we note that equating (1.62) with (1.66) results in Elizalde's formula (1.45).

$$\square$$

Integrating Elizalde's identity (1.45) gives us

$$\int_0^x \varsigma'(-1, a)\, da =$$

$$\frac{1}{72} x \left[ -4x^2 + 9x - 6 + 6(x-1)(2x-1) \log x \right] - \frac{1}{12} x^3 + \frac{1}{12} x - \frac{1}{8\pi^2} \sum_{n=1}^{\infty} \frac{\cos(2n\pi x)}{n^3}$$

$$+ \frac{1}{8\pi^2} \varsigma(3) - \frac{1}{4\pi^3} \sum_{n=1}^{\infty} \frac{1}{n^3} \left[ \sin(2n\pi x) Ci(2n\pi x) - \cos(2n\pi x) Si(2n\pi x) \right]$$

where we have used the integrals (easily derived using integration by parts)

$$\int Ci(x)\,dx = x\,Ci(x) - \sin x \qquad \text{and} \qquad \int Si(x)\,dx = x\,Si(x) + \cos x$$

In evaluating the integral at $a = 0$, we have used the fact that $Si(0) = 0$ and from (1.3) we have

$$\sin x\,Ci(x) = \gamma \sin x + \sin x \log x + \sin x \int_0^x \frac{\cos t - 1}{t}\,dt$$

$$= \gamma \sin x + \frac{\sin x}{x} x \log x + \sin x \int_0^x \frac{\cos t - 1}{t}\,dt$$

We therefore see that

$$\lim_{x \to 0} \sin x\,Ci(x) = 0$$

Since the Hurwitz zeta function is analytic in the whole complex plane except for $s \neq 1$, its partial derivatives commute in the region where the function is analytic: we therefore have

$$\frac{\partial}{\partial t}\frac{\partial}{\partial s}\varsigma(s,t) = \frac{\partial}{\partial s}\frac{\partial}{\partial t}\varsigma(s,t) = -\frac{\partial}{\partial s}[s\varsigma(s+1,t)]$$

$$= -\varsigma(s+1,t) - s\frac{\partial}{\partial s}\varsigma(s+1,t)$$

and upon integrating with respect to $t$ we see that

$$-s\int_0^v \varsigma'(s+1,t)\,dt = \int_0^v \frac{\partial}{\partial t}\frac{\partial}{\partial s}\varsigma(s,t)\,dt + \int_0^v \varsigma(s+1,t)\,dt$$

We therefore get

$$-s\int_0^v \varsigma'(s+1,t)\,dt = \varsigma'(s,v) - \varsigma'(s,0) + \int_0^v \varsigma(s+1,t)\,dt$$

and with $s = -n$ we have

$$n\int_0^v \varsigma'(1-n,t)\,du = \varsigma'(-n,v) - \varsigma'(-n,0) + \int_0^v \varsigma(1-n,t)\,dt$$

Then, using the well-known result [6, p.264]

$$\zeta(1-n,v) = -\frac{B_n(v)}{n} \text{ for } n \geq 1$$

we obtain

(1.66)
$$n\int_0^v \zeta'(1-n,t)\,dt = \frac{B_{n+1} - B_{n+1}(v)}{n(n+1)} + \zeta'(-n,v) - \zeta'(-n,0)$$

We have

$$\lim_{a\to 0}[\zeta(s,a) - a^{-s}] = \zeta(s)$$

$$\lim_{a\to 0}[\zeta(-s,a) - a^{s}] = \zeta(-s)$$

$$\frac{\partial}{\partial s}\lim_{a\to 0}[\zeta(-s,a) - a^{s}] = -\zeta'(-s)$$

$$\frac{\partial}{\partial s}\lim_{a\to 0}[\zeta(-s,a) - a^{s}] = \lim_{a\to 0}[-\zeta'(-s,a) - a^{s}\log a]$$

$$= -\zeta'(-s,0)$$

Therefore we have

$$\zeta'(-s,0) = \zeta'(-s)$$

and thus we obtain

(1.67)
$$n\int_0^v \zeta'(1-n,t)\,dt = \frac{B_{n+1} - B_{n+1}(v)}{n(n+1)} + \zeta'(-n,v) - \zeta'(-n)$$

This integral was originally derived by Adamchik [2] in a different manner in 1998.

We have for $n = 2$

$$\int_0^x \zeta'(-1,a)\,da = -\frac{1}{12}B_3(x) + \frac{1}{2}\zeta'(-2,x) - \frac{1}{2}\zeta'(-2)$$

Therefore we obtain

$$\frac{1}{72}x[-4x^2 + 9x - 6 + 6(x-1)(2x-1)\log x] - \frac{1}{12}x^3 + \frac{1}{12}x - \frac{1}{8\pi^2}\sum_{n=1}^{\infty}\frac{\cos(2n\pi x)}{n^3}$$

$$+ \frac{1}{8\pi^2}\zeta(3) - \frac{1}{4\pi^3}\sum_{n=1}^{\infty}\frac{1}{n^3}[\sin(2n\pi x)Ci(2n\pi x) - \cos(2n\pi x)Si(2n\pi x)]$$

$$= -\frac{1}{12}B_3(x) + \frac{1}{2}\varsigma'(-2,x) - \frac{1}{2}\varsigma'(-2)$$

This is easily simplified to

(1.68)

$$\frac{1}{12}x(x-1)(2x-1)\log x - \frac{5}{36}x^3 + \frac{1}{8}x^2 - \frac{1}{8\pi^2}\sum_{n=1}^{\infty}\frac{\cos(2n\pi x)}{n^3} + \frac{1}{8\pi^2}\varsigma(3)$$

$$-\frac{1}{4\pi^3}\sum_{n=1}^{\infty}\frac{1}{n^3}[\sin(2n\pi x)Ci(2n\pi x) - \cos(2n\pi x)Si(2n\pi x)]$$

$$= -\frac{1}{12}B_3(x) + \frac{1}{2}\varsigma'(-2,x) - \frac{1}{2}\varsigma'(-2)$$

Equation (1.68) could then be integrated to produce an identity involving $\varsigma'(-3,x)$ and so on.

$\square$

Differentiating (1.41) term by term (and boldly assuming that the procedure is valid) easily results in

(1.70)     $\psi(a) = \log a - \frac{1}{2a} + 2\sum_{n=1}^{\infty}[\cos(2n\pi a)Ci(2n\pi a) + \sin(2n\pi a)si(2n\pi a)]$

which appears in Nörlund's book [46, p.108]. Letting $a = 1$ results in

(1.71)     $\frac{1}{2} - \gamma = 2\sum_{n=1}^{\infty}Ci(2n\pi)$

and this corrects the corresponding formula given by Nielsen [44, p.80]. It appears that Nielsen's analysis is incorrect because he effectively used the Fourier series (1.41.1) which does not hold at the point $a = 1$.

This formula was also used in [20]. With $a = 1/2$ we get

$$\psi\left(\frac{1}{2}\right) = -\log 2 - 1 + 2\sum_{n=1}^{\infty}(-1)^n Ci(n\pi)$$

and hence we have

(1.72)
$$2\sum_{n=1}^{\infty}(-1)^n Ci(n\pi) = 1 - \gamma - \log 2$$

With $a \to a/2$ we see that

(1.73)
$$\psi(a/2) = \log(a/2) - \frac{1}{a} + 2\sum_{n=1}^{\infty}[\cos(n\pi a)Ci(n\pi a) + \sin(n\pi a)si(n\pi a)]$$

and we also have

$$\psi(a+1) = \log(a+1) - \frac{1}{2(a+1)}$$

$$+ 2\sum_{n=1}^{\infty}\cos[2n\pi(a+1)]Ci[2n\pi(a+1)] + \sin[2n\pi(a+1)]si[2n\pi(a+1)]$$

Equation (1.70) may be expressed as

$$\psi(a) = \log a - \frac{1}{2a} + 2\sum_{n=1}^{\infty} a_n$$

and we have with $a \to 2a$

$$\psi(2a) = \log(2a) - \frac{1}{4a} + 2\sum_{n=1}^{\infty} a_{2n}$$

Since $2\sum_{n=1}^{\infty} a_{2n} = \sum_{n=1}^{\infty} a_n + \sum_{n=1}^{\infty}(-1)^n a_n$ we have

$$\psi(2a) = \log(2a) - \frac{1}{4a} + \frac{1}{2}\psi(a) - \frac{1}{2}\log a + \frac{1}{4a} + \sum_{n=1}^{\infty}(-1)^n a_n$$

and using the duplication formula for the gamma function [52, p.15] (which may be obtained by differentiating (2.16.2))

$$\psi(2a) = \log 2 + \frac{1}{2}\psi(a) + \frac{1}{2}\psi\left(a + \frac{1}{2}\right)$$

this results in

(1.74)
$$\psi\left(a + \frac{1}{2}\right) = \log a + 2\sum_{n=1}^{\infty}(-1)^n[\cos(n\pi a)Ci(n\pi a) + \sin(n\pi a)si(n\pi a)]$$

Differentiating (1.41.1) results in (which is valid for $0 < a < 1$)

$$(1.75) \quad \psi(a) = -\frac{1}{a} + 2\sum_{n=1}^{\infty} \left[\cos 2n\pi a \, Ci(2n\pi) + \sin 2n\pi a \, si(2n\pi)\right]$$

or equivalently

$$\psi(a+1) = 2\sum_{n=1}^{\infty} \left[\cos 2n\pi a \, Ci(2n\pi) + \sin 2n\pi a \, si(2n\pi)\right]$$

which may be compared with (1.70)

$$\psi(a) = \log a - \frac{1}{2a} + 2\sum_{n=1}^{\infty} \left[\cos(2n\pi a) Ci(2n\pi a) + \sin(2n\pi a) si(2n\pi a)\right]$$

$\square$

Differentiating (1.1) with respect to $a$

$$\int_0^1 \log \Gamma(x+a) \sin 2n\pi x \, dx = -\frac{1}{2\pi n}\left[\log a - \cos(2n\pi a) Ci(2n\pi a) - \sin(2n\pi a) si(2n\pi a)\right]$$

we obtain

$$\int_0^1 \psi(x+a) \sin 2n\pi x \, dx = -\frac{1}{2\pi n}\left[\begin{array}{l} \dfrac{1}{a} - \cos(2n\pi a)\dfrac{\cos(2n\pi a)}{a} + 2n\pi \sin(2n\pi a) Ci(2n\pi a) \\[2mm] -\sin(2n\pi a)\dfrac{\sin(2n\pi a)}{a} - 2n\pi \cos(2n\pi a) si(2n\pi a) \end{array}\right]$$

and therefore we have

$$(1.76) \quad \int_0^1 \psi(x+a) \sin 2n\pi x \, dx = -\sin(2n\pi a) Ci(2n\pi a) + \cos(2n\pi a) si(2n\pi a)$$

which corrects the entry in [33, p.652, 6.467 1].

With $a = 1$ we get

$$(1.77) \quad \int_0^1 \psi(x+1) \sin 2n\pi x \, dx = si(2n\pi)$$

and with $a = 0$ we get

$(1.77.1)$ $\displaystyle\int_0^1 \psi(x)\sin 2n\pi x\,dx = si(0) = -\frac{\pi}{2}$

which is correctly reported in [33, p.652, 6.465 2].

We see that

$$\int_0^1 \psi(x+1)\sin 2n\pi x\,dx = \int_0^1 \psi(x)\sin 2n\pi x\,dx + \int_0^1 \frac{\sin 2n\pi x}{x}\,dx$$

$$= \int_0^1 \psi(x)\sin 2n\pi x\,dx + Si(2n\pi) - Si(0)$$

and this concurs with (1.77) and (1.77.1).

Similarly, differentiating (1.27) we obtain

$(1.78)$ $\displaystyle\int_0^1 \psi(x+a)\cos 2n\pi x\,dx = \sin(2n\pi a)si(2n\pi a) + \cos(2n\pi a)Ci(2n\pi a)$

which corrects the entry in [33, p.652, 6.467 2]. With $a = 1$ we have

$(1.78.1)$ $\displaystyle\int_0^1 \psi(x+1)\cos 2n\pi x\,dx = Ci(2n\pi)$

$\square$

Using integration by parts we see that for $a \geq 0$

$$\int_0^1 \psi(x+a)\sin 2n\pi x\,dx = \sin 2n\pi x \log\Gamma(x+a)\Big|_0^1 - 2n\pi\int_0^1 \Gamma(x+a)\cos 2n\pi x\,dx$$

We have for $a \geq 0$

$$\lim_{x\to 0}\sin 2n\pi x \log\Gamma(x+a) = \lim_{x\to 0}\frac{\sin 2n\pi x}{x}x\log\Gamma(x+a) = 0$$

and hence we have for $a \geq 0$

$(1.79)$ $\displaystyle\int_0^1 \psi(x+a)\sin 2n\pi x\,dx = -2n\pi\int_0^1 \Gamma(x+a)\cos 2n\pi x\,dx$

Using (1.27)

$$\int_0^1 \log\Gamma(x+a)\cos 2n\pi x\,dx = -\frac{1}{2n\pi}\big[\sin(2n\pi a)Ci(2n\pi a)+\cos(2n\pi a)si(2n\pi a)\big]$$

we again obtain

$$\int_0^1 \psi(x+a)\sin 2n\pi x\,dx = \sin(2n\pi a)Ci(2n\pi a)+\cos(2n\pi a)si(2n\pi a)$$

Similarly we have for $a > 0$

(1.80) $\qquad \int_0^1 \psi(x+a)\cos 2n\pi x\,dx = \log\Gamma(1+a)-\log\Gamma(a)+2n\pi\int_0^1 \Gamma(x+a)\sin 2n\pi x\,dx$

and using (1.1)

$$\int_0^1 \log\Gamma(x+a)\sin 2n\pi x\,dx = -\frac{1}{2n\pi}\big[\log a-\cos(2n\pi a)Ci(2n\pi a)-\sin(2n\pi a)si(2n\pi a)\big]$$

we again obtain

$$\int_0^1 \psi(x+a)\cos 2n\pi x\,dx = \cos(2n\pi a)Ci(2n\pi a)+\sin(2n\pi a)si(2n\pi a)$$

Nielsen [44, p.80] reports these integrals in the case $a = 1$ as

(1.81) $\qquad \int_0^1 \psi(x+1)\sin 2n\pi x\,dx = -2n\pi\int_0^1 \log\Gamma(x+1)\cos 2n\pi x\,dx$

(1.82) $\qquad \int_0^1 \psi(x+1)\cos 2n\pi x\,dx = 2n\pi\int_0^1 \log\Gamma(x+1)\sin 2n\pi x\,dx$

$\square$

Applying Parseval's identity with (1.78) and (1.78.1) we obtain

$$2\int_0^1 \psi^2(x+1)\,dx = \sum_{n=1}^{\infty}[Ci^2(2n\pi)+si^2(2n\pi)]$$

With regard to this, see also Lewin's monograph [43, p.24] where he showed that

$$\int_0^\infty e^{-ux} Li_2(-x^2)\, dx = -\frac{2}{u}\Big[ Ci^2(u) + si^2(u) \Big]$$

and with $u = 2n\pi$ we have

$$n\pi \int_0^\infty e^{-2n\pi x} Li_2(-x^2)\, dx = -\Big[ Ci^2(2n\pi) + si^2(2n\pi) \Big]$$

We make the summation

$$\pi \int_0^\infty \sum_{n=1}^\infty n\, e^{-2n\pi x} Li_2(-x^2)\, dx = -\sum_{n=1}^\infty [Ci^2(2n\pi) + si^2(2n\pi)]$$

Using the derivative of the geometric series we have

$$\sum_{n=1}^\infty n\, z^n = \frac{z}{(1-z)^2}$$

so that

$$\pi \int_0^\infty \sum_{n=1}^\infty n\, e^{-2n\pi x} Li_2(-x^2)\, dx = \pi \int_0^\infty \frac{e^{-2\pi x}}{(1-e^{-2\pi x})^2} Li_2(-x^2)\, dx$$

and hence we have the curious equality

$$\pi \int_0^\infty \frac{e^{-2\pi x}}{(1-e^{-2\pi x})^2} Li_2(-x^2)\, dx = -2\int_0^1 \psi^2(x+1)\, dx$$

A different representation of $Ci^2(u) + si^2(u)$ as an integral is given in [44, p.32].

**Applications of Nielsen's representation of $\log\Gamma(x)$**

We recall Nielsen's representation of $\log\Gamma(x)$ in equation (1.41.1)

$$\log\Gamma(x) = \frac{1}{2}\log(2\pi) - 1 - \log x + \frac{1}{\pi}\sum_{n=1}^\infty \frac{1}{n}[\sin 2n\pi x\, Ci(2n\pi) - \cos 2n\pi x\, si(2n\pi)]$$

and multiplying this by $\cos p\pi x$ and integrating gives us

(1.90) $\quad \int\limits_{0}^{1} \log \Gamma(x) \cos p\pi x \, dx = \left[ \frac{1}{2} \log(2\pi) - 1 \right] \frac{\sin p\pi}{p\pi} + \frac{Si(p\pi)}{p\pi}$

$$+ \frac{2(1 - \cos p\pi)}{\pi^2} \sum_{n=1}^{\infty} \frac{Ci(2n\pi)}{4n^2 - p^2} + \frac{p \sin p\pi}{\pi^2} \sum_{n=1}^{\infty} \frac{1}{n} \frac{si(2n\pi)}{4n^2 - p^2}$$

where we have employed the integral

$$\int \log x \cos p\pi x \, dx = \frac{1}{p\pi} [\log x \sin p\pi x - Si(p\pi x)]$$

(1.91) $\quad \int\limits_{0}^{1} \log x \cos p\pi x \, dx = -\frac{Si(p\pi)}{p\pi}$

Letting $p = 0$ in (1.90) results in

$$\int\limits_{0}^{1} \log \Gamma(x) \, dx = \left[ \frac{1}{2} \log(2\pi) - 1 \right] + \lim_{p \to 0} \frac{Si(p\pi)}{p\pi}$$

and using L'Hôpital's rule we see that $\lim\limits_{p \to 0} \dfrac{Si(p\pi)}{p\pi} = 1$ giving us the familiar result

$$\int\limits_{0}^{1} \log \Gamma(x) \, dx = \frac{1}{2} \log(2\pi)$$

With $p = 1$ we have

(1.92) $\quad \int\limits_{0}^{1} \log \Gamma(x) \cos \pi x \, dx = \frac{Si(\pi)}{\pi} + \frac{4}{\pi^2} \sum_{n=1}^{\infty} \frac{Ci(2n\pi)}{4n^2 - 1}$

and $p = 2k$ results in

$$\int\limits_{0}^{1} \log \Gamma(x) \cos 2k\pi x \, dx = \frac{Si(2k\pi)}{2k\pi} - \frac{si(2k\pi)}{2k\pi} = \frac{1}{4k}$$

in agreement with (3.10). The following limits have been used in the above calculation

$$\lim_{p \to 2k} \frac{2(1 - \cos p\pi)}{\pi^2} \sum_{n=1}^{\infty} \frac{Ci(2n\pi)}{4n^2 - p^2} = \lim_{p \to 2k} \frac{2\pi \sin p\pi \, Ci(2k\pi)}{\pi^2(-2p)} = 0$$

$$\lim_{p \to 2k} \frac{p \sin p\pi}{\pi^2} \sum_{n=1}^{\infty} \frac{1}{n} \frac{si(2n\pi)}{4n^2 - p^2} = \lim_{p \to 2k} \frac{(p\pi \cos p\pi + \sin p\pi)si(2k\pi)}{k\pi^2(-2p)} = -\frac{si(2k\pi)}{2k\pi}$$

With $p = 2k + 1$ we have

(1.93) $\displaystyle\int_0^1 \log \Gamma(x) \cos(2k+1)\pi x \, dx = \frac{Si[(2k+1)\pi]}{(2k+1)\pi} + \frac{4}{\pi^2} \sum_{n=1}^{\infty} \frac{Ci(2n\pi)}{4n^2 - (2k+1)^2}$

We see from the definition (1.3) that

$$Ci(2n\pi x) = \gamma + \log(2n\pi x) + \int_0^{2n\pi x} \frac{\cos t - 1}{t} dt$$

and we have the summation

$$\sum_{n=1}^{\infty} \frac{Ci(2n\pi x)}{4n^2 - 1} = [\gamma + \log(2\pi x)] \sum_{n=1}^{\infty} \frac{1}{4n^2 - 1} + \sum_{n=1}^{\infty} \frac{\log n}{4n^2 - 1} + \sum_{n=1}^{\infty} \frac{1}{4n^2 - 1} \left[ \int_0^{2n\pi x} \frac{\cos t - 1}{t} dt \right]$$

$$= \frac{1}{2}[\gamma + \log(2\pi x)] + \sum_{n=1}^{\infty} \frac{\log n}{4n^2 - 1} + \sum_{n=1}^{\infty} \frac{1}{4n^2 - 1} \left[ \int_0^{2n\pi x} \frac{\cos t - 1}{t} dt \right]$$

Using $\displaystyle\int_0^{2n\pi x} \frac{\cos t - 1}{t} dt = \int_0^x \frac{\cos(2n\pi y) - 1}{y} dy$ the last series becomes

$$\sum_{n=1}^{\infty} \frac{1}{4n^2 - 1} \left[ \int_0^{2n\pi x} \frac{\cos t - 1}{t} dt \right] = \sum_{n=1}^{\infty} \frac{1}{4n^2 - 1} \left[ \int_0^x \frac{\cos(2n\pi y) - 1}{y} dy \right]$$

$$= \sum_{n=1}^{\infty} \frac{1}{4n^2 - 1} \lim_{a \to 0} \left[ \int_0^x \frac{\cos(2n\pi y) - 1}{y} dy \right]$$

$$= \lim_{a \to 0} \int_a^x \sum_{n=1}^{\infty} \frac{\cos(2n\pi y) - 1}{(4n^2 - 1)y} dy$$

$$= \lim_{a \to 0} \int_a^x \left[ -\frac{1}{2y} + \sum_{n=1}^{\infty} \frac{\cos(2n\pi y)}{(4n^2 - 1)y} \right] dy$$

It is an exercise in Apostol's book [7, p.337] to show that for $0 < x < \pi$

$$(1.94) \qquad \sin x = \frac{2}{\pi} - \frac{4}{\pi} \sum_{n=1}^{\infty} \frac{\cos(2nx)}{4n^2 - 1}$$

and in fact, since $\sum_{n=1}^{\infty} \frac{1}{4n^2 - 1} = \frac{1}{2}$, this holds true for $0 \le x < \pi$. It seems likely that (1.94) may also be determined using the differentiation theorems for Fourier series appearing in Tolstov's book [55, p.142].

We have

$$\frac{1}{4}(2 - \pi \sin \pi y) = \sum_{n=1}^{\infty} \frac{\cos(2n\pi y)}{4n^2 - 1}$$

so that

$$\lim_{a \to 0} \int_a^x \left[ -\frac{1}{2y} + \sum_{n=1}^{\infty} \frac{\cos(2n\pi y)}{(4n^2 - 1)y} \right] dy = -\lim_{a \to 0} \frac{1}{4} \int_a^x \frac{\pi \sin \pi y}{y} dy = -\frac{\pi}{4} Si(\pi x)$$

Hence we see that

$$(1.95) \qquad \sum_{n=1}^{\infty} \frac{Ci(2n\pi x)}{4n^2 - 1} = \frac{1}{2} [\gamma + \log(2\pi x)] + \sum_{n=1}^{\infty} \frac{\log n}{4n^2 - 1} - \frac{\pi}{4} Si(\pi x)$$

and referring to (1.93)

$$\int_0^1 \log \Gamma(x) \cos \pi x \, dx = \frac{Si(\pi)}{\pi} + \frac{4}{\pi^2} \sum_{n=1}^{\infty} \frac{Ci(2n\pi)}{4n^2 - 1}$$

we obtain

$$\int_0^1 \log \Gamma(x) \cos \pi x \, dx = \frac{2}{\pi^2} \left[ \log(2\pi) + \gamma + 2 \sum_{n=1}^{\infty} \frac{\log n}{4n^2 - 1} \right]$$

which we shall derive in a different manner in (4.5.2).

We may write (1.95) as

$$\sum_{n=1}^{\infty} \frac{Ci(2n\pi x) - \log(2n\pi x)}{4n^2 - 1} + \sum_{n=1}^{\infty} \frac{\log(2n\pi x)}{4n^2 - 1} - \log(2\pi x) \sum_{n=1}^{\infty} \frac{1}{4n^2 - 1} + \frac{\pi}{4} Si(\pi x) = \frac{1}{2}\gamma + \sum_{n=1}^{\infty} \frac{\log n}{4n^2 - 1}$$

and take the limit as $x \to 0$. Having regard to (1.3) we see that the limit equates to the right-hand side.

Comparing (1.93) with (4.6) results in

$(1.95.1)$ $\qquad \dfrac{4}{\pi^2}\displaystyle\sum_{n=1}^{\infty}\dfrac{\gamma+\log(2\pi n)}{4n^2-(2k+1)^2}=\dfrac{Si[(2k+1)\pi]}{(2k+1)\pi}+\dfrac{4}{\pi^2}\sum_{n=1}^{\infty}\dfrac{Ci(2n\pi)}{4n^2-(2k+1)^2}$

Integrating (1.95) results in

$$x\sum_{n=1}^{\infty}\frac{Ci(2n\pi x)}{4n^2-1}-\frac{1}{2\pi}\sum_{n=1}^{\infty}\frac{1}{n}\frac{\sin(2n\pi x)}{4n^2-1}=\frac{x}{2}[\gamma+\log(2\pi x)]-\frac{x}{2}+x\sum_{n=1}^{\infty}\frac{\log n}{4n^2-1}$$

$$-\frac{\pi}{4}x\,Si(\pi x)-\frac{1}{4}[\cos(\pi x)-1]$$

and hence we obtain

$(1.96)$ $\qquad \displaystyle\sum_{n=1}^{\infty}\dfrac{1}{n}\dfrac{\sin(2n\pi x)}{4n^2-1}=\dfrac{\pi}{2}[2x+\cos(\pi x)-1]$

which concurs with the result directly obtained by integrating (1.94).

In passing, we note that it is an exercise in Bromwich's book [15, p.371] to show that for $0<x<\pi$

$(1.96.1)$ $\qquad \displaystyle\sum_{n=0}^{\infty}\dfrac{\cos nx}{(2n+1)^2-a^2}=\dfrac{\pi\sin a(\frac{\pi}{2}-x)}{4a\cos\frac{a\pi}{2}}$

$\square$

Multiplying (1.41.1) by $\sin p\pi x$ and integrating gives us

$(1.97)$ $\displaystyle\int_{0}^{1}\log\Gamma(x)\sin p\pi x\,dx=\left[\dfrac{1}{2}\log(2\pi)-1\right]\dfrac{1-\cos p\pi}{p\pi}-\dfrac{Ci(p\pi)}{p\pi}+\dfrac{\gamma+\log(p\pi)}{p\pi}$

$$-\frac{2\sin p\pi}{\pi^2}\sum_{n=1}^{\infty}\frac{Ci(2n\pi)}{4n^2-p^2}+\frac{p(1-\cos p\pi)}{\pi^2}\sum_{n=1}^{\infty}\frac{1}{n}\frac{si(2n\pi)}{4n^2-p^2}$$

where we have employed the integral (1.10.3)

$$\int_{0}^{1}\log x\sin p\pi x\,dx=\frac{Ci(p\pi)-\gamma-\log(p\pi)}{p\pi}$$

With $p=1$ we have

$$(1.98) \quad \int_0^1 \log \Gamma(x) \sin \pi x \, dx = \left[\frac{1}{2}\log(2\pi) - 1\right]\frac{2}{\pi} - \frac{Ci(\pi)}{\pi} + \frac{\gamma + \log \pi}{\pi} + \frac{2}{\pi^2}\sum_{n=1}^{\infty}\frac{1}{n}\frac{si(2n\pi)}{4n^2 - 1}$$

and comparing this with (5.7) we deduce that

$$(1.99) \quad \sum_{n=1}^{\infty}\frac{1}{n}\frac{si(2n\pi)}{4n^2 - 1} = \frac{\pi}{2}\left[3 + Ci(\pi) - \gamma - \log(4\pi)\right]$$

Letting $p = 2k$ results in

$$\int_0^1 \log \Gamma(x) \sin 2k\pi x \, dx$$

$$= -\frac{Ci(2k\pi)}{2k\pi} + \frac{\gamma + \log(2k\pi)}{2k\pi} - \lim_{p \to 2k}\frac{2\pi \cos p\pi \, Ci(2k\pi)}{\pi^2(-2p)} + \lim_{p \to 2k}\frac{(1 + p\pi \sin p\pi - \cos p\pi)si(2k\pi)}{\pi^2 k(-2p)}$$

and hence we obtain another derivation of

$$\int_0^1 \log \Gamma(x) \sin 2k\pi x \, dx = \frac{\gamma + \log(2\pi k)}{2\pi k}$$

With $p = 2k + 1$ we have

$(1.100)$

$$\int_0^1 \log \Gamma(x) \sin(2k+1)\pi x \, dx = \left[\frac{1}{2}\log(2\pi) - 1\right]\frac{2}{(2k+1)\pi} - \frac{Ci[(2k+1)\pi]}{(2k+1)\pi} + \frac{\gamma + \log[(2k+1)\pi]}{(2k+1)\pi}$$

$$+ \frac{2(2k+1)}{\pi^2}\sum_{n=1}^{\infty}\frac{1}{n}\frac{si(2n\pi)}{4n^2 - (2k+1)^2}$$

which may be compared with (5.6)

$$\int_0^1 \log \Gamma(x) \sin(2k+1)\pi x \, dx = \frac{1}{(2k+1)\pi}\left[\log\left(\frac{\pi}{2}\right) + \frac{1}{2k+1} + 2\sum_{j=0}^{k-1}\frac{1}{2j+1}\right]$$

Further identities may be obtained by multiplying (1.41.1) by $x^n \sin p\pi x$ or $x^n \cos p\pi x$ and integrating as before.

Equation (1.99) may also be derived as follows. Referring to (1.91) we have

$$\int\limits_0^1 \log x \cos 2n\pi x \, dx = -\frac{Si(2n\pi)}{2n\pi}$$

and we make the summation

$$\sum_{n=1}^\infty \int\limits_0^1 \log x \frac{\cos 2n\pi x}{4n^2-1} \, dx = -\frac{1}{2\pi}\sum_{n=1}^\infty \frac{1}{n}\frac{Si(2n\pi)}{4n^2-1}$$

$$= -\frac{1}{2\pi}\sum_{n=1}^\infty \frac{1}{n}\frac{si(2n\pi)}{4n^2-1} - \frac{1}{4}\sum_{n=1}^\infty \frac{1}{n}\frac{1}{4n^2-1}$$

and using (2.15) this becomes

$$= -\frac{1}{2\pi}\sum_{n=1}^\infty \frac{1}{n}\frac{si(2n\pi)}{4n^2-1} - \frac{1}{2}\log 2 + \frac{1}{4}$$

Assuming that the interchange of summation and integration is valid

$$\sum_{n=1}^\infty \int\limits_0^1 \log x \frac{\cos 2n\pi x}{4n^2-1} \, dx = \int\limits_0^1 \log x \sum_{n=1}^\infty \frac{\cos 2n\pi x}{4n^2-1} \, dx$$

Employing (1.94) this becomes

$$= \frac{1}{4}\int\limits_0^1 \log x (2 - \pi \sin \pi x) \, dx$$

and using (1.10.3) we have

$$= -\frac{1}{4}[2 + Ci(\pi) - \gamma - \log \pi]$$

Therefore we obtain

(1.101) $$\sum_{n=1}^\infty \frac{1}{n}\frac{si(2n\pi)}{4n^2-1} = \frac{\pi}{2}[3 + Ci(\pi) - \gamma - \log(4\pi)]$$

Alternatively we consider

$$\sum_{n=1}^\infty \frac{1}{n}\frac{si(2n\pi)}{4n^2-1} = \sum_{n=1}^\infty \frac{1}{n}\frac{Si(2n\pi)}{4n^2-1} - \frac{\pi}{2}\sum_{n=1}^\infty \frac{1}{n}\frac{1}{4n^2-1}$$

$$= \sum_{n=1}^{\infty} \frac{1}{n} \frac{Si(2n\pi)}{4n^2-1} - \frac{\pi}{2}(2\log 2 - 1)$$

where we have used (2.15).

$$\sum_{n=1}^{\infty} \frac{1}{n} \frac{Si(2n\pi)}{4n^2-1} = \int_0^{2\pi} \sum_{n=1}^{\infty} \frac{1}{n} \frac{\sin nt}{4n^2-1} \frac{dt}{t}$$

$$= \int_0^1 \sum_{n=1}^{\infty} \frac{1}{n} \frac{\sin(2n\pi t)}{4n^2-1} \frac{dt}{t}$$

and using (1.96.1) this becomes

$$= \frac{\pi}{2} \int_0^1 \frac{\cos(\pi t) - 1 + 2t}{t} dt$$

We note from (1.3) that

$$Ci(\pi x) = \gamma + \log \pi + \log x + \int_0^{\pi x} \frac{\cos t - 1}{t} dt$$

$$= \gamma + \log \pi + \log x + \int_0^x \frac{\cos(\pi t) - 1}{t} dt$$

Therefore we obtain another derivation of

$$\sum_{n=1}^{\infty} \frac{1}{n} \frac{Si(2n\pi)}{4n^2-1} = \frac{\pi}{2}[3 + Ci(\pi) - \gamma - \log(4\pi)]$$

$\square$

Equating (1.90) and (4.4.1) gives us

$$\left[\frac{1}{2}\log(2\pi) - 1\right]\frac{\sin p\pi}{p\pi} + \frac{Si(p\pi)}{p\pi} + \frac{2(1-\cos p\pi)}{\pi^2}\sum_{n=1}^{\infty}\frac{Ci(2n\pi)}{4n^2-p^2} + \frac{p\sin p\pi}{\pi^2}\sum_{n=1}^{\infty}\frac{1}{n}\frac{si(2n\pi)}{4n^2-p^2}$$

$$= \frac{\sin p\pi[\gamma + \log(2\pi)]}{2p\pi} + \frac{\sin p\pi}{4p\pi}\left[\psi\left(\frac{p}{2}\right) + \psi\left(-\frac{p}{2}\right)\right]$$

$$+ \frac{2(1 - \cos p\pi)[\gamma + \log(2\pi)]}{\pi^2}\left[\frac{1}{2p^2} - \frac{\pi}{4p}\cot\left(\frac{p\pi}{2}\right)\right] + \frac{2(1 - \cos p\pi)}{\pi^2}\sum_{n=1}^{\infty}\frac{\log n}{4n^2 - p^2}$$

Equating (1.97) and (5.5.1) results in

$$\left[\frac{1}{2}\log(2\pi) - 1\right]\frac{1 - \cos p\pi}{p\pi} - \frac{Ci(p\pi)}{p\pi} + \frac{\gamma + \log(p\pi)}{p\pi} - \frac{2\sin p\pi}{\pi^2}\sum_{n=1}^{\infty}\frac{Ci(2n\pi)}{4n^2 - p^2}$$

$$+ \frac{p(1 - \cos p\pi)}{\pi^2}\sum_{n=1}^{\infty}\frac{1}{n}\frac{si(2n\pi)}{4n^2 - p^2}$$

$$= \frac{(1 - \cos p\pi)[\gamma + \log(2\pi)]}{2p\pi} + \frac{(1 - \cos p\pi)}{4p\pi}\left[\psi\left(1 + \frac{p}{2}\right) + \psi\left(1 - \frac{p}{2}\right)\right]$$

$$- \frac{2\sin p\pi[\gamma + \log(2\pi)]}{\pi^2}\left[\frac{1}{2p^2} - \frac{\pi}{4p}\cot\left(\frac{p\pi}{2}\right)\right] - \frac{2\sin p\pi}{\pi^2}\sum_{n=1}^{\infty}\frac{\log n}{4n^2 - p^2}$$

and hence we have two simultaneous equations involving $\sum_{n=1}^{\infty}\dfrac{Ci(2n\pi)}{4n^2 - p^2}$ and

$\sum_{n=1}^{\infty}\dfrac{1}{n}\dfrac{si(2n\pi)}{4n^2 - p^2}$ .

$\square$

We showed in [21] that

$$(1.102) \qquad \int_a^b p(x)\cot x\, dx = 2\sum_{n=1}^{\infty}\int_a^b p(x)\sin 2nx\, dx$$

which is valid for a wide class of suitably behaved functions. Specifically we require that $p(x)$ is a twice continuously differentiable function and that either (i) both $\sin(x/2)$ and $\cos(x/2)$ have no zero in $[a,b]$ or (ii) if either $\sin(a/2)$ or $\cos(a/2)$ is equal to zero then $p(a)$ must also be zero. Condition (i) is equivalent to the requirement that $\sin x$ has no zero in $[a,b]$.

Since $\log\Gamma(1) = \log\Gamma(2) = 0$, it is clear that $p(x) = \log\Gamma(x+1)$ satisfies the necessary conditions on the interval $[0,1]$. Therefore we have

$$\int_0^1 \log\Gamma(x+1)\cot\pi x\, dx = 2\sum_{n=1}^{\infty}\int_0^1 \log\Gamma(x+1)\sin 2n\pi x\, dx$$

and using (1.11) this becomes

$$(1.103) \qquad \int_0^1 \log \Gamma(x+1) \cot \pi x \, dx = \sum_{n=1}^{\infty} \frac{Ci(2n\pi)}{n\pi}$$

This is an interesting result if only for the fact that *Mathematica* is not able to evaluate even the seemingly simpler integral $\int_0^1 \log x \cot \pi x \, dx$.

We note that integration by parts gives us

$$\int_0^1 \log \Gamma(x+1) \cot \pi x \, dx = \frac{1}{\pi} \log \Gamma(x+1) \log \sin \pi x \Big|_0^1 - \frac{1}{\pi} \int_0^1 \psi(1+x) \log \sin \pi x \, dx$$

and using L'Hôpital's rule we obtain

$$\lim_{x \to 1} [\log \Gamma(x+1) \log \sin \pi x] = \lim_{x \to 1} \left[ \frac{\log \sin \pi x}{1/\log \Gamma(x+1)} \right]$$

$$= \lim_{x \to 1} \left[ -\frac{\pi \cot \pi x \log^2 \Gamma(x+1)}{\psi(x+1)} \right]$$

$$= \lim_{x \to 1} \left[ -\frac{\pi \log^2 \Gamma(x+1)}{\psi(2) \sin \pi x} \right]$$

We also have

$$\lim_{x \to 1} \left[ \frac{\log \Gamma(x+1)}{\sin \pi x} \log \Gamma(x+1) \right] = \lim_{x \to 1} \left[ \frac{\psi(x+1)}{\pi \cos \pi x} \right] \lim_{x \to 1} \log \Gamma(x+1) = 0$$

Similarly we have

$$\lim_{x \to 0} \left[ \frac{\log \Gamma(x+1)}{\sin \pi x} \log \Gamma(x+1) \right] = \lim_{x \to 0} \left[ \frac{\pi x}{\sin \pi x} \frac{\log \Gamma(x+1)}{\pi x} \log \Gamma(x+1) \right]$$

and since $\lim_{x \to 0} \left[ \frac{\log \Gamma(x+1)}{\pi x} \right] = \lim_{x \to 0} \left[ \frac{\psi(x+1)}{\pi} \right]$ we obtain $\lim_{x \to 0} [\log \Gamma(x+1) \log \sin \pi x] = 0$

Hence we have

$$\int_0^1 \log \Gamma(x+1) \cot \pi x \, dx = -\frac{1}{\pi} \int_0^1 \psi(1+x) \log \sin \pi x \, dx$$

$\square$

It was also shown in [21] that under the same conditions we have

(1.104)
$$\int_a^b p(x) \, dx = 2 \sum_{n=0}^\infty \int_a^b p(x) \cos \alpha n x \, dx$$

and we employ this with $p(x) = \log \Gamma(x+1)$ and (1.27) to give us

$$-\frac{1}{2} \int_0^1 \log \Gamma(x+1) \, dx = \sum_{n=1}^\infty \int_0^1 \log \Gamma(x+1) \cos 2n\pi x \, dx = -\frac{1}{2\pi} \sum_{n=1}^\infty \frac{si(2n\pi)}{n}$$

We therefore get

(1.105)
$$\sum_{n=1}^\infty \frac{si(2n\pi)}{n} = \frac{\pi}{2} \log(2\pi) - \pi$$

which we have previously seen in (1.54).

A number of other closed form Dirichlet series involving the sine and cosine integrals are given in my earlier paper [21].

$\square$

Motivated by a problem posed by Furdui [31] in 2009, we find a more direct derivation of (1.63) which is set out below.

We see from (1.3) that

$$Ci(nx) = \gamma + \log nx + \int_0^{nx} \frac{\cos t - 1}{t} \, dt$$

and we make the summation

$$\sum_{n=1}^\infty \frac{Ci(nx)}{n^2} = (\gamma + \log x) \sum_{n=1}^\infty \frac{1}{n^2} + \sum_{n=1}^\infty \frac{\log n}{n^2} + \sum_{n=1}^\infty \frac{1}{n^2} \left[ \int_0^{nx} \frac{\cos t - 1}{t} \, dt \right]$$

$$= \varsigma(2)(\gamma + \log x) - \varsigma'(2) + \sum_{n=1}^\infty \frac{1}{n^2} \left[ \int_0^{nx} \frac{\cos t - 1}{t} \, dt \right]$$

We see that

$$\int_0^{nx} \frac{\cos t - 1}{t}\, dt = \int_0^x \frac{\cos ny - 1}{y}\, dy$$

and thus

$$\sum_{n=1}^\infty \frac{1}{n^2}\left[\int_0^{nx} \frac{\cos t - 1}{t}\, dt\right] = \sum_{n=1}^\infty \frac{1}{n^2}\left[\int_0^x \frac{\cos ny - 1}{y}\, dy\right]$$

$$= \sum_{n=1}^\infty \frac{1}{n^2} \lim_{a\to 0}\left[\int_a^x \frac{\cos ny - 1}{y}\, dy\right]$$

$$= \lim_{a\to 0}\int_a^x \sum_{n=1}^\infty \frac{\cos ny - 1}{n^2 y}\, dy$$

We have the well known Fourier series [55, p.148]

$$\sum_{n=1}^\infty \frac{\cos ny}{n^2} = \frac{\pi^2}{6} - \frac{\pi y}{2} + \frac{y^2}{4}$$

so that

$$\sum_{n=1}^\infty \frac{\cos ny}{n^2 y} = \frac{\pi^2}{6 y} - \frac{\pi}{2} + \frac{y}{4}$$

and

$$\sum_{n=1}^\infty \frac{\cos ny - 1}{n^2 y} = -\frac{\pi}{2} + \frac{y}{4}$$

We then easily see that

(1.106)
$$\sum_{n=1}^\infty \frac{Ci(nx)}{n^2} = \varsigma(2)(\gamma + \log x) - \varsigma'(2) - \frac{\pi x}{2} + \frac{x^2}{8}$$

and with $x = 2\pi$ we obtain (1.63) again

(1.107)
$$\sum_{n=1}^\infty \frac{Ci(2n\pi x)}{n^2} = \varsigma(2)[\gamma + \log(2\pi)] - \varsigma'(2) - \frac{1}{2}\pi^2$$

More generally we have

$$(1.108) \qquad \sum_{n=1}^{\infty} \frac{Ci(nx)}{n^s} = \varsigma(s)(\gamma + \log x) - \varsigma'(s) + \lim_{a \to 0} \int_{a}^{x} \sum_{n=1}^{\infty} \frac{\cos ny - 1}{n^s y} dy$$

Similarly we may consider the sine integral function

$$Si(x) = \int_{0}^{x} \frac{\sin t}{t} dt$$

$$Si(nx) = \int_{0}^{nx} \frac{\sin t}{t} dt$$

$$\sum_{n=1}^{\infty} \frac{Si(nx)}{n^s} = \sum_{n=1}^{\infty} \frac{1}{n^s} \left[ \int_{0}^{nx} \frac{\sin t}{t} dt \right]$$

$$\sum_{n=1}^{\infty} \frac{Si(nx)}{n^s} = \int_{0}^{x} \sum_{n=1}^{\infty} \frac{\sin nt}{n^s t} dt$$

and we can certainly evaluate this for $s = 3$ using the Fourier series [55, p.148]

$$\sum_{n=1}^{\infty} \frac{\sin nt}{n^3} = \frac{1}{12}(t^3 - 3\pi t^2 + 2\pi^2 t)$$

In this regard, see also the problems proposed by Choulakian [18] in 1998.

$\square$

Using (1.10.3) we make the summation

$$\int_{0}^{1} \log x \sum_{k=0}^{\infty} \frac{\sin(2k+1)\pi x}{2k+1} dx = \sum_{k=0}^{\infty} \frac{Ci[(2k+1)\pi] - \gamma - \log[(2k+1)\pi]}{(2k+1)^2 \pi}$$

and substituting the Fourier series [51, p.149]

$$\sum_{k=0}^{\infty} \frac{\sin(2k+1)\pi x}{2k+1} = \frac{\pi}{4}$$

we obtain

$$-\frac{\pi^2}{4} = \sum_{k=0}^{\infty} \frac{Ci[(2k+1)\pi]}{(2k+1)^2} - (\gamma + \log \pi) \sum_{k=0}^{\infty} \frac{1}{(2k+1)^2} - \sum_{k=0}^{\infty} \frac{\log(2k+1)}{(2k+1)^2}$$

or equivalently

$$-\frac{\pi^2}{4} = \sum_{k=0}^{\infty} \frac{Ci[(2k+1)\pi]}{(2k+1)^2} - \frac{\pi^2}{8}(\gamma + \log \pi) - \sum_{k=0}^{\infty} \frac{\log(2k+1)}{(2k+1)^2}$$

which may be compared with (B.10)

$$\frac{1}{4}(\gamma + \log 2\pi) - \log 2 = \frac{2}{\pi^2} \sum_{n=0}^{\infty} \frac{\left[Ci[2(2n+1)\pi] - \log(2n+1)\right]}{(2n+1)^2}$$

$\square$

As a matter of interest, Abramowitz and Stegun [1, p.232] define auxiliary functions

$$f(x) = -\cos x \, si(x) + \sin x \, Ci(x) = \int_0^{\infty} \frac{\sin y}{y+x} dy$$

$$g(x) = -\cos x \, Ci(x) - \sin x \, si(x) = \int_0^{\infty} \frac{\cos y}{y+x} dy$$

and report that for $\mathrm{Re}(x) > 0$

$$f(x) = \int_0^{\infty} \frac{e^{-xu}}{1+u^2} du$$

$$g(x) = \int_0^{\infty} \frac{u e^{-xu}}{1+u^2} du$$

The above results may be derived by considering the double integral

$$I = \int_0^{\infty} \int_0^{\infty} e^{-(a+y)x} \sin y \, dx dy$$

where integrating with respect to $x$ gives us

$$\int_0^{\infty} e^{-(a+y)x} dx = \frac{1}{a+y}$$

and thus we have

$$I = \int\limits_0^\infty \frac{\sin y}{a+y}\,dy$$

Similarly, integrating with respect to $y$ gives us

$$\int\limits_0^\infty e^{-(a+y)x}\sin y\,dy = \frac{e^{-ax}}{2i}\int\limits_0^\infty e^{-yx}(e^{iy}-e^{-iy})\,dy = \frac{e^{-ax}}{1+x^2}$$

Therefore we see that

$$\int\limits_0^\infty \frac{e^{-ax}}{1+x^2}\,dx = \int\limits_0^\infty \frac{\sin y}{a+y}\,dy$$

and the validity of the operation

$$\int\limits_0^\infty dx \int\limits_0^\infty e^{-(a+y)x}\sin y\,dy = \int\limits_0^\infty dy \int\limits_0^\infty e^{-(a+y)x}\sin y\,dx$$

is confirmed by [8, p.282]. The formula

$$\int\limits_0^\infty \frac{xe^{-ax}}{1+x^2}\,dx = \int\limits_0^\infty \frac{\cos y}{a+y}\,dy$$

may be derived in a similar fashion.

Letting $t = xy$ we see that

$$\int\limits_0^\infty \frac{\sin(xy)}{x+a}\,dx = \int\limits_0^\infty \frac{\sin t}{t+ay}\,dt$$

In accordance with the above, we have from [33, p.338]

(1.109)  $$\int\limits_0^\infty \frac{ve^{-nv}}{a^2+v^2}\,dv = -[\cos(na)Ci(na) + \sin(na)si(na)]$$

and we make the summation

$$\sum_{n=1}^\infty \int\limits_0^\infty \frac{ve^{-nv}}{a^2+v^2}\,dv = -\sum_{n=1}^\infty [\cos(na)Ci(na) + \sin(na)si(na)]$$

The geometric series gives us

$$\sum_{n=1}^{\infty} \int_0^{\infty} \frac{v e^{-nv}}{a^2 + v^2} \, dv = \int_0^{\infty} \frac{v e^{-v}}{(a^2 + v^2)(1 - e^{-v})} \, dv$$

and hence we have

$$-\sum_{n=1}^{\infty} [Ci(na)\cos(na) + si(na)\sin(na)] = \int_0^{\infty} \frac{v}{(a^2 + v^2)(e^v - 1)} \, dv$$

We now let $a \to 2\pi a$ and $v \to 2\pi v$ to obtain

$$(1.110) \quad -\sum_{n=1}^{\infty} [\cos(2\pi na)Ci(2\pi na) + \sin(2\pi na)si(2\pi na)] = \int_0^{\infty} \frac{t}{(a^2 + t^2)(e^{2\pi t} - 1)} \, dt$$

Using (1.70)

$$\psi(a) = \log a - \frac{1}{2a} + 2\sum_{n=1}^{\infty} [\cos(2n\pi a)Ci(2n\pi a) + \sin(2n\pi a)si(2n\pi a)]$$

we then see that

$$(1.111) \quad \psi(a) = \log a - \frac{1}{2a} - 2\int_0^{\infty} \frac{t}{(a^2 + t^2)(e^{2\pi t} - 1)} \, dt$$

Equation (1.111) is well known and, inter alia, is reported in [56, p.251].

Integrating (1.110) gives us

$$\int_0^x da \int_0^{\infty} \frac{av}{(a^2 + v^2)(e^{2\pi v} - 1)} \, dv = -\sum_{n=1}^{\infty} \int_0^x a[\cos(2n\pi a)Ci(2n\pi a) + \sin(2n\pi a)si(2n\pi a)] \, da$$

We have in a scintilla temporis using the Wolfram Online Integrator

$$\int a[\cos(2n\pi a)Ci(2n\pi a) + \sin(2n\pi a)Si(2n\pi a)] \, da$$

$$= \frac{[2n\pi \sin(2n\pi a) + \cos(2n\pi a)]Ci(2n\pi a) - [2n\pi a\cos(2n\pi a) - \sin(2n\pi a)]Si(2n\pi a)] - \log a}{(2n\pi)^2}$$

where we note that *Mathematica* defines *SinIntegral*[$x$] as *Si(x)*. This integral may of course be easily obtained using integration by parts.

We have

$$\int_0^x a[\cos(2n\pi a)Ci(2n\pi a) + \sin(2n\pi a)Si(2n\pi a)]\,da$$

$$= \frac{\cos(2n\pi x)Ci(2n\pi x) + \sin(2n\pi x)Si(2n\pi x) - \log x - \gamma + \log(2n\pi)}{(2n\pi)^2}$$

$$+ \frac{x\sin(2n\pi x)Ci(2n\pi x) - x\cos(2n\pi x)Si(2n\pi x)}{2n\pi}$$

where we have used (1.9.1).

We have

$$\int_0^x a[\cos(2n\pi a)Ci(2n\pi a) + \sin(2n\pi a)Si(2n\pi a)]\,da$$

$$= \int_0^x a[\cos(2n\pi a)Ci(2n\pi a) + \sin(2n\pi a)si(2n\pi a)]\,da + \frac{\pi}{2}\int_0^x a\sin(2n\pi a)\,da$$

and

$$\int_0^x a\sin(2n\pi a)\,da = \frac{\sin(2n\pi x) - 2n\pi x\cos(2n\pi x)}{(2n\pi)^2}$$

$$\int_0^x a[\cos(2n\pi a)Ci(2n\pi a) + \sin(2n\pi a)si(2n\pi a)]\,da$$

$$= \frac{\cos(2n\pi x)Ci(2n\pi x) + \sin(2n\pi x)Si(2n\pi x) - \log x - \gamma + \log(2n\pi)}{(2n\pi)^2}$$

$$+ \frac{x\sin(2n\pi x)Ci(2n\pi x) - x\cos(2n\pi x)Si(2n\pi x)}{2n\pi}$$

$$- \frac{\pi}{2}\frac{\sin(2n\pi x)}{(2n\pi)^2} + \frac{\pi}{2}\frac{x\cos(2n\pi x)}{2n\pi}$$

$$= \frac{\cos(2n\pi x)Ci(2n\pi x) + \sin(2n\pi x)si(2n\pi x) - \log x - \gamma + \log(2n\pi)}{(2n\pi)^2}$$

$$+ \frac{x\sin(2n\pi x)Ci(2n\pi x) - x\cos(2n\pi x)si(2n\pi x)}{2n\pi}$$

Therefore we have

$$\sum_{n=1}^{\infty} \int_0^x a[\cos(2n\pi a)Ci(2n\pi a) + \sin(2n\pi a)si(2n\pi a)]\,da$$

$$= \sum_{n=1}^{\infty} \frac{\cos(2n\pi x)Ci(2n\pi x) + \sin(2n\pi x)si(2n\pi x)}{(2n\pi)^2} - \frac{1}{24}(\log x - \log(2\pi) + \gamma) - \frac{\varsigma'(2)}{4\pi^2}$$

$$+ x\sum_{n=1}^{\infty} \frac{\sin(2n\pi x)Ci(2n\pi x) - \cos(2n\pi x)si(2n\pi x)}{2n\pi}$$

Using (1.41) this becomes

$$= \sum_{n=1}^{\infty} \frac{\cos(2n\pi x)Ci(2n\pi x) + \sin(2n\pi x)si(2n\pi x)}{(2n\pi)^2} - \frac{1}{24}(\log x + \gamma) - \frac{\varsigma'(2)}{4\pi^2}$$

$$+ \frac{1}{2}x\log\Gamma(x) - \frac{x}{4}\log(2\pi) - \frac{x}{2}\left(x - \frac{1}{2}\right)\log x + \frac{1}{2}x^2$$

Integrating (1.110) gives us

$$I = \int_0^x da \int_0^{\infty} \frac{av}{(a^2 + v^2)(e^{2\pi v} - 1)}\,dv = \frac{1}{2}\int_0^{\infty} \frac{v\,dv}{e^{2\pi v} - 1} \int_0^x \frac{2a}{a^2 + v^2}\,da$$

$$= \frac{1}{2}\int_0^{\infty} \frac{v\log(x^2 + v^2)}{e^{2\pi v} - 1}\,dv - \int_0^{\infty} \frac{v\log v}{e^{2\pi v} - 1}\,dv$$

We have [24]

(1.112) $$\log G(1+x) = \frac{1}{2}x^2\left[\log x - \frac{3}{2}\right] + \frac{1}{2}x\log(2\pi) + \varsigma'(-1) - \int_0^{\infty} \frac{v\log\left(v^2 + x^2\right)}{e^{2\pi v} - 1}\,dv$$

where $G(x)$ is the Barnes double gamma function. This result was originally obtained by Adamchik [3] in 2004. With $x = 0$ we have

(1.113) $$\int_0^{\infty} \frac{v\log v}{e^{2\pi v} - 1}\,dv = \frac{1}{2}\varsigma'(-1)$$

Therefore we obtain

$$I = \frac{1}{4}x^2\left[\log x - \frac{3}{2}\right] + \frac{1}{4}x\log(2\pi) - \frac{1}{2}\log G(1+x)$$

Hence we obtain

(1.114) $\displaystyle\sum_{n=1}^{\infty} \frac{\cos(2n\pi x)Ci(2n\pi x) + \sin(2n\pi x)si(2n\pi x)}{(2n\pi)^2}$

$$= \frac{1}{24}(\log x - \log(2\pi) + \gamma) + \frac{1}{4}x^2\left[\log x - \frac{3}{2}\right] + \frac{1}{2}x\log(2\pi) - \frac{1}{2}\log G(1+x) - \frac{\varsigma'(2)}{4\pi^2}$$

$$- \frac{1}{2}x\log\Gamma(x) + \frac{x}{2}\left(x - \frac{1}{2}\right)\log x - \frac{1}{2}x^2$$

Alternatively, we multiply (1.111) by $a$ and integrate to obtain

$$\int_0^x a\psi(a)da = \frac{1}{2}x^2\log x - \frac{1}{4}x^2 - \frac{1}{2}x - \frac{1}{2}x^2\left[\log x - \frac{3}{2}\right] - \frac{1}{2}x\log(2\pi) + \log G(1+x)$$

Integration by parts results in

$$\int_0^x a\psi(a)da = x\log\Gamma(x) - \int_0^x \log\Gamma(a)da$$

and we obtain Alexeiewsky's theorem [52, p.32]

(1.115) $\displaystyle\int_0^x \log\Gamma(a)da = x\log\Gamma(x) - \log G(1+x) - \frac{1}{2}x^2 + \frac{1}{2}x + \frac{1}{2}x\log(2\pi)$

$\square$

We also have from [33, p.338]

$$\int_0^{\infty} \frac{e^{-\mu v}}{n^2 + v^2}dv = \frac{1}{n}[\sin(n\mu)Ci(n\mu) - \cos(n\mu)si(n\mu)]$$

and we make the summation

$$\sum_{n=1}^{\infty}\int_0^{\infty} \frac{e^{-\mu v}}{n^2 + v^2}dv = \sum_{n=1}^{\infty}\frac{1}{n}[\sin(n\mu)Ci(n\mu) - \cos(n\mu)si(n\mu)]$$

Assuming that interchanging the order of summation and integration is valid

$$\sum_{n=1}^{\infty} \int_0^{\infty} \frac{e^{-\mu v}}{n^2+v^2}\, dv = \int_0^{\infty} e^{-\mu v} \sum_{n=1}^{\infty} \frac{1}{n^2+v^2}\, dv$$

and using (obtained by letting $\alpha = iv$ in (4.7))

(1.116) $$\pi \coth \pi v = \frac{1}{v} + 2v \sum_{n=1}^{\infty} \frac{1}{n^2+v^2}$$

we have

(1.117) $$\frac{1}{2} \int_0^{\infty} \frac{\pi v \coth \pi v - 1}{v^2} e^{-\mu v}\, dv = \sum_{n=1}^{\infty} \frac{1}{n}[\sin(n\mu)Ci(n\mu) - \cos(n\mu)si(n\mu)]$$

*Mathematica* cannot evaluate this integral.

With $\mu = \pi$ we have

$$\frac{1}{2} \int_0^{\infty} \frac{\pi v \coth \pi v - 1}{v^2} e^{-\pi v}\, dv = \sum_{n=1}^{\infty} \frac{(-1)^n}{n} si(n\pi)$$

and referring to (1.49) we obtain

(1.118) $$\frac{1}{2} \int_0^{\infty} \frac{\pi v \coth \pi v - 1}{v^2} e^{-\pi v}\, dv = \frac{\pi}{2} \log 2 - \frac{\pi}{2}$$

With $\mu = 2\pi$ we have

$$\frac{1}{2} \int_0^{\infty} \frac{\pi v \coth \pi v - 1}{v^2} e^{-2\pi v}\, dv = \sum_{n=1}^{\infty} \frac{si(2n\pi)}{n}$$

and referring to (1.54) we obtain

(1.119) $$\frac{1}{2} \int_0^{\infty} \frac{\pi v \coth \pi v - 1}{v^2} e^{-2\pi v}\, dv = \frac{\pi}{2} \log(2\pi) - \pi$$

$\square$

We note [33, p.336, 3.341] that

$$\int_0^{\pi/2} \exp(-p\tan x)dx = \sin p\, Ci(p) - \cos p\, si(p)$$

so that

$$\int_0^{\pi/2} \exp(-2n\pi a\tan x)dx = \sin 2n\pi a\, Ci(2n\pi a) - \cos 2n\pi a\, si(2n\pi a)$$

Referring to (1.41)

$$\log\Gamma(a) = \frac{1}{2}\log(2\pi) + \left(a - \frac{1}{2}\right)\log a - a + \frac{1}{\pi}\sum_{n=1}^{\infty}\frac{1}{n}[\sin(2n\pi a)Ci(2n\pi a) - \cos(2n\pi a)si(2n\pi a)]$$

and using

$$\sum_{n=1}^{\infty}\frac{1}{n}\exp(-2n\pi a\tan x) = -\log[1 - \exp(-2\pi a\tan x)]$$

we obtain the integral representation

$$(1.119.1)\qquad \log\Gamma(a) = \frac{1}{2}\log(2\pi) + \left(a - \frac{1}{2}\right)\log a - a - \frac{1}{\pi}\int_0^{\pi/2}\log[1 - \exp(-2\pi a\tan x)]dx$$

and we also note that

$$\int_0^{\pi/2}\frac{2\pi a\sec^2 x}{\exp(-2\pi a\tan x) - 1}dx = \int_0^{\pi/2}\log[1 - \exp(-2\pi a\tan x)]dx$$

$\square$

As noted by Nielsen [44, p.72] and Bartle [8, p.345] we have the Fourier series (provided $x$ is not an integer)

$$(1.120)\qquad \cos tx = \frac{2x\sin\pi x}{\pi}\left[\frac{1}{2x^2} + \sum_{n=1}^{\infty}\frac{(-1)^{n+1}}{n^2 - x^2}\cos nt\right]$$

$$= \frac{\sin\pi x}{\pi}\left[\frac{1}{x} + 2x\sum_{n=1}^{\infty}\frac{(-1)^{n+1}}{n^2 - x^2}\cos nt\right]$$

and with $t = 2\pi$ we obtain

$$\frac{\pi}{\sin \pi x} = \frac{1}{x} + 2x \sum_{n=1}^{\infty} \frac{(-1)^{n+1}}{n^2 - x^2}$$

We then have

$$\cos tx = \frac{\sin \pi x}{\pi} \left[ \frac{\pi}{\sin \pi x} - 2x \sum_{n=1}^{\infty} \frac{(-1)^{n+1}}{n^2 - x^2} + 2x \sum_{n=1}^{\infty} \frac{(-1)^{n+1}}{n^2 - x^2} \cos nt \right]$$

which we write as

$$\cos tx - 1 = \frac{2x \sin \pi x}{\pi} \sum_{n=1}^{\infty} \frac{(-1)^{n+1}}{n^2 - x^2} (\cos nt - 1)$$

We now divide by $t$ and integrate to obtain

$$\int_0^u \frac{\cos tx - 1}{t} dt = \frac{2x \sin \pi x}{\pi} \sum_{n=1}^{\infty} \frac{(-1)^{n+1}}{n^2 - x^2} \int_0^u \frac{\cos nx - 1}{t} dt$$

It is easily seen that

$$\int_0^u \frac{\cos tx - 1}{t} dt = \int_0^{ux} \frac{\cos y - 1}{y} dy$$

and hence we have

$$\int_0^u \frac{\cos tx - 1}{t} dt = Ci(ux) - \gamma - \log(ux)$$

Therefore we obtain

$$Ci(ux) - \gamma - \log(ux) = \frac{2x \sin \pi x}{\pi} \sum_{n=1}^{\infty} \frac{(-1)^{n+1}}{n^2 - x^2} [Ci(nu) - \gamma - \log(nu)]$$

$$= \frac{2x \sin \pi x}{\pi} \sum_{n=1}^{\infty} \frac{(-1)^{n+1}}{n^2 - x^2} [Ci(nu) - \log n] - (\gamma + \log u) \frac{2x \sin \pi x}{\pi} \sum_{n=1}^{\infty} \frac{(-1)^{n+1}}{n^2 - x^2}$$

$$= \frac{2x \sin \pi x}{\pi} \sum_{n=1}^{\infty} \frac{(-1)^{n+1}}{n^2 - x^2} [Ci(nu) - \log n] - (\gamma + \log u) \left[ 1 - \frac{\sin \pi x}{\pi x} \right]$$

where we have used

$$1 - \frac{\sin \pi x}{\pi x} = \frac{2x \sin \pi x}{\pi} \sum_{n=1}^{\infty} \frac{(-1)^{n+1}}{n^2 - x^2}$$

Hence we obtain

(1.121) $\quad Ci(ux) - \log x = \frac{2x \sin \pi x}{\pi} \sum_{n=1}^{\infty} \frac{(-1)^{n+1}}{n^2 - x^2} [Ci(nu) - \log n] + \frac{(\gamma + \log u) \sin \pi x}{\pi x}$

which corrects a misprint in Nielsen's book [44, p.72].

In examining the limit as $x \to 0$, using (1.19) we see that both sides of the equation are equal.

Using L'Hôpital's rule we easily see that both sides are also equivalent in the case where $x = k$ where $k$ is an integer.

In passing, we note that we may write (1.120) as

(1.122) $\quad \frac{\pi x \cos tx - \sin \pi x}{2x^2 \sin \pi x} = \sum_{n=1}^{\infty} \frac{(-1)^{n+1}}{n^2 - x^2} \cos nt$

and then take the limit as $x \to 0$. Applying L'Hôpital's rule three times gives us the well known Fourier series [55, p.148]

$$\sum_{n=1}^{\infty} \frac{(-1)^{n+1} \cos nt}{n^2} = \frac{1}{12}(\pi^2 - 3t^2)$$

With $t = \pi$ in (1.122) we obtain

$$\frac{\pi x \cos \pi x - \sin \pi x}{2x^2 \sin \pi x} = -\sum_{n=1}^{\infty} \frac{1}{n^2 - x^2}$$

which is the same as (4.7).

Differentiating (1.122) with respect to $t$ results in

$$\frac{\pi \sin tx}{2 \sin \pi x} = \sum_{n=1}^{\infty} \frac{(-1)^{n+1} n}{n^2 - x^2} \sin nt$$

and applying L'Hôpital's rule as $x \to 0$ gives us

$$\lim_{x \to 0} \frac{\pi \sin tx}{2 \sin \pi x} = \lim_{x \to 0} \frac{\pi t \cos tx}{2\pi \cos \pi x} = \frac{t}{2}$$

and we end up with the well-known Fourier series [55, p.148]

$$\frac{t}{2} = \sum_{n=1}^{\infty} \frac{(-1)^{n+1}}{n} \sin nt$$

This could also be employed to derive (1.49).

Integrating (1.122) with respect to $t$ results in

$$\frac{\pi \sin tx - t \sin \pi x}{2x^2 \sin \pi x} = \sum_{n=1}^{\infty} \frac{1}{n} \frac{(-1)^{n+1}}{n^2 - x^2} \sin nt$$

$\square$

We note that [37a] gives a closed form for the following Fourier series in response to a problem posed by Fettis and Glasser for $|x| < \pi$, $(p,q) = 1$, $q \geq 2$

$$\sum_{n=1}^{\infty} \frac{(-1)^n}{n^2 - (p/q)^2} \sin nt = \frac{1}{pq} \sum_{k=1}^{q-1} \sin\left[\frac{p}{q}(x + k\pi)\right] \frac{\sin(k\pi p/q)}{\sin(\pi p/q)} \cdot \log\left|\frac{\sin(1/2q)(x + 2k + 1)\pi)}{\sin(1/2q)(x + 2k - 1)\pi)}\right|$$

$\square$

We have [52, p.14]

$$\psi(x) = \log x + \sum_{k=0}^{\infty} \left[\log\left(1 + \frac{1}{x + k}\right) - \frac{1}{x + k}\right]$$

which could then be used to evaluate $\int_0^1 \psi(x) \sin p\pi x \, dx$ using (1.8), (1.10.3) and the integral

$$\int \frac{\sin p\pi x}{x + k} dx = \cos(kp) Si[p(x + k)] - \sin(kp) Ci[p(x + k)]$$

## 2. Another approach to the $\log \Gamma(x)$ integrals

Part of the following analysis is taken from the 1862 book "Exposé de la théorie, propriétés, des formules de transformation, et méthodes d'évaluation des intégrales définies" by David Bierens de Haan [35, p.269], a copy of which is available on the internet courtesy of the University of Michigan Historical Mathematics Collection. The voluminous integrals evaluated by Bierens de Haan were an important source for the table of integrals subsequently compiled by Gradshteyn and Ryzhik [33]. Bierens de Haan (1822-1895) was, inter alia, a mathematician, an actuary and a historian; a short interesting biography of him appears in Talvila's paper [53] and a more detailed presentation is given in [49].

Let us designate the following integrals as

$$(2.1) \qquad I(p) = \int_0^1 \log \Gamma(x) \cos p\pi x\, dx$$

$$(2.2) \qquad K(p) = \int_0^1 \log \Gamma(x) \sin p\pi x\, dx$$

Then we have by letting $x \to 1 - x$

$$(2.3) \qquad I(p) = \cos p\pi \int_0^1 \log \Gamma(1-x) \cos p\pi x\, dx + \sin p\pi \int_0^1 \log \Gamma(1-x) \sin p\pi x\, dx$$

and we see that

$$(2.4) \qquad I(1/2) = \int_0^1 \log \Gamma(1-x) \sin(\pi x/2)\, dx$$

We have from the definition

$$I(1-p) = \int_0^1 \log \Gamma(x) [\cos \pi x \cos p\pi x + \sin \pi x \sin p\pi x]\, dx$$

The integral

$$J(p) = \int_0^1 \log \Gamma(1-x) \cos p\pi x\, dx$$

becomes by letting $x \to 1 - x$

$$= \int_0^1 \log \Gamma(x) [\cos p\pi \cos p\pi x + \sin p\pi \sin p\pi x]\, dx$$

Similarly, with $L(p)$ defined by

$$L(p) = \int_0^1 \log \Gamma(1-x) \sin p\pi x\, dx$$

becomes by letting $x \to 1-x$

$$= \int\limits_0^1 \log \Gamma(x)[\sin p\pi \cos p\pi x - \cos p\pi \sin p\pi x]\,dx$$

In the same way as above we find that

(2.5)    $K(p) = \sin p\pi \int\limits_0^1 \log \Gamma(1-x) \cos p\pi x\,dx - \cos p\pi \int\limits_0^1 \log \Gamma(1-x) \sin p\pi x\,dx$

and we see that

(2.5.1)    $K(1/2) = \int\limits_0^1 \log \Gamma(1-x) \cos(\pi x/2)\,dx$

We have Euler's reflection formula [10] for all $z \in \mathbf{C}$ except $z = 0, \pm n$ where $n \in \mathbf{N}$

$$\Gamma(z)\Gamma(1-z) = \frac{\pi}{\sin \pi z}$$

and therefore for $0 < x < 1$

(2.6)    $\log \Gamma(x) + \log \Gamma(1-x) = \log(2\pi) - \log(2\sin \pi x)$

Multiplying this equation by $\cos p\pi x$ and integrating, we obtain

$$\int\limits_0^1 \log \Gamma(x) \cos p\pi x\,dx + \int\limits_0^1 \log \Gamma(1-x) \cos p\pi x\,dx$$

$$= \log(2\pi) \int\limits_0^1 \cos p\pi x\,dx - \int\limits_0^1 \log(2\sin \pi x) \cos p\pi x\,dx$$

Then, using the Fourier series [16, p.241]

(2.7)    $\log\left[2\sin(t/2)\right] = -\sum\limits_{n=1}^{\infty} \frac{\cos nt}{n}$      $0 < t < 2\pi$

and, assuming that changing the order of integration and summation is valid, we get

(2.8)    $\int\limits_0^1 \log \Gamma(x) \cos p\pi x\,dx + \int\limits_0^1 \log \Gamma(1-x) \cos p\pi x\,dx$

$$= \frac{\sin p\pi}{p\pi} \log(2\pi) + \sum_{n=1}^{\infty} \frac{1}{2n} \int_0^1 \Big[ \cos\{(2n+p)\pi x\} + \cos\{(2n-p)\pi x\} \Big] dx$$

$$= \frac{\sin p\pi}{p\pi} \log(2\pi) + \sum_{n=1}^{\infty} \frac{1}{2n} \left[ \frac{\sin\{(2n+p)\pi\}}{(2n+p)\pi} + \frac{\sin\{(2n-p)\pi\}}{(2n-p)\pi} \right]$$

Hence, provided $p \neq 2n$, as originally determined by Bierens de Haan [35, p.270] we have (where I have corrected a typographical error)

(2.9)
$$\int_0^1 \log\Gamma(x) \cos p\pi x \, dx + \int_0^1 \log\Gamma(1-x) \cos p\pi x \, dx = \frac{\sin p\pi}{p\pi} \log(2\pi) - \frac{p \sin p\pi}{\pi} \sum_{n=1}^{\infty} \frac{1}{n} \frac{1}{4n^2 - p^2}$$

Rather confusingly, Bierens de Haan [35, p.270] reports the left-hand side of (2.9) as being equal to $I(p) + I(1-p)$.

From (2.9) we see that with $p = 0$

(2.10)
$$\int_0^1 \log\Gamma(x) \, dx + \int_0^1 \log\Gamma(1-x) \, dx = \log(2\pi)$$

and by a simple change of variable we have

$$\int_0^1 \log\Gamma(1-x) \, dx = \int_0^1 \log\Gamma(x) \, dx$$

Hence we obtain Raabe's integral [56, p.261]

(2.11)
$$\int_0^1 \log\Gamma(x) \, dx = \frac{1}{2} \log(2\pi)$$

From (2.9) we see that with $p = 1$

(2.12)
$$\int_0^1 \log\Gamma(x) \cos\pi x \, dx + \int_0^1 \log\Gamma(1-x) \cos\pi x \, dx = 0$$

which may of course also be obtained by a simple change of variable.

With $p = 1$ we have

$$\int_0^1 \log \Gamma(x) \cos \pi x \, dx + \int_0^1 \log \Gamma(1-x) \cos \pi x \, dx$$

$$= \log(2\pi) \int_0^1 \cos \pi x \, dx - \int_0^1 \log(2 \sin \pi x) \cos \pi x \, dx$$

$$= -\int_0^1 \log(2 \sin \pi x) \cos \pi x \, dx$$

and, with the substitution $y = \sin \pi x$, we see that the integral on the right-hand side vanishes. This is easily confirmed because, with the substitution $u = 1 - x$, we see that

$$\int_0^1 \log \Gamma(1-x) \cos \pi x \, dx = -\int_0^1 \log \Gamma(u) \cos \pi u \, du$$

We are therefore unable to evaluate $\int_0^1 \log \Gamma(x) \cos \pi x \, dx$ using this methodology.

We note in passing that Raabe's integral (2.11) may also be easily obtained as follows by integrating (2.6), i.e. with $p = 0$

$$\int_0^1 \log \Gamma(x) \, dx + \int_0^1 \log \Gamma(1-x) \, dx = \log(2\pi) - \int_0^1 \log(2 \sin \pi x) \, dx$$

where we also use the well-known integral due to Euler

$$\int_0^1 \log \sin \pi x \, dx = -\log 2$$

to determine that

$$\int_0^1 \log \Gamma(x) \, dx + \int_0^1 \log \Gamma(1-x) \, dx = \log(2\pi)$$

and finally we see that

$$\int_0^1 \log \Gamma(x) \, dx = \int_0^1 \log \Gamma(1-x) \, dx$$

□

In passing, we note that an interesting derivation of Euler's integral appears in an 1864 paper by Jeffery [38, p.94]. We consider

$$\int_0^{\pi/2} \log \sin x \, dx = \frac{1}{2} \int_0^{\pi/2} \log \sin^2 x \, dx$$

and with the substitution $t = \sin^2 x$ we have

$$= \frac{1}{4} \int_0^1 \frac{\log t}{t^{\frac{1}{2}}(1-t)^{\frac{1}{2}}} \, dt$$

$$= \frac{1}{4} \frac{\partial}{\partial u} B(u, v) \Big|_{u=v=1/2}$$

$$= \frac{1}{4} \left\{ \Gamma'\left(\frac{1}{2}\right) \Gamma\left(\frac{1}{2}\right) - \Gamma^2\left(\frac{1}{2}\right) \Gamma'(1) \right\}$$

$$= -\frac{\pi}{2} \log 2$$

$\square$

Letting $p = 1/2$ in (2.9) we see that

$$\int_0^1 \log \Gamma(x) \cos(\pi x/2) \, dx + \int_0^1 \log \Gamma(1-x) \cos(\pi x/2) \, dx = \frac{2}{\pi} \log(2\pi) - \frac{1}{2\pi} \sum_{n=1}^\infty \frac{1}{n} \frac{1}{4n^2 - (1/4)}$$

$$= \frac{2}{\pi} \log(2\pi) - \frac{2}{\pi} \sum_{n=1}^\infty \frac{1}{n} \frac{1}{16n^2 - 1}$$

With the obvious substitution we see that

$$\int_0^1 \log \Gamma(1-x) \cos(\pi x/2) \, dx = \int_0^1 \log \Gamma(x) \sin(\pi x/2) \, dx$$

Hence we have

(2.13) $$\int_0^1 \log \Gamma(x) [\cos(\pi x/2) + \sin(\pi x/2)] \, dx = \frac{2}{\pi} \log(2\pi) - \frac{2}{\pi} \sum_{n=1}^\infty \frac{1}{n} \frac{1}{16n^2 - 1}$$

With regard to the above infinite series, Ramanujan [11, p.29] showed that

(2.14) $$\sum_{n=1}^{\infty} \frac{1}{n} \frac{1}{16n^2 - 1} = 3\log 2 - 2$$

In fact, as reported in [11, p.29], this summation was the **first** problem submitted by Ramanujan in 1911 to the Journal of the Indian Mathematical Society.

Ramanujan [11, p.26] also showed that

(2.15) $$\sum_{n=1}^{\infty} \frac{1}{n} \frac{1}{4n^2 - 1} = 2\log 2 - 1$$

This series is contained in [33, p.9,] where it is attributed to Bromwich's text [15, p.51]. A proof is also given in Ramanujan's Notebooks [11, p.26]. In [15] the proof is based on the definition of Euler's constant $\gamma$ and it is therefore quite apt that the series should be connected with the gamma function.

Therefore we obtain from (2.13) and (2.14)

(2.16) $$\int_0^1 \log \Gamma(x)[\cos(\pi x / 2) + \sin(\pi x / 2)]\,dx = \frac{2}{\pi}\big[\log \pi - 2\log 2 + 2\big]$$

or equivalently

(2.16.1) $$\int_0^1 \log \Gamma(x) \cos\left(\frac{\pi}{4}[2x-1]\right)dx = \frac{\sqrt{2}}{\pi}\big[\log \pi - 2\log 2 + 2\big]$$

$\square$

By definition we have

$$I\left(\frac{1}{2}\right) = \int_0^1 \log \Gamma(x)\cos(\pi x / 2)\,dx$$

$$= 2\int_0^{1/2} \log \Gamma(2t)\cos(\pi t)\,dt$$

We have Legendre's duplication formula for the gamma function (see [56, p.240] and [25] for example)

(2.16.2) $$2^{2x-1}\Gamma(x)\Gamma\left(x + \frac{1}{2}\right) = \sqrt{\pi}\,\Gamma(2x)$$

which may be expressed as

$$\log \Gamma(2x) = \log \Gamma(x) + \log \Gamma\left(x + \frac{1}{2}\right) + (2x-1)\log 2 - \frac{1}{2}\log \pi$$

Hence we have

$$\frac{1}{2}I(1/2) = \int_0^{1/2} \log \Gamma(t)\cos(\pi t)\,dt + \int_0^{1/2} \log \Gamma\left(t + \frac{1}{2}\right)\cos(\pi t)\,dt$$

$$+ \int_0^{1/2}\left[(2t-1)\log 2 - \frac{1}{2}\log \pi\right]\cos(\pi t)\,dt$$

Using the obvious change of variables we have

$$\int_0^{1/2} \log \Gamma\left(t + \frac{1}{2}\right)\cos(\pi t)\,dt = \int_{1/2}^1 \log \Gamma(u)\cos \pi\left(\frac{1}{2} - u\right)du$$

$$= \int_{1/2}^1 \log \Gamma(u)\sin \pi u\,du$$

$$= \int_0^1 \log \Gamma(u)\sin \pi u\,du - \int_0^{1/2} \log \Gamma(u)\sin \pi u\,du$$

$$= \frac{1}{\pi}\left[\log \frac{\pi}{2} + 1\right] - \int_0^{1/2} \log \Gamma(u)\sin \pi u\,du$$

where we used (2.26) below. Therefore we have

$$(2.17)\quad \int_0^{1/2} \log \Gamma(t)[\cos \pi t - \sin \pi t]\,dt + \left[\frac{1}{\pi} - \frac{2}{\pi^2}\right]\log 2 - \frac{1}{\pi}\left[\log 2 + \frac{1}{2}\log \pi\right] + \frac{1}{\pi}\left[\log \frac{\pi}{2} + 1\right]$$

$$\square$$

We define

$$\phi(a,k) = 1 + 2\sum_{n=1}^{k} \frac{1}{an} \frac{1}{(an)^2 - 1}$$

and

(2.18)
$$\phi(a) = \lim_{k \to \infty} \phi(a,k) = 1 + 2\sum_{n=1}^{\infty} \frac{1}{an} \frac{1}{(an)^2 - 1}$$

Sitaramachandrarao [50] has shown that if $p, q$ are positive integers and if $\phi(p,q)$ exists, then

(2.19)
$$\phi\left(\frac{p}{q}\right) = \frac{2q}{p} \left[ \log 2p - 2\sum \cos\left(\frac{2\pi qj}{p}\right) \log \sin\left(\frac{\pi j}{p}\right) \right] + 2q \sum_{j=1}^{\lfloor q/p \rfloor} \frac{1}{pj - q}$$

where the (corrected) first summation runs over $0 < j < p/2$, and $\lfloor x \rfloor$ denotes the largest integer $\leq x$.

Hence we have from (2.9)

$$I\left(\frac{p}{q}\right) + J\left(1 - \frac{p}{q}\right) = \frac{q \sin(p\pi/q)}{p\pi} \log(2\pi) - \frac{p \sin(p\pi/q)}{q\pi} \sum_{n=1}^{\infty} \frac{1}{n} \frac{1}{4n^2 - (p/q)^2}$$

$$= \frac{q \sin(p\pi/q)}{p\pi} \log(2\pi) - \frac{q^2 \sin(p\pi/q)}{p^2\pi} \sum_{n=1}^{\infty} \frac{2}{2qn/p} \frac{1}{(2qn/p)^2 - 1}$$

and thus

(2.20)
$$I\left(\frac{p}{q}\right) + J\left(1 - \frac{p}{q}\right) = \frac{q \sin(p\pi/q)}{p\pi} \log(2\pi) - \frac{q^2 \sin(p\pi/q)}{p^2\pi} \left[ \phi\left(\frac{2q}{p}\right) - 1 \right]$$

As mentioned by Sitaramachandrarao, (2.19) is clearly related to Gauss's formula for the digamma function (see [15, p.522] and [52, p.18].

$$\psi\left(\frac{p}{q}\right) = -\gamma - \frac{\pi}{2} \cot\left(\frac{p\pi}{q}\right) - \log q + \sum_{j=1}^{q-1} \cos\left(\frac{2\pi pj}{q}\right) \log\left[ 2\sin\left(\frac{\pi j}{q}\right) \right]$$

With $p = 2$ and $q = 2k + 1$ in (2.19) we obtain

$$\phi\left(\frac{2}{2k+1}\right) = 2(2k+1)\log 2 + 2(2k+1) \sum_{j=1}^{k} \frac{1}{2j - (2k+1)}$$

or equivalently

$$\phi\left(\frac{2}{2k+1}\right) = 2(2k+1)\log 2 - 2(2k+1)\sum_{j=1}^{k}\frac{1}{2j-1}$$

By definition we have

$$\phi\left(\frac{2}{2k+1}\right) = 1 + (2k+1)^3 \sum_{n=1}^{\infty}\frac{1}{n}\frac{1}{4n^2 - (2k+1)^2}$$

and hence we have

(2.21) $$\sum_{n=1}^{\infty}\frac{1}{n}\frac{1}{4n^2-(2k+1)^2} = \frac{1}{(2k+1)^2}\left[2\log 2 - \frac{1}{2k+1} - 2\sum_{j=0}^{k-1}\frac{1}{2j+1}\right]$$

This identity is employed in (5.6).

For example, with $k=0$ in (2.21) we recover Ramanujan's formula (2.15).

We immediately notice the strong similarity with the following integral which is evaluated in Section 6

$$\int_0^1 \log \sin \pi x \sin(2k+1)\pi x\, dx = \frac{2}{(2k+1)\pi}\left[\log 2 - \frac{1}{2k+1} - 2\sum_{j=0}^{k-1}\frac{1}{2j+1}\right]$$

and, with a minor adjustment, we see that

$$\int_0^1 \log(2\sin \pi x)\sin(2k+1)\pi x\, dx = \frac{2}{(2k+1)\pi}\left[2\log 2 - \frac{1}{2k+1} - 2\sum_{j=0}^{k-1}\frac{1}{2j+1}\right]$$

The very fact that Ramanujan deals with $\phi(a)$ in Chapter 8 of [11], which also involves connections with Kummer's formula (3.11) and related matters, strongly suggests that he was well aware of equation (2.9).

We note from equation (A.3) in the Appendix that

(2.22) $$\psi(1+x) + \psi(1-x) + 2\gamma = -2x^2\sum_{n=1}^{\infty}\frac{1}{n}\frac{1}{n^2-x^2}$$

which concurs with Prudnikov et al. [48, 5.1.25-13].

In another paper [27] we have derived rapidly converging series for Catalan's constant $G = \beta(2)$ and for Apéry's constant $\varsigma(3)$ by employing (2.22). These are set out below:

$$G = 3(1 - \log 2) - \sum_{n=1}^{\infty} \frac{1}{n} \frac{1}{(16n^2 - 1)^2}$$

$$7\zeta(3) = 12(5 - 4\log 2) - 20G + 16\sum_{n=1}^{\infty} \frac{1}{n} \frac{1}{(16n^2 - 1)^3}$$

It may be noted that *Mathematica* confirms these identities both numerically and analytically.

The method may be easily generalised to produce new series representations for values of the Riemann zeta function $\zeta(2n+1)$ and the Dirichlet beta function $\beta(2n)$.

With $x \to 1/x$ in (2.22) we get

$$\psi\left(\frac{1}{x}\right) + \psi\left(1 - \frac{1}{x}\right) = -2\gamma - x\left[1 + 2\sum_{n=1}^{\infty} \frac{1}{nx} \frac{1}{(nx)^2 - 1}\right]$$

as reported by Ramanujan [11, p.186].

With $x = p/2$ in (2.22) we have

(2.23) $$\psi\left(1 + \frac{p}{2}\right) + \psi\left(1 - \frac{p}{2}\right) + 2\gamma = -2p^2 \sum_{n=1}^{\infty} \frac{1}{n} \frac{1}{4n^2 - p^2}$$

with the result that we may write (2.9) as

(2.24) $$\int_0^1 \log\Gamma(x) \cos p\pi x \, dx + \int_0^1 \log\Gamma(1-x) \cos p\pi x \, dx$$

$$= \frac{\sin p\pi}{p\pi}\left[\gamma + \log(2\pi) + \frac{1}{2}\psi\left(1 + \frac{p}{2}\right) + \frac{1}{2}\psi\left(1 - \frac{p}{2}\right)\right]$$

For example, with $p = 1/2$ in (2.23) we have

$$= \frac{1}{\pi}\left[\gamma + \log(2\pi) + \frac{1}{2}\psi\left(\frac{5}{4}\right) + \frac{1}{2}\psi\left(\frac{3}{4}\right)\right]$$

and using [52, p.20]

$$\psi\left(\frac{3}{4}\right) = -\gamma + \frac{1}{2}\pi - 3\log 2 \qquad \psi\left(\frac{5}{4}\right) = -\gamma - \frac{1}{2}\pi - 3\log 2 + 4$$

we therefore have another derivation of (2.16).

□

We now consider integrals of the form $\int_0^1 \log \Gamma(x) \sin p\pi x \, dx$. Multiplying (2.6) by $\sin p\pi x$ and integrating, we obtain

$$\int_0^1 \log \Gamma(x) \sin p\pi x \, dx + \int_0^1 \log \Gamma(1-x) \sin p\pi x \, dx$$

$$= \frac{1-\cos p\pi}{p\pi} \log(2\pi) + \sum_{n=1}^{\infty} \frac{1}{2n} \int_0^1 \Big[ \sin\{(2n+p)\pi x\} - \sin\{(2n-p)\pi x\} \Big] dx$$

$$= \frac{1-\cos p\pi}{p\pi} \log(2\pi) + \sum_{n=1}^{\infty} \frac{1}{2n} \left[ \frac{1-\cos\{(2n+p)\pi\}}{(2n+p)\pi} - \frac{1-\cos\{(2n-p)\pi\}}{(2n-p)\pi} \right]$$

Hence, provided $p \neq 2n$, we have

(2.25) $\int_0^1 \log \Gamma(x) \sin p\pi x \, dx + \int_0^1 \log \Gamma(1-x) \sin p\pi x \, dx$

$$= \frac{1-\cos p\pi}{p\pi} \log(2\pi) - \frac{p(1-\cos p\pi)}{\pi} \sum_{n=1}^{\infty} \frac{1}{n} \frac{1}{4n^2 - p^2}$$

With $p = 1$ in (2.25) we obtain

$$\int_0^1 \log \Gamma(x) \sin \pi x \, dx + \int_0^1 \log \Gamma(1-x) \sin \pi x \, dx = \frac{2}{\pi} \log(2\pi) - \frac{2}{\pi} \sum_{n=1}^{\infty} \frac{1}{n} \frac{1}{4n^2 - 1}$$

and substituting (2.15) we have

$$= \frac{2}{\pi} \log(2\pi) - \frac{2}{\pi} [2\log 2 - 1]$$

$$= \frac{2}{\pi} \left[ \log\left(\frac{\pi}{2}\right) + 1 \right]$$

It is easily seen that

$$\int_0^1 \log \Gamma(1-x) \sin \pi x \, dx = \int_0^1 \log \Gamma(x) \sin \pi x \, dx$$

and we therefore obtain

$$(2.26) \quad \int_0^1 \log \Gamma(x) \sin \pi x \, dx = \frac{1}{\pi} \left[ \log\left(\frac{\pi}{2}\right) + 1 \right]$$

which we shall also obtain in (5.7).

Another approach is set out below. We have using Euler's reflection formula

$$\int_0^1 [\log \Gamma(x) + \log \Gamma(1-x)] \sin \pi x \, dx = \int_0^1 [\log \pi - \log \sin \pi x] \sin \pi x \, dx$$

$$= \frac{2}{\pi} \log \pi - \int_0^1 \log \sin \pi x \sin \pi x \, dx$$

so that

$$2 \int_0^1 \log \Gamma(x) \sin \pi x \, dx = \frac{2}{\pi} \log \pi - \int_0^1 \log \sin \pi x \sin \pi x \, dx$$

Dwilewicz and Mináč [28] have shown that (see also equation (6.4) in Section 6 below)

$$(2.27) \quad \frac{2p\pi}{1 - \cos 2p\pi} \int_0^1 \log(\sin \pi x) \sin 2p\pi x \, dx = -\frac{1}{2} [\psi(1+p) + \psi(1-p) + 2\gamma + 2\log 2]$$

and with $p = 1/2$ this becomes

$$\frac{\pi}{2} \int_0^1 \log(\sin \pi x) \sin \pi x \, dx = -\frac{1}{2} [\psi(3/2) + \psi(1/2) + 2\gamma + 2\log 2]$$

$$= \log 2 - 1$$

and this also results in

$$\int_0^1 \log \Gamma(x) \sin \pi x \, dx = \frac{1}{\pi} \left[ \log\left(\frac{\pi}{2}\right) + 1 \right]$$

which we will also see in (5.7).

$\square$

With $p = 1/2$ in (2.25) we obtain

$$\int_0^1 \log \Gamma(x) \sin(\pi x / 2)\, dx + \int_0^1 \log \Gamma(1-x) \sin(\pi x / 2)\, dx = \frac{2}{\pi} \log(2\pi) - \frac{2}{\pi} \sum_{n=1}^{\infty} \frac{1}{n} \frac{1}{16n^2 - 1}$$

It is easily shown that

$$\int_0^1 \log \Gamma(1-x) \sin(\pi x / 2)\, dx = \int_0^1 \log \Gamma(x) \cos(\pi x / 2)\, dx$$

and we simply obtain (2.13) again.

□

Letting $p = 2k$ in (2.8) where $k$ is an integer we get

$$\int_0^1 \log \Gamma(x) \cos 2k\pi x\, dx + \int_0^1 \log \Gamma(1-x) \cos 2k\pi x\, dx$$

$$= \frac{\sin 2k\pi}{2k\pi} \log(2\pi) + \sum_{n=1}^{\infty} \frac{1}{2n} \left[ \frac{\sin\{(n+k)2\pi\}}{(n+k)2\pi} \right] + \sum_{n=1}^{k-1} \frac{\sin\{(k-n)2\pi\}}{(k-n)2\pi}$$

$$+ \frac{1}{2k} \int_0^1 \cos(0)\, dx + \sum_{n=k+1}^{\infty} \frac{\sin\{(k-n)2\pi\}}{(k-n)2\pi}$$

We have by a simple substitution

$$\int_0^1 \log \Gamma(1-x) \cos 2k\pi x\, dx = \int_0^1 \log \Gamma(x) \cos 2k\pi x\, dx$$

so that

(2.28)    $$\int_0^1 \log \Gamma(x) \cos 2k\pi x\, dx = \frac{1}{4k}$$

Letting $p = 2k+1$ in (2.8) does not provide any information because both sides of the equation vanish.

This result (2.28) may also be obtained more directly from (2.9) by taking the limit as $p \to 2k$. It may be seen from

$$\frac{p \sin p\pi}{\pi} \sum_{n=1}^{\infty} \frac{1}{n} \frac{1}{4n^2 - p^2} = \frac{p \sin p\pi}{\pi} \left[ \frac{1}{k} \frac{1}{4k^2 - p^2} + \sum_{\substack{n=1 \\ n \neq k}}^{\infty} \frac{1}{n} \frac{1}{4n^2 - p^2} \right]$$

that all, bar one, of the terms automatically vanish and we have the limit via L'Hôpital's rule

$$\lim_{p \to 2k} \frac{1}{\pi k} \frac{p \sin p\pi}{4k^2 - p^2} = -\frac{1}{2k}$$

$\square$

We may also write

$$\int_0^1 [\log \Gamma(x) + \log \Gamma(1-x)] \cos p\pi x \, dx$$

$$= \int_0^1 [\log \pi - \log \sin \pi x] \cos p\pi x \, dx$$

$$= \frac{\sin p\pi}{p\pi} \log \pi - \int_0^1 \log \sin \pi x \cos p\pi x \, dx$$

and similarly we have

$$\int_0^1 [\log \Gamma(x) + \log \Gamma(1-x)] \sin p\pi x \, dx = \frac{\cos p\pi}{p\pi} \log \pi - \int_0^1 \log \sin \pi x \sin p\pi x \, dx$$

We note that the Fourier coefficients of $\log \sin \pi x$ are given by [33, p.577, 4.384]

(2.29)     $$\int_0^1 \log \sin \pi x \cos 2k\pi x \, dx = -\frac{1}{2k} \quad k > 0$$

(2.30)     $$\int_0^1 \log \sin \pi x \cos(2k+1)\pi x \, dx = 0$$

(2.31)     $$\int_0^1 \log \sin \pi x \, dx = -\log 2 \qquad k = 0$$

(2.32)     $$\int_0^1 \log \sin \pi x \sin 2k\pi x \, dx = 0 \qquad k \geq 0$$

$$(2.33) \qquad \int_0^1 \log \sin \pi x \sin(2k+1)\pi x\, dx = \frac{2}{(2k+1)\pi}\left[\log 2 - \frac{1}{2k+1} - 2\sum_{j=0}^{k-1}\frac{1}{2j+1}\right]$$

It should be noted that the above integral (2.33) is incorrectly reported in Nielsen's book [45, p.202]. Derivations of these integrals are given in Section 6.

We could accordingly use these Fourier coefficients of $\log \sin \pi x$ to evaluate the corresponding integrals involving the log gamma function

$\square$

Denoting $f(p)$ by

$$(2.34) \qquad f(p) = \frac{\sin p\pi}{p\pi}$$

we have

$$pf(p) = \frac{\sin p\pi}{\pi}$$

and thus we get

$$(2.35) \qquad pf'(p) + f(p) = \cos p\pi$$

or equivalently

$$f'(p) = \frac{\cos p\pi - f(p)}{p}$$

We see that $f(0) = 1$ and applying L'Hôpital's rule we find that

$$f'(0) = \lim_{p\to 0}\frac{\cos p\pi - f(p)}{p} = \lim_{p\to 0}[-\pi \sin p\pi - f'(p)] = -f'(0)$$

and thus $f'(0) = 0$. This could of course have been more readily obtained by using the elementary Maclaurin series for $f(p)$.

Differentiating (2.9) results in

$$(2.36) \qquad I'(p) + J'(p) = f'(p)\log(2\pi) - \frac{p\pi \cos p\pi + \sin p\pi}{\pi}\sum_{n=1}^{\infty}\frac{1}{n}\frac{1}{4n^2 - p^2}$$

$$-\frac{2p^2 \sin p\pi}{\pi} \sum_{n=1}^{\infty} \frac{1}{n} \frac{1}{(4n^2 - p^2)^2}$$

where we have

$$I(p) = \int_0^1 \log \Gamma(x) \cos p\pi x \, dx$$

$$J(p) = \int_0^1 \log \Gamma(1-x) \cos p\pi x \, dx$$

$$I'(p) + J'(p) = -\pi \int_0^1 x[\log \Gamma(x) + \log \Gamma(1-x)] \sin p\pi x \, dx$$

We have

$$f'(p) = \frac{\pi^2 p \cos p\pi - \pi \sin p\pi}{\pi^2 p^2}$$

and thus

$$f'(1) = -1$$

With $p = 1$ we have

$$I'(1) + J'(1) = -\pi \int_0^1 x \log \Gamma(x) \sin \pi x \, dx - \pi \int_0^1 x \log \Gamma(1-x) \sin \pi x \, dx$$

$$\int_0^1 x \log \Gamma(1-x) \sin \pi x \, dx = \int_0^1 (1-u) \log \Gamma(u) \sin[\pi(1-u)] \, du$$

$$= \int_0^1 (1-u) \log \Gamma(u) \sin \pi u \, du$$

$$= \int_0^1 \log \Gamma(u) \sin \pi u \, du - \int_0^1 u \log \Gamma(u) \sin \pi u \, du$$

Unfortunately, the terms involving $\int_0^1 x \log \Gamma(x) \sin \pi x \, dx$ cancel out and hence we have

$$I'(1) + J'(1) = -\pi \int_0^1 \log \Gamma(x) \sin \pi x \, dx$$

The right-hand side of (2.36) at $p = 1$ is

$$= -\log(2\pi) + \sum_{n=1}^{\infty} \frac{1}{n} \frac{1}{4n^2 - 1}$$

$$= -\log(2\pi) + 2\log 2 - 1$$

where we have used (2.15). Therefore we obtain another derivation of

$$\int_0^1 \log \Gamma(x) \sin \pi x \, dx = \frac{1}{\pi}\left[\log \frac{\pi}{2} + 1\right]$$

We also have

$$I'(1) + J'(1) = -\pi \log \pi \int_0^1 x \sin \pi x \, dx + \pi \int_0^1 x \log \sin \pi x . \sin \pi x \, dx$$

and therefore we obtain

(2.36.1)     $$\int_0^1 x \log \sin \pi x . \sin \pi x \, dx = \frac{1}{\pi}\left[\log 2 - 1\right]$$

□

Differentiating (2.35) gives us

$$p f''(p) + 2 f'(p) = -\pi \sin p\pi$$

or equivalently

$$f''(p) = -\frac{\pi \sin p\pi + 2 f'(p)}{p}$$

Using L'Hôpital's rule then gives us

$$f''(0) = -\lim_{p \to 0} \frac{\pi \sin p\pi + 2 f'(p)}{p} = -\lim_{p \to 0}[\pi^2 \cos p\pi + 2 f''(p)]$$

so that  we have  $f''(0) = -\frac{1}{3}\pi^2$.

Differentiating (2.36) results in

$$I''(p) + J''(p) = f''(p)\log(2\pi) - \frac{6p(p\pi\cos p\pi + \sin p\pi)}{\pi}\sum_{n=1}^{\infty}\frac{1}{n}\frac{1}{(4n^2 - p^2)^2}$$

$$-\frac{p\pi^2\sin p\pi + 2\pi\cos p\pi}{\pi}\sum_{n=1}^{\infty}\frac{1}{n}\frac{1}{4n^2 - p^2} - \frac{8p^3\sin p\pi}{\pi}\sum_{n=1}^{\infty}\frac{1}{n}\frac{1}{(4n^2 - p^2)^3}$$

Therefore we obtain with $p = 0$

$$I''(0) + J''(0) = -\frac{1}{3}\pi^2\log(2\pi) - \frac{1}{2}\varsigma(3)$$

We have from the definitions of $I(p)$ and $J(p)$

$$I''(p) + J''(p) = -\pi^2\int_0^1 x^2[\log\Gamma(x) + \log\Gamma(1-x)]\cos p\pi x\,dx$$

$$I''(0) + J''(0) = -\pi^2\int_0^1 x^2[\log\Gamma(x) + \log\Gamma(1-x)]\,dx$$

$$= -\pi^2\int_0^1 x^2[\log\pi - \log\sin\pi x]\,dx$$

$$= -\frac{1}{3}\pi^2\log\pi + \pi^2\int_0^1 x^2\log\sin\pi x\,dx$$

Hence we deduce that

$$(2.37)\qquad \int_0^1 x^2\log\sin\pi x\,dx = -\frac{1}{3}\log 2 - \frac{\varsigma(3)}{2\pi^2}$$

as previously reported in [30].

Alternatively, we have

$$I''(0) + J''(0) = -\pi^2\int_0^1 x^2\log\Gamma(x)\,dx - \pi^2\int_0^1 x^2\log\Gamma(1-x)\,dx$$

We see that

$$\int\limits_0^1 x^2 \log \Gamma(1-x)\, dx = \int\limits_0^1 (1-u)^2 \log \Gamma(u)\, du$$

so that we have

$$I''(0) + J''(0) = -2\pi^2 \int\limits_0^1 x^2 \log \Gamma(x)\, dx + 2\pi^2 \int\limits_0^1 x \log \Gamma(x)\, dx - \pi^2 \int\limits_0^1 \log \Gamma(x)\, dx$$

We have the well known integral (see for example (3.25) and Appendix C)

(2.38) $$\int\limits_0^1 x \log \Gamma(x)\, dx = \frac{1}{4}\log(2\pi) - \log A$$

so that

$$-2\pi^2 \int\limits_0^1 x^2 \log \Gamma(x)\, dx + 2\pi^2 \left[\frac{1}{4}\log(2\pi) - \log A\right] - \frac{1}{2}\pi^2 \log(2\pi) = -\frac{1}{3}\pi^2 \log(2\pi) - \frac{1}{2}\varsigma(3)$$

Hence we deduce the known integral

(2.39) $$\int\limits_0^1 x^2 \log \Gamma(x)\, dx = \frac{1}{6}\log(2\pi) - \log A + \frac{\varsigma(3)}{4\pi^2}$$

which is also shown in (C.6) in Appendix C.

□

Alternatively, we now rewrite (2.9) as

$$\sum_{n=1}^{\infty} \frac{1}{n} \frac{1}{4n^2 - p^2} = \frac{\log(2\pi)\sin p\pi - \pi p[I(p) + J(p)]}{p^2 \sin p\pi}$$

and we have the limit as $p \to 0$

$$\frac{1}{4}\varsigma(3) = \lim_{p \to 0} \frac{\log(2\pi)\sin p\pi - \pi p[I(p) + J(p)]}{p^2 \sin p\pi}$$

Applying L'Hôpital's rule three times we have

$$= \lim_{p \to 0} \frac{\pi \log(2\pi)\cos p\pi - \pi p[I'(p) + J'(p)] - \pi[I(p) + J(p)]}{\pi p^2 \cos p\pi + 2p \sin p\pi}$$

$$= \lim_{p \to 0} \frac{-\pi^2 \log(2\pi) \sin p\pi - \pi p[I''(p) + J''(p)] - 2\pi[I'(p) + J'(p)]}{-\pi^2 p^2 \sin p\pi + 4p\pi \cos p\pi + 2\sin p\pi}$$

$$= \lim_{p \to 0} \frac{-\pi^3 \log(2\pi) \cos p\pi - \pi p[I'''(p) + J'''(p)] - 3\pi[I''(p) + J''(p)]}{-\pi^3 p^2 \cos p\pi - 6\pi^2 p \sin p\pi + 6\pi \cos p\pi}$$

so that

$$\frac{1}{4}\varsigma(3) = -\frac{1}{6}\pi^2 \log(2\pi) - \frac{1}{2}[I''(0) + J''(0)]$$

and this gives us another derivation of (2.39). See also Section 6.

$\square$

We recall (2.25)

$$\int_0^1 \log\Gamma(x) \sin p\pi x \, dx + \int_0^1 \log\Gamma(1-x) \sin p\pi x \, dx$$

$$= \frac{1-\cos p\pi}{p\pi} \log(2\pi) - \frac{p(1-\cos p\pi)}{\pi} \sum_{n=1}^{\infty} \frac{1}{n} \frac{1}{4n^2 - p^2}$$

and differentiating gives us

(2.40) $\pi \int_0^1 x \log\Gamma(x) \cos p\pi x \, dx + \pi \int_0^1 x \log\Gamma(1-x) \cos p\pi x \, dx$

$$= \frac{p\pi^2 \sin p\pi - \pi(1-\cos p\pi)}{(p\pi)^2} \log(2\pi) - \frac{2p^2(1-\cos p\pi)}{\pi} \sum_{n=1}^{\infty} \frac{1}{n} \frac{1}{(4n^2 - p^2)^2}$$

$$- \frac{p\pi \sin p\pi + 1 - \cos p\pi}{\pi} \sum_{n=1}^{\infty} \frac{1}{n} \frac{1}{4n^2 - p^2}$$

Letting $p = 0$ in (2.40) we easily determine that the left-hand side becomes $\frac{\pi}{2} \log(2\pi)$, which, using L'Hôpital's rule, is also seen to be the only surviving term on the right-hand side.

With $p = 1$ in (2.40) we have

$$(2.41) \qquad \pi \int_0^1 x \log \Gamma(x) \cos \pi x \, dx + \pi \int_0^1 x \log \Gamma(1-x) \cos \pi x \, dx$$

$$= -\frac{2}{\pi} \log(2\pi) - \frac{4}{\pi} \sum_{n=1}^\infty \frac{1}{n} \frac{1}{(4n^2-1)^2} - \frac{2}{\pi} \sum_{n=1}^\infty \frac{1}{n} \frac{1}{4n^2-1}$$

Differentiating (2.23) results in

$$(2.42) \qquad \frac{1}{2} \psi'\left(1+\frac{p}{2}\right) - \frac{1}{2} \psi'\left(1-\frac{p}{2}\right) = -4p^3 \sum_{n=1}^\infty \frac{1}{n} \frac{1}{(4n^2-p^2)^2} - 4p \sum_{n=1}^\infty \frac{1}{n} \frac{1}{4n^2-p^2}$$

and with $p=1$ we obtain

$$\frac{1}{2} \psi'\left(1+\frac{1}{2}\right) - \frac{1}{2} \psi'\left(\frac{1}{2}\right) = -4 \sum_{n=1}^\infty \frac{1}{n} \frac{1}{(4n^2-1)^2} - 4 \sum_{n=1}^\infty \frac{1}{n} \frac{1}{4n^2-1}$$

Since $\psi(1+x) - \psi(x) = \frac{1}{x}$ we have $\psi'(1+x) - \psi'(x) = -\frac{1}{x^2}$ and thus

$$\frac{1}{2} \psi'\left(1+\frac{1}{2}\right) - \frac{1}{2} \psi'\left(\frac{1}{2}\right) = -2$$

Hence using (2.15) we obtain

$$(2.43) \qquad \sum_{n=1}^\infty \frac{1}{n} \frac{1}{(4n^2-1)^2} = \frac{3}{2} - 2\log 2$$

which is reported in [33, p.9].

The right-hand side of (2.41) becomes using (2.15) and (2.43)

$$= -\frac{2}{\pi} \log(2\pi) - \frac{4}{\pi} \left[\frac{3}{2} - 2\log 2\right] - \frac{2}{\pi}(2\log 2 - 1)$$

$$= -\frac{2}{\pi} \left[\log \frac{\pi}{2} + 2\right]$$

Having regard to (2.40), we see that with $p=1$

$$\int_0^1 x \log \Gamma(1-x) \cos \pi x \, dx = \int_0^1 (1-u) \log \Gamma(u) \cos \pi(1-u) \, du$$

$$= \int\limits_0^1 (u-1)\log\Gamma(u)\cos\pi u\,du$$

which gives us

$$2\pi\int\limits_0^1 x\log\Gamma(x)\cos\pi x\,dx - \pi\int\limits_0^1 \log\Gamma(x)\cos\pi x\,dx = 2\pi\int\limits_0^1 \left(x-\frac{1}{2}\right)\log\Gamma(x)\cos\pi x\,dx$$

and we have the result

(2.44) $$\int\limits_0^1 \left(x-\frac{1}{2}\right)\log\Gamma(x)\cos\pi x\,dx = -\frac{1}{\pi^2}\left[\log\frac{\pi}{2}+2\right]$$

We note from (4.5.2) that

$$\int\limits_0^1 \log\Gamma(x)\cos\pi x\,dx = \frac{4}{\pi^2}\sum_{n=1}^\infty \frac{\gamma+\log(2\pi n)}{4n^2-1}$$

so that we have

(2.45) $$\int\limits_0^1 x\log\Gamma(x)\cos\pi x\,dx = \frac{2}{\pi^2}\sum_{n=1}^\infty \frac{\gamma+\log(2\pi n)}{4n^2-1} - \frac{1}{\pi^2}\left[\log\frac{\pi}{2}+2\right]$$

This may also be obtained by differentiating (5.5) with respect to $p$. This gives us

$$\int\limits_0^1 x\log\Gamma(x)\cos\pi x\,dx$$

$$= -\frac{1}{\pi^2}\log(2\pi) - \frac{2}{\pi^2}\sum_{n=1}^\infty \frac{1}{n}\frac{1}{(4n^2-1)^2} - \frac{1}{\pi^2}\sum_{n=1}^\infty \frac{1}{n}\frac{1}{4n^2-1} + \frac{2}{\pi^2}\sum_{n=1}^\infty \frac{\gamma+\log(2\pi n)}{4n^2-1}$$

$$= -\frac{1}{\pi^2}\log(2\pi) - \frac{2}{\pi^2}\left[\frac{3}{2}-2\log 2\right] - \frac{1}{\pi^2}(2\log 2-1) + \frac{2}{\pi^2}\sum_{n=1}^\infty \frac{\gamma+\log(2\pi n)}{4n^2-1}$$

$$= \frac{2}{\pi^2}\sum_{n=1}^\infty \frac{\gamma+\log(2\pi n)}{4n^2-1} - \frac{1}{\pi^2}\left[\log\frac{\pi}{2}+2\right]$$

$\square$

As noted in Section 4 we have

$$\frac{\partial}{\partial s}\int_0^1 \varsigma(s,x)\varsigma(0,x)\cos \pi x \, dx \bigg|_{s=0} = \frac{1}{2}\log(2\pi)\int_0^1 B_1(x)\cos \pi x \, dx$$

$$-\int_0^1 B_1(x)\log\Gamma(x)\cos \pi x \, dx$$

where $B_n(x)$ are the Bernoulli polynomials and in particular we have

$$B_1(x) = x - \frac{1}{2}$$

□

Differentiating (2.42) results in

$$\frac{1}{4}\psi''\left(1+\frac{p}{2}\right)+\frac{1}{4}\psi''\left(1-\frac{p}{2}\right) = -16p^4\sum_{n=1}^{\infty}\frac{1}{n}\frac{1}{(4n^2-p^2)^3} - 20p^2\sum_{n=1}^{\infty}\frac{1}{n}\frac{1}{(4n^2-p^2)^2}$$

$$-4\sum_{n=1}^{\infty}\frac{1}{n}\frac{1}{4n^2-p^2}$$

and with $p=1$ we obtain

$$\psi''\left(1+\frac{1}{2}\right)+\psi''\left(\frac{1}{2}\right) = -64\sum_{n=1}^{\infty}\frac{1}{n}\frac{1}{(4n^2-1)^3} - 80\sum_{n=1}^{\infty}\frac{1}{n}\frac{1}{(4n^2-1)^2} - 16\sum_{n=1}^{\infty}\frac{1}{n}\frac{1}{4n^2-1}$$

It is well known that (see [1, 6.4.4] and [41])

$$\psi^{(k)}\left(\frac{1}{2}\right) = (-1)^{k+1}k!(2^{k+1}-1)\varsigma(k+1)$$

$$\psi^{(2)}\left(\frac{1}{2}\right) = -14\varsigma(3)$$

$$\psi''(1+x) - \psi''(x) = \frac{2}{x^3}$$

$$\psi''\left(1+\frac{1}{2}\right)+\psi''\left(\frac{1}{2}\right) = 16 - 28\varsigma(3)$$

$$28\zeta(3) - 16 = 64\sum_{n=1}^{\infty}\frac{1}{n}\frac{1}{(4n^2-1)^3} + 80\sum_{n=1}^{\infty}\frac{1}{n}\frac{1}{(4n^2-1)^2} + 16\sum_{n=1}^{\infty}\frac{1}{n}\frac{1}{4n^2-1}$$

and using (2.15) and (2.43) this becomes

$$28\zeta(3) - 16 = 64\sum_{n=1}^{\infty}\frac{1}{n}\frac{1}{(4n^2-1)^3} + 80\left[\frac{3}{2} - 2\log 2\right] + 16\left[2\log 2 - 1\right]$$

which gives us

$$(2.46) \qquad 16\sum_{n=1}^{\infty}\frac{1}{n}\frac{1}{(4n^2-1)^3} = 7\zeta(3) + 32\log 2 - 30$$

$\square$

We recall (2.24)

$$\int_0^1 \log\Gamma(x)\cos p\pi x\,dx + \int_0^1 \log\Gamma(1-x)\cos p\pi x\,dx$$

$$= \frac{\sin p\pi}{p\pi}\left[\gamma + \log(2\pi) + \frac{1}{2}\psi\left(1+\frac{p}{2}\right) + \frac{1}{2}\psi\left(1-\frac{p}{2}\right)\right]$$

and differentiation results in

$$-\pi\int_0^1 x\log\Gamma(x)\sin p\pi x\,dx - \pi\int_0^1 x\log\Gamma(1-x)\sin p\pi x\,dx$$

$$= \frac{\sin p\pi}{p\pi}\left[\frac{1}{2^2}\psi'\left(1+\frac{p}{2}\right) - \frac{1}{2^2}\psi'\left(1-\frac{p}{2}\right)\right]$$

$$+ \frac{\pi^2 p\cos p\pi - \pi\sin p\pi}{\pi^2 p^2}\left[\gamma + \log(2\pi) + \frac{1}{2}\psi\left(1+\frac{p}{2}\right) + \frac{1}{2}\psi\left(1-\frac{p}{2}\right)\right]$$

With $p = 1$ this becomes

$$-\pi\int_0^1 x\log\Gamma(x)\sin\pi x\,dx - \pi\int_0^1 x\log\Gamma(1-x)\sin\pi x\,dx$$

$$= -\left[\gamma + \log(2\pi) + \frac{1}{2}\psi\left(1+\frac{1}{2}\right) + \frac{1}{2}\psi\left(\frac{1}{2}\right)\right]$$

$$= -\left[\log \frac{\pi}{2} + 1\right]$$

We see that

$$\int_0^1 x \log \Gamma(1-x) \sin \pi x \, dx = \int_0^1 (1-x) \log \Gamma(x) \sin \pi x \, dx$$

$$-\pi \int_0^1 x \log \Gamma(x) \sin \pi x \, dx - \pi \int_0^1 x \log \Gamma(1-x) \sin \pi x \, dx = -\pi \int_0^1 \log \Gamma(x) \sin \pi x \, dx$$

and hence we obtain

$$\int_0^1 \log \Gamma(x) \sin \pi x \, dx = \frac{1}{\pi}\left[\log \frac{\pi}{2} + 1\right]$$

**An observation by Glasser**

Part of the following is based on an observation made by Glasser [32] in 1966. Let us consider the integral

$$I = \int_0^1 f(x) \psi(x) dx$$

where $f(x) = -f(1-x)$ and $f(x)$ is selected so that the integral converges. We then have

$$I = \int_0^1 f(1-t) \psi(1-t) dt = -\int_0^1 f(t) \psi(1-t) dt$$

and hence we see that

$$2I = \int_0^1 f(x)[\psi(x) - \psi(1-x)] dx$$

Therefore, since $\psi(x) - \psi(1-x) = -\pi \cot \pi x$, we have

(2.47) $$\int_0^1 f(x) \psi(x) dx = -\frac{\pi}{2} \int_0^1 f(x) \cot \pi x \, dx = -\pi \int_0^{1/2} f(x) \cot \pi x \, dx$$

Glasser [32] used the function $f(x) = x(1-x)\cos\pi x$ to show that

(2.48) $$\int_0^1 \psi(x)x(1-x)\cos\pi x\,dx = \frac{1}{\pi^2}\left[2 - \frac{7}{2}\varsigma(3)\right]$$

Integration by parts gives us

$$\int_0^1 \psi(x)x(1-x)\cos\pi x\,dx = \log\Gamma(x)x(1-x)\cos\pi x\big|_0^1$$
$$+\int_0^1 \log\Gamma(x)[x(1-x)\pi\sin\pi x + (2x-1)\cos\pi x]\,dx$$

Since

$$x\log\Gamma(x) = x\log\Gamma(1+x) - x\log x$$

we see that the integrated part vanishes and we obtain

$$\int_0^1 \log\Gamma(x)[x(1-x)\pi\sin\pi x + (2x-1)\cos\pi x]\,dx = \frac{1}{\pi^2}\left[2 - \frac{7}{2}\varsigma(3)\right]$$

Unfortunately, if we substitute the relevant component integrals (which are listed in Appendix D) we find that the interesting summations cancel out.

Some functions which satisfy the condition $f(x) = -f(1-x)$ are listed below:

$$B_{2n+1}(x)\,,\,B_{2n}(x)\cos(2m+1)\pi x\,,\,B_{2n+1}(x)\sin 2m\pi x\text{ and }h\big(x(1-x)\big)\cos(2m+1)\pi x$$

since $B_n(x) = (-1)^n B_n(1-x)$ where $n$ and $m$ are positive integers.

Similarly, if $g(x) = g(1-x)$, we obtain

$$\int_0^1 g(x)[\psi(x) - \psi(1-x)]\,dx = \int_0^1 g(x)\cot\pi x\,dx = 0$$

$$I = \int_0^1 f(x)\psi(x)\,dx = f(x)\log\Gamma(x)\big|_0^1 - \int_0^1 f'(x)\log\Gamma(x)\,dx$$

We have for suitably behaved functions

$$\int\limits_0^1 f(x)\psi(x)\,dx = -\int\limits_0^1 f'(x)\log\Gamma(x)\,dx$$

and an example of this is

$$\int\limits_0^1 B_{2n+1}(x)\psi(x)\,dx = B_{2n+1}(x)\log\Gamma(x)\big|_0^1 - (2n+1)\int\limits_0^1 B_{2n}(x)\log\Gamma(x)\,dx$$

$$= -(2n+1)\int\limits_0^1 B_{2n}(x)\log\Gamma(x)\,dx$$

Glasser's formula gives us

(2.49) $$\int\limits_0^1 B_{2n+1}(x)\psi(x)\,dx = -\frac{\pi}{2}\int\limits_0^1 B_{2n+1}(x)\cot\pi x\,dx$$

and hence we have

(2.50) $$\int\limits_0^1 B_{2n}(x)\log\Gamma(x)\,dx = \frac{\pi}{2(2n+1)}\int\limits_0^1 B_{2n+1}(x)\cot\pi x\,dx$$

The following integral appears in Abramowitz and Stegun [1, p.807]

(2.51) $$\int\limits_0^1 B_{2n+1}(x)\cot\pi x\,dx = \frac{(-1)^{n+1}2(2n+1)!\varsigma(2n+1)}{(2\pi)^{2n+1}}$$

and there are many derivations; for example, see the recent one by Dwilewicz and Mináč [28]. A similar identity was also derived by Espinosa and Moll [30] in the form

(2.52) $$\int\limits_0^1 B_{2n}(x)\log\sin\pi x\,dx = (-1)^n\frac{(2n)!\varsigma(2n+1)}{(2\pi)^{2n}}$$

and it is easily shown that equation (2.52) above is equivalent to (2.51) following a simple integration by parts.

 Hence we obtain

(2.53) $$\int\limits_0^1 B_{2n}(x)\log\Gamma(x)\,dx = \frac{(-1)^{n+1}(2n)!\varsigma(2n+1)}{2(2\pi)^{2n+1}}$$

in agreement with [30].

□

Some of the above analysis could be replicated by employing the well-known relation

$$\Gamma\left(x+\frac{1}{2}\right)\Gamma\left(\frac{1}{2}-x\right) = \frac{\pi}{\cos\pi x}$$

which results from Euler's reflection formula with the substitution $x \rightarrow x+1/2$ .

□

Using the Maclaurin series for $\cos\pi x$ we have

$$\int_0^1 \log\Gamma(x)\cos\pi x\,dx = \sum_{n=0}^{\infty}(-1)^n\frac{\pi^{2n}}{(2n)!}\int_0^1 x^{2n}\log\Gamma(x)\,dx$$

From [30] we have the inversion formula for the Bernoulli polynomials

$$x^j = \frac{1}{j+1}\sum_{k=0}^{j}\binom{j+1}{k}B_k(x)$$

$$\int_0^1 x^{2n}\log\Gamma(x)\,dx = \frac{1}{2n+1}\sum_{k=0}^{2n}\binom{2n+1}{k}\int_0^1 B_k(x)\log\Gamma(x)\,dx$$

We have from equation (6.10) of [30]

$$\int_0^1 B_k(x)\log\Gamma(x)\,dx = H_k\varsigma(-k) + (-1)^{k+1}\varsigma'(-k)$$

and thus we have

$$\int_0^1 \log\Gamma(x)\cos\pi x\,dx = \sum_{n=0}^{\infty}(-1)^n\frac{\pi^{2n}}{(2n+1)!}\sum_{k=0}^{2n}\binom{2n+1}{k}\left[H_k\varsigma(-k) + (-1)^{k+1}\varsigma'(-k)\right]$$

We note that via integration by parts this links up with equation (12.9) of [30]

$$\int_0^1 x^{2n+1}\psi(x)\,dx = \sum_{k=0}^{2n}\binom{2n+1}{k}\left[H_k\varsigma(-k) + (-1)^{k+1}\varsigma'(-k)\right]$$

## 3. Some applications of the Fourier series for the Hurwitz zeta function

In Titchmarsh's treatise [54, p.37] we have Hurwitz's formula for the Fourier series of the Hurwitz zeta function $\varsigma(s, x)$

$$(3.1) \qquad \varsigma(s, x) = 2\Gamma(1-s)\left[\sin\left(\frac{\pi s}{2}\right)\sum_{k=1}^{\infty}\frac{\cos 2k\pi x}{(2\pi k)^{1-s}} + \cos\left(\frac{\pi s}{2}\right)\sum_{n=1}^{\infty}\frac{\sin 2k\pi x}{(2\pi k)^{1-s}}\right]$$

where $\text{Re}(s) < 0$ and $0 < x \leq 1$. In 2000, Boudjelkha [14] showed that this formula also applies in the region $\text{Re}(s) < 1$. It may be noted that when $x = 1$ this reduces to Riemann's functional equation for $\varsigma(s)$.

From the above are immediately found the corresponding Fourier coefficients

$$(3.2) \qquad \int_0^1 \varsigma(s, x)\sin 2n\pi x\, dx = \frac{\Gamma(1-s)}{(2\pi n)^{1-s}}\cos\left(\frac{\pi s}{2}\right)$$

$$(3.3) \qquad \int_0^1 \varsigma(s, x)\cos 2n\pi x\, dx = \frac{\Gamma(1-s)}{(2\pi n)^{1-s}}\sin\left(\frac{\pi s}{2}\right)$$

Using Euler's reflection formula (2.6) these may be written as (see [30])

$$(3.4) \qquad \int_0^1 \varsigma(s, x)\sin 2n\pi x\, dx = \frac{(2\pi)^s n^{s-1}}{4\Gamma(s)}\csc\left(\frac{\pi s}{2}\right)$$

$$(3.5) \qquad \int_0^1 \varsigma(s, x)\cos 2n\pi x\, dx = \frac{(2\pi)^s n^{s-1}}{4\Gamma(s)}\sec\left(\frac{\pi s}{2}\right)$$

Having regard to (3.2) we consider

$$f(s) = \frac{\Gamma(1-s)}{(2\pi n)^{1-s}}\cos\left(\frac{\pi s}{2}\right)$$

and logarithmic differentiation gives us

$$\frac{f'(s)}{f(s)} = \log(2\pi n) - \psi(1-s) - \frac{\pi}{2}\tan\left(\frac{\pi s}{2}\right)$$

Therefore differentiating (3.2) results in (this idea was inspired by a paper by Espinosa and Moll [30])

(3.6)    $\int_0^1 \varsigma'(s,x) \sin 2n\pi x \, dx = \left[ \log(2\pi n) - \psi(1-s) - \frac{\pi}{2} \tan\left(\frac{\pi s}{2}\right) \right] \frac{\Gamma(1-s)}{(2\pi n)^{1-s}} \cos\left(\frac{\pi s}{2}\right)$

Applying Lerch's identity [10]

(3.7)    $\varsigma'(0,x) = \log \Gamma(x) - \frac{1}{2} \log(2\pi)$

we then have as $s \to 0$

$$\int_0^1 \left[ \log \Gamma(x) - \frac{1}{2} \log(2\pi) \right] \sin 2n\pi x \, dx = \frac{\log(2\pi n) + \gamma}{2\pi n}$$

and this becomes the well-known result which we derived in a different manner in (1.9.2)

(3.8)    $\int_0^1 \log \Gamma(x) \sin 2n\pi x \, dx = \frac{\log(2\pi n) + \gamma}{2\pi n}$

Similarly differentiating (3.3) gives us

(3.9)    $\int_0^1 \varsigma'(s,x) \cos 2n\pi x \, dx = \left[ \log(2\pi n) - \psi(1-s) + \frac{\pi}{2} \cot\left(\frac{\pi s}{2}\right) \right] \frac{\Gamma(1-s)}{(2\pi n)^{1-s}} \sin\left(\frac{\pi s}{2}\right)$

and as $s \to 0$ we obtain

$$\int_0^1 \left[ \log \Gamma(x) - \frac{1}{2} \log(2\pi) \right] \cos 2n\pi x \, dx = \frac{1}{4n}$$

This then gives us (1.9)

(3.10)    $\int_0^1 \log \Gamma(x) \cos 2n\pi x \, dx = \frac{1}{4n}$

We have thereby obtained the Fourier coefficients for $\log \Gamma(x)$ and, in the process, have thus rediscovered Kummer's Fourier series for $\log \Gamma(x)$

(3.11)    $\log \Gamma(x) = \frac{1}{2} \log \pi - \frac{1}{2} \log \sin \pi x + \sum_{n=1}^{\infty} \frac{(\gamma + \log 2\pi n) \sin 2\pi n x}{\pi n}$   $(0 < x < 1)$

(this formula was originally derived by Kummer in 1847 [42]). Reference to (2.7) confirms that (3.11) is properly described as a Fourier series expansion for $\log \Gamma(x)$.

Kummer's formula may also be written as (cf. Nielsen [45, p.201])

$$(3.12) \qquad \log \Gamma(x) = \frac{1}{2} \log \frac{\pi}{\sin \pi x} + \left( \frac{1}{2} - x \right)(\gamma + \log 2\pi) + \frac{1}{\pi} \sum_{n=1}^{\infty} \frac{\log n}{n} \sin 2\pi n x$$

Nörlund [46, p.118] reports that for $0 < x < 1$

$$(3.13)$$

$$\log \Gamma \left( \frac{x}{2} \right) - \log \Gamma \left( \frac{1+x}{2} \right) = \frac{1}{2}(\gamma + \log 2\pi) + \frac{1}{2} \log \cot \left( \frac{\pi x}{2} \right) + \frac{2}{\pi} \sum_{n=0}^{\infty} \frac{\log(2n+1)}{2n+1} \sin(2n+1)\pi x$$

(this was also derived by Williams and Zhang [57] using Kummer's formula (3.12)).

Integrating this results in

$$\int_0^1 \log \Gamma \left( \frac{x}{2} \right) dx - \int_0^1 \log \Gamma \left( \frac{1+x}{2} \right) dx = \frac{1}{2}(\gamma + \log 2\pi) + \frac{1}{2} \int_0^1 \log \cot \left( \frac{\pi x}{2} \right) dx + \frac{4}{\pi^2} \sum_{n=0}^{\infty} \frac{\log(2n+1)}{(2n+1)^2}$$

We see that

$$\int_0^1 \log \cot \left( \frac{\pi x}{2} \right) dx = \int_0^1 \log \cos \left( \frac{\pi x}{2} \right) dx - \int_0^1 \log \sin \left( \frac{\pi x}{2} \right) dx$$

and we have the integrals [33, p.526]

$$\int_0^1 \log \cos \left( \frac{\pi x}{2} \right) dx = \frac{2}{\pi} \int_0^{\pi/2} \log \cos u \, du = -\log 2$$

$$\int_0^1 \log \sin \left( \frac{\pi x}{2} \right) dx = \frac{2}{\pi} \int_0^{\pi/2} \log \sin u \, du = -\log 2$$

We have

$$\int_0^1 \log \Gamma \left( \frac{x}{2} \right) dx = 2 \int_0^{1/2} \log \Gamma(u) \, du$$

and using (1.64) this becomes

$$\int_0^1 \log \Gamma\left(\frac{x}{2}\right) dx = \frac{5}{12}\log 2 + \frac{1}{2}\log \pi + 3\log A$$

Hence we obtain

$$(3.13.1) \quad \int_0^1 \log \Gamma\left(\frac{1+x}{2}\right) dx = 3\log A - \frac{1}{12}\log 2 - \frac{1}{2}\gamma - \frac{4}{\pi^2}\sum_{n=0}^{\infty}\frac{\log(2n+1)}{(2n+1)^2}$$

We have for $\text{Re}(s) > 1$

$$(1-2^{-s})\varsigma(s) = \sum_{n=0}^{\infty}\frac{1}{(2n+1)^s}$$

and differentiation results in

$$(1-2^{-s})\varsigma'(s) + 2^{-s}\varsigma(s)\log 2 = -\sum_{n=0}^{\infty}\frac{\log(2n+1)}{(2n+1)^s}$$

We may therefore write (3.13.1) as

$$(3.14) \quad \int_0^1 \log \Gamma\left(\frac{1+x}{2}\right) dx = 3\log A - \frac{1}{12}\log 2 - \frac{1}{2}\gamma - \frac{1}{\pi^2}[3\varsigma'(2) + \varsigma(2)\log 2]$$

**Another proof of Hurwitz's formula for the Fourier series expansion of the Hurwitz zeta function**

The following proof of Hurwitz's formula for the Fourier series expansion of the Hurwitz zeta function is based on Kummer's method which is outlined in [5, p.29].

Since $\log \Gamma(x)$ is differentiable in $(0,1)$ it has a Fourier expansion

$$\log \Gamma(x) = C_0 + 2\sum_{k=1}^{\infty}C_k \cos 2k\pi x + 2\sum_{k=1}^{\infty}D_k \sin 2k\pi x$$

where $C_k = \int_0^1 \log \Gamma(x)\cos 2k\pi x\, dx$ and $D_k = \int_0^1 \log \Gamma(x)\sin 2k\pi x\, dx$

The $C_k$ coefficients are relatively easy to determine using Euler's reflection formula for the gamma function

$$\log \Gamma(x) + \log \Gamma(1-x) = \log 2\pi - \log\left[2\sin \pi x\right]$$

$$= \log 2\pi + \sum_{n=1}^{\infty} \frac{\cos 2n\pi x}{n}$$

The corresponding Fourier series is

$$\log \Gamma(x) + \log \Gamma(1-x) = 2C_0 + 4\sum_{k=1}^{\infty} C_k \cos 2k\pi x$$

and equating coefficients in the last two equations gives

(3.15) $$C_0 = \frac{1}{2}\log 2\pi \quad \text{and} \quad C_k = \frac{1}{4k} \text{ for } k \geq 1$$

We now refer to Malmstén's formula [52, p.16] (which is also derived in equation (E.22g) of [22])

(3.16) $$\log \Gamma(x) = \int_0^{\infty} \left[ (x-1)e^{-\alpha} - \frac{e^{-\alpha} - e^{-\alpha x}}{1 - e^{-\alpha}} \right] \frac{d\alpha}{\alpha}$$

and a change of variables $u = e^{-\alpha}$ gives us

(3.17) $$\log \Gamma(x) = \int_0^1 \left[ \frac{1 - u^{x-1}}{1 - u} - x + 1 \right] \frac{du}{\log u}$$

We therefore obtain

$$D_k = \int_0^1 \log \Gamma(x) \sin 2k\pi x \, dx$$

$$= \int_0^1 \int_0^1 \left[ \frac{1 - u^{x-1}}{1 - u} - x + 1 \right] \frac{\sin 2k\pi x \, du \, dx}{\log u}$$

We have

$$\int_0^1 \sin 2k\pi x \, dx = 0$$

$$\int_0^1 x \sin 2k\pi x \, dx = -\frac{1}{2k\pi}$$

$$\int_0^1 u^{x-1} \sin 2k\pi x \, dx = \frac{1}{u} \operatorname{Im} \int_0^1 \exp\left[x(\log u + 2k\pi i)\right] dx$$

$$= \frac{1}{u} \operatorname{Im} \frac{u-1}{\log u + 2k\pi i}$$

Therefore we get

(3.18) $$D_k = \int_0^1 \left[ \frac{-2k\pi}{u\left(\log^2 u + 4k^2\pi^2\right)} + \frac{1}{2k\pi} \right] \frac{du}{\log u}$$

With $u = e^{-2k\pi t}$ we obtain

$$D_k = \frac{1}{2\pi k} \int_0^\infty \left[ \frac{1}{\left(1+t^2\right)} - e^{-2k\pi t} \right] \frac{dt}{t}$$

Taking $k = 1$ we have

$$D_1 = \frac{1}{2\pi} \int_0^\infty \left[ \frac{1}{\left(1+t^2\right)} - e^{-2\pi t} \right] \frac{dt}{t}$$

With $x = 1$ in Malmstén's formula (3.16) we get

$$-\frac{\gamma}{2\pi} = \frac{1}{2\pi} \int_0^\infty \left[ e^{-t} - \frac{1}{1+t} \right] \frac{dt}{t}$$

Hence we have

(3.19) $$D_1 - \frac{\gamma}{2\pi} = \frac{1}{2\pi} \int_0^\infty \frac{e^{-t} - e^{-2\pi t}}{t} dt + \frac{1}{2\pi} \int_0^\infty \left[ \frac{1}{\left(1+t^2\right)} - \frac{1}{\left(1+t\right)} \right] \frac{dt}{t}$$

We have $1/r = \int_0^\infty e^{-rx} dx$ ($r > 0$) and, integrating that expression, we obtain Frullani's integral

$$\int_1^b \frac{dr}{r} = \int_1^n dr \int_0^\infty e^{-rx} dx = \int_0^\infty dx \int_1^b e^{-rx} dr$$

which implies that

$$(3.20) \qquad \log b = \int\limits_{0}^{\infty} \frac{e^{-x} - e^{-bx}}{x} \, dx$$

By (3.20) the first integral in (3.19) is equal to $\log 2\pi$ and a change of variables $t \to 1/t$ shows that the second integral vanishes. Therefore we have

$$D_1 = \frac{\gamma}{2\pi} + \frac{1}{2\pi} \log 2\pi$$

$D_k$ is found by noting that

$$kD_k - D_1 = \frac{1}{2\pi} \int\limits_{0}^{\infty} \frac{e^{-2\pi t} - e^{-2k\pi t}}{t} \, dt = \frac{1}{2\pi} \log k$$

where the integral is also evaluated with (3.20). Thus we have

$$(3.21) \qquad D_k = \frac{\gamma + \log 2k\pi}{2\pi k} = \int\limits_{0}^{1} \log \Gamma(x) \sin 2k\pi x \, dx$$

and Kummer's formula thereby follows.

$\square$

We see from (3.18) and (3.21) that

$$\int\limits_{0}^{1} \left[ \frac{1}{(2\pi k)^2} - \frac{1}{u\left(\log^2 u + 4k^2\pi^2\right)} \right] \frac{du}{\log u} = \frac{\gamma + \log(2k\pi)}{(2\pi k)^2}$$

and we make the summation

$$\sum_{k=1}^{\infty} \int\limits_{0}^{1} \left[ \frac{1}{(2\pi k)^2} - \frac{1}{u\left(\log^2 u + 4k^2\pi^2\right)} \right] \frac{du}{\log u} = \sum_{k=1}^{\infty} \frac{\gamma + \log(2k\pi)}{(2\pi k)^2}$$

Assuming that interchanging of the order of integration and summation is valid we have

$$\int\limits_{0}^{1} \left[ \frac{\varsigma(2)}{(2\pi)^2} - \frac{1}{u} \sum_{k=1}^{\infty} \frac{1}{\log^2 u + 4k^2\pi^2} \right] \frac{du}{\log u} = \sum_{k=1}^{\infty} \frac{\gamma + \log(2k\pi)}{(2\pi k)^2}$$

Letting $\alpha \to i\alpha$ in (4.7) gives us

$$\pi \coth \alpha \pi = \frac{1}{\alpha} + 2\alpha \sum_{k=1}^{\infty} \frac{1}{k^2 + \alpha^2}$$

and we obtain

$$\frac{1}{(2\pi)^2} \sum_{k=1}^{\infty} \frac{1}{k^2 + (\log u / 2\pi)^2} = \frac{1}{4\pi} \frac{\frac{1}{2} \log u \coth[(\log u)/2] - 1}{\log u}$$

$$= \frac{\log u \coth[(\log u)/2] - 2}{8\pi \log u}$$

This gives us

$$\int_0^1 \left[ \frac{\varsigma(2)}{(2\pi)^2} - \frac{1}{u} \sum_{k=1}^{\infty} \frac{1}{\log^2 u + 4k^2\pi^2} \right] \frac{du}{\log u} = \int_0^1 \left[ \frac{\varsigma(2)}{(2\pi)^2} - \frac{\log u \coth[(\log u)/2] - 2}{8\pi u \log u} \right] \frac{du}{\log u}$$

and thus

$$\int_0^1 \left[ \frac{\varsigma(2)}{(2\pi)^2} - \frac{\log u \coth[(\log u)/2] - 2}{8\pi u \log u} \right] \frac{du}{\log u} = \sum_{k=1}^{\infty} \frac{\gamma + \log(2k\pi)}{(2\pi k)^2}$$

Since $\coth[(\log u)/2] = \dfrac{u+1}{u-1}$ we may write this as

$$\frac{1}{4\pi} \int_0^1 \left[ \frac{\pi}{6} + \frac{1}{2} \frac{1+u}{u(1-u)} + \frac{1}{u \log u} \right] \frac{du}{\log u} = \sum_{k=1}^{\infty} \frac{\gamma + \log(2k\pi)}{(2\pi k)^2}$$

It may be noted that this bears a structural similarity to the integral reported in Gradshteyn and Ryzhik [33, p.550, 4.283 6]

$$(3.23) \qquad \int_0^1 \left[ \frac{1}{2} \frac{1+u}{1-u} + \frac{1}{\log u} - \log u \right] \frac{du}{\log u} = \frac{1}{2} \log(2\pi)$$

$\square$

We also have

$$D_k = \frac{1}{2\pi k} \int_0^{\infty} \frac{e^{-2\pi t} - e^{-2k\pi t}}{t} dt = \frac{\gamma + \log(2\pi k)}{2\pi k}$$

and we have the summation

$$\sum_{k=1}^{\infty} \frac{D_k}{k} = \frac{1}{2\pi} \sum_{k=1}^{\infty} \frac{1}{k^2} \int_0^{\infty} \frac{e^{-2\pi t} - e^{-2k\pi t}}{t} dt = \frac{1}{2\pi} \sum_{k=1}^{\infty} \frac{\gamma + \log(2k\pi)}{k^2}$$

We see that

$$\sum_{k=1}^{\infty} \frac{1}{k^2} \int_0^{\infty} \frac{e^{-2\pi t} - e^{-2k\pi t}}{t} dt = \int_0^{\infty} \frac{\varsigma(2)e^{-2\pi t} - Li_2(e^{-2\pi t})}{t} dt$$

in terms of the dilogarithm function. With the substitution $x = e^{-2\pi t}$ the integral becomes

$$= 4\pi^2 \int_0^1 [Li_2(x) - \varsigma(2)x] \frac{x}{\log x} dx$$

so that we have

$$(3.24) \qquad \int_0^1 [Li_2(x) - \varsigma(2)x] \frac{x}{\log x} dx = \sum_{k=1}^{\infty} \frac{\gamma + \log(2k\pi)}{(2k\pi)^2}$$

$\square$

Kummer's Fourier series expansion (3.12) may also be written as

$$\log \Gamma(x) = -\frac{1}{2}\log x - \frac{1}{2}\log \frac{\sin \pi x}{\pi x} + \sum_{n=1}^{\infty} \frac{(\gamma + \log 2\pi n)\sin 2\pi nx}{\pi n} \quad (0 < x < 1)$$

and hence we have

$$\frac{1}{2}\log \frac{\Gamma(x)}{\Gamma(1-x)} = \sum_{n=1}^{\infty} \frac{(\gamma + \log 2\pi n)\sin 2\pi nx}{\pi n}$$

We have $D_k = \int_0^1 \log \Gamma(x) \sin 2k\pi x\, dx = \frac{1}{2\pi k}\left(\gamma + \log 2k\pi\right)$ for $k \geq 1$ and hence we get

$$\sum_{k=1}^{\infty} \frac{D_k}{k} = \sum_{k=1}^{\infty} \frac{1}{2\pi k^2}\left(\gamma + \log 2k\pi\right)$$

$$= \frac{\left(\gamma + \log 2\pi\right)}{2\pi}\varsigma(2) - \frac{1}{2\pi}\varsigma'(2)$$

We also have

$$\sum_{k=1}^{\infty} \frac{D_k}{k} = \sum_{k=1}^{\infty} \frac{1}{k} \int_0^1 \log \Gamma(x) \sin 2k\pi x \, dx = \int_0^1 \log \Gamma(x) \sum_{k=1}^{\infty} \frac{\sin 2k\pi x}{k} \, dx$$

Now using the Fourier series which is valid for $x \in (0,1)$

$$\sum_{k=1}^{\infty} \frac{\sin 2k\pi x}{k} = \frac{\pi}{2}(1-2x)$$

we have

$$\sum_{k=1}^{\infty} \frac{D_k}{k} = \frac{\pi}{2} \int_0^1 \log \Gamma(x)(1-2x) \, dx$$

Using $\int_0^1 \log \Gamma(x) \, dx = \frac{1}{2} \log 2\pi$ we obtain

$$\frac{\pi}{4} \log 2\pi - \pi \int_0^1 x \log \Gamma(x) \, dx = \frac{(\gamma + \log 2\pi)}{2\pi} \varsigma(2) - \frac{1}{2\pi} \varsigma'(2)$$

Therefore we obtain

(3.25) $$\int_0^1 x \log \Gamma(x) \, dx = \frac{1}{6} \log 2\pi - \frac{\gamma}{12} + \frac{1}{2\pi^2} \varsigma'(2)$$

which we also derived in (2.38).

With reference to series of the form $\sum_{k=1}^{\infty} \frac{D_k}{k^p}$ and $\sum_{k=1}^{\infty} (-1)^k \frac{D_k}{k^p}$ we may easily obtain further identities using the Fourier series reported, for example, in [55, p.148].

We have $C_k = \int_0^1 \log \Gamma(x) \cos 2k\pi x \, dx = \frac{1}{4k}$ for $k \geq 1$ and hence we get

$$\sum_{k=1}^{\infty} (-1)^{k+1} \frac{C_k}{k} = \frac{1}{4} \varsigma_a(2)$$

We also have

$$\sum_{k=1}^{\infty}(-1)^{k+1}\frac{C_k}{k}=\sum_{k=1}^{\infty}\frac{1}{k}\int_0^1\log\Gamma(x)(-1)^{k+1}\cos 2k\pi x\,dx=\int_0^1\log\Gamma(x)\sum_{k=1}^{\infty}(-1)^{k+1}\frac{\cos 2k\pi x}{k}\,dx$$

Now using [55, p.148] we have

$$\sum_{k=1}^{\infty}(-1)^{k+1}\frac{\cos 2k\pi x}{k}=\log[2\cos\pi x]$$

and accordingly

$$\sum_{k=1}^{\infty}(-1)^{k+1}\frac{C_k}{k}=\int_0^1\log\Gamma(x)\log[2\cos\pi x]\,dx=\frac{1}{4}\varsigma_a(2)$$

We therefore have

(3.26) $$\int_0^1\log\Gamma(x)\log\cos\pi x\,dx=-\frac{1}{2}\log 2\log 2\pi+\frac{\pi^2}{48}$$

It is noted in [30, Eq. (7.4)] that

(3.27) $$\int_0^1\log\Gamma(x)\log\sin\pi x\,dx=-\frac{1}{2}\log 2\log 2\pi-\frac{\pi^2}{24}$$

which was derived in a different manner in [26].

Addition of the two integrals results in

(3.28) $$\int_0^1\log\Gamma(x)\log\sin 2\pi x\,dx=-\frac{1}{2}\log 2\log 2\pi-\frac{\pi^2}{48}$$

From equation (6.92a) of [21] we have

$$\int_0^1\log\Gamma(x+1)\log\big[2\sin(\pi x)\big]\,dx=\frac{1}{2\pi}\sum_{n=1}^{\infty}\frac{si(2n\pi)}{n^2}$$

and we see that

$$\int_0^1\log\Gamma(x+1)\log\big[2\sin(\pi x)\big]\,dx=\int_0^1\big[\log x+\log\Gamma(x)\big]\big[\log 2+\log\sin(\pi x)\big]\,dx$$

$$= \log 2 \int_0^1 \log x \, dx + \int_0^1 \log x \log \sin(\pi x) dx + \log 2 \int_0^1 \log \Gamma(x) dx + \int_0^1 \log \Gamma(x) \log \sin(\pi x) \, dx$$

$$= -\log 2 + \int_0^1 \log x \log \sin(\pi x) dx + \frac{1}{2} \log 2 \log 2\pi - \frac{1}{2} \log 2 \log 2\pi - \frac{\pi^2}{24}$$

Hence we have

$$(3.29) \qquad \int_0^1 \log x \log \sin(\pi x) dx = \frac{1}{2\pi} \sum_{n=1}^{\infty} \frac{si(2n\pi)}{n^2} + \log 2 + \frac{\pi^2}{24}$$

and $\displaystyle\sum_{n=1}^{\infty} \frac{si(2n\pi)}{n^2}$ may be determined from (6.117j) in [21]

$$\frac{1}{2\pi^2} \sum_{n=1}^{\infty} \frac{Si(2n\pi)}{n^2} = \log A - \frac{1}{4}$$

so that we have

$$(3.29.1) \quad \int_0^1 \log x \log \sin(\pi x) dx = \pi \log(2A) - \frac{\pi}{4}$$

$\square$

Using Euler's reflection formula we see that

$$\int_0^1 \log \Gamma(x) \log \sin \pi x \, dx + \int_0^1 \log \Gamma(1-x) \log \sin \pi x \, dx$$

$$= \log \pi \int_0^1 \log \sin \pi x \, dx - \int_0^1 \log^2 \sin \pi x \, dx$$

and we have

$$\int_0^1 \log \Gamma(1-x) \log \sin \pi x \, dx = \int_0^1 \log \Gamma(u) \log \sin \pi u \, du$$

Therefore using (3.27) we obtain

$$(3.30) \quad \int_0^1 \log^2 \sin \pi x \, dx = \log 2 \log 2\pi + \frac{\pi^2}{12} - \log 2 \log \pi = \log^2 2 + \frac{\pi^2}{12}$$

which is derived in a different manner in (6.3).

$\square$

In a similar manner we multiply by $\log \Gamma(x)$ and integrate to obtain

$$\int_0^1 \log^2 \Gamma(x)\,dx + \int_0^1 \log \Gamma(1-x) \log \Gamma(x)\,dx$$

$$= \log \pi \int_0^1 \log \Gamma(x)\,dx - \int_0^1 \log \Gamma(x) \log \sin \pi x\,dx$$

Espinosa and Moll [30] proved that

$$(3.32) \quad \int_0^1 \log^2 \Gamma(x)\,dx = \frac{\gamma^2}{12} + \frac{\pi^2}{48} + \frac{1}{6}\gamma \log(2\pi) + \frac{1}{3}\log^2(2\pi) - [\gamma + \log(2\pi)]\frac{\varsigma'(2)}{\pi^2} + \frac{\varsigma''(2)}{2\pi^2}$$

and a different derivation is given in [26]. Therefore we obtain

$$\int_0^1 \log \Gamma(1-x) \log \Gamma(x)\,dx = \frac{1}{2}\log \pi \log(2\pi) + \frac{1}{2}\log 2 \log 2\pi + \frac{\pi^2}{24}$$

$$- \left[\frac{\gamma^2}{12} + \frac{\pi^2}{48} + \frac{1}{6}\gamma \log(2\pi) + \frac{1}{3}\log^2(2\pi) - [\gamma + \log(2\pi)]\frac{\varsigma'(2)}{\pi^2} + \frac{\varsigma''(2)}{2\pi^2}\right]$$

which simplifies to

$$(3.33)$$
$$\int_0^1 \log \Gamma(1-x) \log \Gamma(x)\,dx = -\left[\frac{\gamma^2}{12} - \frac{\pi^2}{48} + \frac{1}{6}\gamma \log(2\pi) - \frac{1}{6}\log^2(2\pi) - [\gamma + \log(2\pi)]\frac{\varsigma'(2)}{\pi^2} + \frac{\varsigma''(2)}{2\pi^2}\right]$$

We may also use Parseval's theorem to evaluate $\int_0^1 \log^2 \Gamma(x)\,dx$ in a very elementary manner. Referring to (3.8) and (3.10) we see that

$$\int_0^1 \log^2 \Gamma(x)\,dx = \frac{1}{2}\int_0^1 \log \Gamma(x)\,dx + 2\sum_{k=1}^{\infty}\left(\frac{1}{16k^2} + \frac{[\gamma + \log(2k\pi)]^2}{4\pi^2 k^2}\right)$$

$$= \frac{1}{4}\log(2\pi) + 2\sum_{k=1}^{\infty}\left(\frac{1}{16k^2} + \frac{\gamma^2 + \log^2(2\pi) + 2\gamma \log(2\pi) + 2[\gamma + \log(2\pi)]\log k + \log^2 k}{4\pi^2 k^2}\right)$$

$$= \frac{1}{4}\log(2\pi) + \frac{1}{8}\varsigma(2) + \frac{1}{12}[\gamma + \log(2\pi)]^2 - [\gamma + \log(2\pi)]\frac{\varsigma'(2)}{\pi^2} + \frac{\varsigma''(2)}{2\pi^2}$$

in agreement with (3.32) above.

Amdeberhan et al. [4] have recently given further consideration to integrals of the form $\int_0^1 \log^n \Gamma(x)\,dx$.

$\square$

Berndt [10] used Kummer's formula (3.12) to derive Euler's reflection formula in the following elementary manner. We have

$$\log \Gamma(x) = \frac{1}{2}\log \pi - \frac{1}{2}\log \sin \pi x + \sum_{n=1}^{\infty} \frac{(\gamma + \log 2\pi n)\sin 2\pi nx}{\pi n} \quad (0 < x < 1)$$

Letting $x \to 1 - x$ we get

$$\log \Gamma(1-x) = \frac{1}{2}\log \pi - \frac{1}{2}\log \sin \pi(1-x) + \sum_{n=1}^{\infty} \frac{(\gamma + \log 2\pi n)\sin 2\pi n(1-x)}{\pi n}$$

$$= \frac{1}{2}\log \pi - \frac{1}{2}\log \sin \pi x - \sum_{n=1}^{\infty} \frac{(\gamma + \log 2\pi n)\sin 2\pi nx}{\pi n}$$

and hence

$$\log \Gamma(x) + \log \Gamma(1-x) = \log \frac{\pi}{\sin \pi x}$$

We then obtain Euler's reflection formula $\Gamma(x)\Gamma(1-x) = \frac{\pi}{\sin \pi x}$.

Alternatively, we may subtract the above two equations to obtain the following formula which appears in Ramanujan's Notebooks [11, Vol.1, p.199]

$$\log \frac{\Gamma(x)}{\Gamma(1-x)} + [\gamma + \log(2\pi)](2x-1) = \frac{2}{\pi}\sum_{n=1}^{\infty} \frac{\log n}{n}\sin 2\pi nx$$

where we have used the Fourier series $\sum_{n=1}^{\infty} \frac{\sin 2\pi nx}{n\pi} = \frac{1}{2} - x$ which is valid for $0 < x < 1$.

Integration then results in

$$\int\limits_0^u \log \Gamma(x)dx - \int\limits_0^u \log \Gamma(1-x)dx + [\gamma + \log(2\pi)](u^2 - u)$$

$$= \frac{1}{\pi^2}\sum_{n=1}^\infty \frac{\log n}{n^2} - \frac{1}{\pi^2}\sum_{n=1}^\infty \frac{\log n}{n^2}\cos 2\pi nu$$

We have

$$\int\limits_0^u \log \Gamma(1-x)dx = -\int\limits_0^{-u} \log \Gamma(1+x)dx$$

and we note Alexeiewsky's theorem [52, p.32] (a further derivation of which is contained in (1.115) and also in equation (4.3.85) of [20])

$$\int\limits_0^u \log \Gamma(1+x)dx = \frac{1}{2}\big[\log(2\pi)-1\big]u - \frac{u^2}{2} + u\log \Gamma(1+u) - \log G(1+u)$$

in terms of the Barnes double gamma function $G(u)$ (see Appendix C).

Hence we have

$$\int\limits_0^u \log \Gamma(x)dx - \int\limits_0^u \log \Gamma(1-x)dx = \int\limits_0^u \log \Gamma(x)dx + \int\limits_0^{-u} \log \Gamma(1+x)dx$$

$$= \int\limits_0^u \log \Gamma(1+x)dx + \int\limits_0^{-u} \log \Gamma(1+x)dx - \int\limits_0^u \log x\, dx$$

We see that

$$\int\limits_0^u \log \Gamma(1+x)dx + \int\limits_0^{-u} \log \Gamma(1+x)dx$$

$$= -u^2 + u\log \Gamma(1+u) - u\log \Gamma(1-u) - \log G(1+u) - \log G(1-u)$$

and with $u = 1/2$ we get

$$\int\limits_0^{1/2} \log \Gamma(1+x)dx + \int\limits_0^{-1/2} \log \Gamma(1+x)dx$$

$$= -\frac{1}{4} + \frac{1}{2}\log\Gamma(3/2) - \frac{1}{2}\log\Gamma(1/2) - \log G(3/2) - \log G(1/2)$$

$$= -\frac{1}{4} - \frac{1}{2}\log 2 - \log\Gamma(1/2) - 2\log G(1/2)$$

Therefore we obtain

$$-\frac{1}{4} - \frac{1}{2}\log 2 - \log\Gamma(1/2) - 2\log G(1/2) - \frac{1}{4}[\gamma + \log(2\pi)] = \frac{1}{\pi^2}\sum_{n=1}^{\infty}\frac{\log n}{n^2} - \frac{1}{\pi^2}\sum_{n=1}^{\infty}(-1)^n\frac{\log n}{n^2}$$

$$= -\frac{1}{\pi^2}[\varsigma'(2) + \varsigma_a'(2)]$$

Since $\varsigma_a'(s) = \left(1 - 2^{1-s}\right)\varsigma'(s) + 2^{1-s}\log 2.\varsigma(s)$ we have

$$\varsigma_a'(2) = \frac{1}{2}\varsigma'(2) + \frac{1}{2}\log 2.\varsigma(2)$$

and this gives us

$$-\frac{1}{4} - \frac{1}{2}\log 2 - \log\Gamma(1/2) - 2\log G(1/2) - \frac{1}{4}[\gamma + \log(2\pi)] = -\frac{1}{\pi^2}\left[\frac{3}{2}\varsigma'(2) + \frac{1}{2}\log 2.\varsigma(2)\right]$$

$$= -\frac{3}{2\pi^2}\varsigma'(2) - \frac{1}{12}\log 2$$

The derivative of the Riemann functional equation gives us

$$\log A - \frac{1}{12}[\gamma + \log(2\pi)] = -\frac{1}{2\pi^2}\varsigma'(2)$$

and we then obtain

(3.34)   $\log G(1/2) = -\frac{3}{2}\log A - \frac{1}{4}\log\pi + \frac{1}{8} + \frac{1}{24}\log 2$

as originally determined by Barnes [9] in 1899.

In the following two sections we generalise the analysis to address integrals of the form $\int_0^1 \log \Gamma(x) e^{i\pi px} dx$.

## 4. The $\int_0^1 \log \Gamma(x) \cos p\pi x \, dx$ integral

In this part we first of all consider the integral $\int_0^1 \log \Gamma(x) \cos p\pi x \, dx$; the corresponding integral with $\sin p\pi x$ in the integrand is considered later in Section 5 of this paper.

Provided $a \neq b$ we readily determine that

$$\int \sin ax \cos bx \, dx = -\frac{1}{2}\left[\frac{\cos(a+b)x}{a+b} + \frac{\cos(a-b)x}{a-b}\right] + c$$

and hence we have

$$\int_0^1 \sin 2\pi nx \cos p\pi x \, dx = \frac{1}{2}\left[\frac{1}{(2n+p)\pi} + \frac{1}{(2n-p)\pi}\right] - \frac{1}{2}\left[\frac{\cos(2n+p)\pi}{(2n+p)\pi} + \frac{\cos(2n-p)\pi}{(2n-p)\pi}\right]$$

$$(4.1) \qquad = \frac{2n}{\pi}\frac{1 - \cos p\pi}{4n^2 - p^2}$$

As before, provided $a \neq b$ we readily determine that

$$\int \cos ax \cos bx \, dx = \frac{1}{2}\left[\frac{\sin(a+b)x}{a+b} + \frac{\sin(a-b)x}{a-b}\right] + c$$

and hence we have

$$\int_0^1 \cos 2\pi nx \cos p\pi x \, dx = \frac{1}{2}\left[\frac{\sin(2n+p)\pi}{(2n+p)\pi} + \frac{\sin(2n-p)\pi}{(2n-p)\pi}\right]$$

$$(4.2) \qquad = -\frac{p}{\pi}\frac{\sin p\pi}{4n^2 - p^2}$$

We now multiply Hurwitz's identity (3.1) by $\cos p\pi x$ and integrate to obtain

$$(4.3) \qquad \int_0^1 \varsigma(s,x) \cos p\pi x \, dx$$

$$= 2\Gamma(1-s)\left[ -\sin\left(\frac{\pi s}{2}\right) \sum_{n=1}^{\infty} \frac{1}{(2\pi n)^{1-s}} \frac{p}{\pi} \frac{\sin p\pi}{4n^2-p^2} + \cos\left(\frac{\pi s}{2}\right) \sum_{n=1}^{\infty} \frac{1}{(2\pi n)^{1-s}} \frac{2n}{\pi} \frac{1-\cos p\pi}{4n^2-p^2} \right]$$

Differentiation with respect to $s$ gives us

$$(4.4) \qquad \int_0^1 \varsigma'(s,x) \cos p\pi x \, dx$$

$$= 2\Gamma(1-s)\left[ -\sin\left(\frac{\pi s}{2}\right) \sum_{n=1}^{\infty} \frac{\log(2\pi n)}{(2\pi n)^{1-s}} \frac{p}{\pi} \frac{\sin p\pi}{4n^2-p^2} - \frac{\pi}{2}\cos\left(\frac{\pi s}{2}\right) \sum_{n=1}^{\infty} \frac{1}{(2\pi n)^{1-s}} \frac{p}{\pi} \frac{\sin p\pi}{4n^2-p^2} \right]$$

$$+ 2\Gamma(1-s)\left[ \cos\left(\frac{\pi s}{2}\right) \sum_{n=1}^{\infty} \frac{\log(2\pi n)}{(2\pi n)^{1-s}} \frac{2n}{\pi} \frac{1-\cos p\pi}{4n^2-p^2} - \frac{\pi}{2}\sin\left(\frac{\pi s}{2}\right) \sum_{n=1}^{\infty} \frac{1}{(2\pi n)^{1-s}} \frac{2n}{\pi} \frac{1-\cos p\pi}{4n^2-p^2} \right]$$

$$- 2\Gamma'(1-s)\left[ -\sin\left(\frac{\pi s}{2}\right) \sum_{n=1}^{\infty} \frac{1}{(2\pi n)^{1-s}} \frac{p}{\pi} \frac{\sin p\pi}{4n^2-p^2} + \cos\left(\frac{\pi s}{2}\right) \sum_{n=1}^{\infty} \frac{1}{(2\pi n)^{1-s}} \frac{2n}{\pi} \frac{1-\cos p\pi}{4n^2-p^2} \right]$$

and with $s=0$ we have using Lerch's identity (1.46.1)

$$(4.5) \qquad \int_0^1 \log\Gamma(x) \cos p\pi x \, dx$$

$$= \frac{1}{2}\log(2\pi)\frac{\sin p\pi}{p\pi} - \frac{p\sin p\pi}{2\pi} \sum_{n=1}^{\infty} \frac{1}{n} \frac{1}{4n^2-p^2} + \frac{2(1-\cos p\pi)}{\pi^2} \sum_{n=1}^{\infty} \frac{\gamma+\log(2\pi n)}{4n^2-p^2}$$

and, as expected, with $p=0$ we recover Raabe's integral (2.11).

Using (2.15) and (4.7) this may be expressed as

$$(4.5.1) \qquad \int_0^1 \log\Gamma(x) \cos p\pi x \, dx$$

$$= \frac{1}{2}\log(2\pi)\frac{\sin p\pi}{p\pi} + \frac{\sin p\pi}{4p\pi}\left[ \psi\left(1+\frac{p}{2}\right) + \psi\left(1-\frac{p}{2}\right) + 2\gamma \right]$$

$$+ \frac{2(1-\cos p\pi)[\gamma + \log(2\pi)]}{\pi^2}\left[\frac{1}{2p^2} - \frac{\pi}{4p}\cot\left(\frac{p\pi}{2}\right)\right] + \frac{2(1-\cos p\pi)}{\pi^2}\sum_{n=1}^{\infty}\frac{\log n}{4n^2 - p^2}$$

$$= \frac{\sin p\pi[\gamma + \log(2\pi)]}{2p\pi} + \frac{\sin p\pi}{4p\pi}\left[\psi\left(\frac{p}{2}\right) + \psi\left(-\frac{p}{2}\right)\right]$$

$$+ \frac{2(1-\cos p\pi)[\gamma + \log(2\pi)]}{\pi^2}\left[\frac{1}{2p^2} - \frac{\pi}{4p}\cot\left(\frac{p\pi}{2}\right)\right] + \frac{2(1-\cos p\pi)}{\pi^2}\sum_{n=1}^{\infty}\frac{\log n}{4n^2 - p^2}$$

and with $p = 1$ we have

$$(4.5.2) \qquad \int_0^1 \log\Gamma(x)\cos\pi x\,dx = \frac{2}{\pi^2}\left[\log(2\pi) + \gamma + 2\sum_{n=1}^{\infty}\frac{\log n}{4n^2 - 1}\right]$$

With $p = 2k + 1$ in (4.5) we obtain

$$(4.6) \qquad \int_0^1 \log\Gamma(x)\cos(2k+1)\pi x\,dx = \frac{4}{\pi^2}\sum_{n=1}^{\infty}\frac{\gamma + \log(2\pi n)}{4n^2 - (2k+1)^2}$$

We have the well known identity [8, p.345] which is easily derived via Fourier series (in fact, a different derivation of this is obtained as a by-product in Section 6 below)

$$(4.7) \qquad \pi\cot\alpha\pi = \frac{1}{\alpha} - 2\alpha\sum_{n=1}^{\infty}\frac{1}{n^2 - \alpha^2}$$

and with $2\alpha = 2k + 1$ we obtain [47, No. 5.1.25.4]

$$(4.8) \qquad \sum_{n=1}^{\infty}\frac{1}{4n^2 - (2k+1)^2} = \frac{1}{2(2k+1)^2}$$

which then gives us

$$(4.9) \qquad \int_0^1 \log\Gamma(x)\cos(2k+1)\pi x\,dx = \frac{2}{\pi^2}\left[\frac{\log(2\pi) + \gamma}{(2k+1)^2} + 2\sum_{n=1}^{\infty}\frac{\log n}{4n^2 - (2k+1)^2}\right]$$

as originally determined by Kölbig [40] using Kummer's Fourier series for $\log\Gamma(x)$.

As noted by Kölbig [40], this integral is incorrectly reported in Nielsen's book [45, p.203]. The result was corrected in the 5th edition of Gradshteyn and Ryzhik [33, p.650, 6.443 4].

With $p = 2k$ in (4.5), and applying L'Hôpital's rule to $\dfrac{p \sin p\pi}{2\pi} \sum_{n=1}^{\infty} \dfrac{1}{n} \dfrac{1}{4n^2 - p^2}$, we obtain (2.28) again.

$$\int_0^1 \log \Gamma(x) \cos 2k\pi x \, dx = \frac{1}{4k}$$

Letting $k = 0$ in (4.9) results in (4.5.2) again.

$\square$

Letting $p = 1/2$ in (4.5) results in

$$\int_0^1 \log \Gamma(x) \cos(\pi x / 2) \, dx = \frac{1}{\pi} \log(2\pi) - \frac{1}{\pi} \sum_{n=1}^{\infty} \frac{1}{n} \frac{1}{16n^2 - 1} + \frac{8[\gamma + \log(2\pi)]}{\pi^2} \sum_{n=1}^{\infty} \frac{1}{16n^2 - 1}$$

$$+ \frac{8}{\pi^2} \sum_{n=1}^{\infty} \frac{\log n}{16n^2 - 1}$$

With $\alpha = 1/4$ in (4.7) we obtain

$$\sum_{n=1}^{\infty} \frac{1}{16n^2 - 1} = \frac{4 - \pi}{8}$$

and using (2.14)

$$\sum_{n=1}^{\infty} \frac{1}{n} \frac{1}{16n^2 - 1} = 3\log 2 - 2$$

we obtain

(4.10) $$\int_0^1 \log \Gamma(x) \cos(\pi x / 2) \, dx = \frac{1}{\pi} \log(2\pi) - \frac{3\log 2 - 2}{\pi} + \frac{[\gamma + \log(2\pi)](4 - \pi)}{\pi^2}$$

$$+ \frac{8}{\pi^2} \sum_{n=1}^{\infty} \frac{\log n}{16n^2 - 1}$$

$\square$

Using integration by parts we have

$$\int_0^1 \log\Gamma(x)\cos(2k+1)\pi x\, dx = \frac{\log\Gamma(x)\sin(2k+1)\pi x}{(2k+1)\pi}\Bigg|_0^1 - \frac{1}{(2k+1)\pi}\int_0^1 \psi(x)\sin(2k+1)\pi x\, dx$$

We see that

$$\log\Gamma(x)\sin(2k+1)\pi x = [\log\Gamma(1+x) - \log x]\sin(2k+1)\pi x$$

$$= [x\log\Gamma(1+x) - x\log x]\frac{\sin(2k+1)\pi x}{x}$$

and it is therefore seen that the integrated part vanishes.

Hence we have

(4.11)     $$\int_0^1 \log\Gamma(x)\cos(2k+1)\pi x\, dx = -\frac{1}{(2k+1)\pi}\int_0^1 \psi(x)\sin(2k+1)\pi x\, dx$$

and, as pointed out by Kölbig [40], letting $k=0$ we obtain

(4.12)     $$\int_0^1 \psi(x)\sin\pi x\, dx = -\frac{2}{\pi}\left[\log(2\pi) + \gamma + 2\sum_{n=1}^{\infty}\frac{\log n}{4n^2-1}\right]$$

Kölbig [40] states that the infinite series in (4.12) "does not seem to be expressible in terms of well-known functions".

Letting $k=0$ in (4.11) gives us

$$\int_0^1 \log\Gamma(x)\cos\pi x\, dx = -\frac{1}{\pi}\int_0^1 \psi(x)\sin\pi x\, dx$$

An alternative derivation of (4.12) is shown below.

We recall Lerch's trigonometric series expansion for the digamma function for $0 < x < 1$ (see for example Gronwall's paper [34, p.105] and Nielsen's book [45, p.204])

(4.13)     $$\psi(x)\sin\pi x + \frac{\pi}{2}\cos\pi x + (\gamma + \log 2\pi)\sin\pi x = -\sum_{n=1}^{\infty}\sin(2n+1)\pi x \cdot \log\frac{n+1}{n}$$

Integration of (4.13) gives us

(4.14)     $$\int_0^1 \psi(x)\sin\pi x\, dx = -\frac{2}{\pi}\left[\gamma + \log(2\pi) + \sum_{n=1}^{\infty}\frac{1}{2n+1}\log\frac{n+1}{n}\right]$$

We have

$$2\sum_{n=1}^{\infty}\frac{\log n}{4n^2-1}=\sum_{n=1}^{\infty}\left[\frac{1}{2n-1}-\frac{1}{2n+1}\right]\log n$$

and consider the finite sum

$$\sum_{n=1}^{N}\frac{\log n}{2n-1}=\sum_{m=0}^{N-1}\frac{\log(m+1)}{2m+1}$$

$$=\sum_{n=1}^{N}\frac{\log(n+1)}{2n+1}-\frac{\log(N+1)}{2N+1}$$

Hence we have

$$\sum_{n=1}^{N}\left[\frac{1}{2n-1}-\frac{1}{2n+1}\right]\log n=\sum_{n=1}^{N}\left[\frac{\log(n+1)}{2n+1}-\frac{\log n}{2n+1}\right]-\frac{\log(N+1)}{2N+1}$$

Therefore, as $N\rightarrow\infty$ we see that

(4.15)  $$2\sum_{n=1}^{\infty}\frac{\log n}{4n^2-1}=\sum_{n=1}^{\infty}\frac{1}{2n+1}\log\frac{n+1}{n}$$

and we have another proof of (4.12).

□

In passing, letting $x=1/2$ in (4.13) we obtain

$$\psi(1/2)+\gamma+\log(2\pi)=\sum_{n=1}^{\infty}(-1)^{n+1}\log\frac{n+1}{n}$$

and, since [52, p.20] $\psi(1/2)=-\gamma-2\log 2$, we have

$$\sum_{n=1}^{\infty}(-1)^{n+1}\log\frac{n+1}{n}=\log\frac{\pi}{2}$$

as previously reported by Sondow [51].

□

Integration of (4.13) gives us

$$\int_{0}^{1/2}\psi(x)\sin\pi x\,dx=\frac{1}{2}-\frac{1}{\pi}[\gamma+\log(2\pi)]-\frac{1}{\pi}\sum_{n=1}^{\infty}\frac{1}{2n+1}\log\frac{n+1}{n}$$

We may also write (4.13) as

$$(4.16) \qquad \psi(x) + \frac{\pi}{2}\cot\pi x + (\gamma + \log 2\pi) = -\sum_{n=1}^{\infty}\frac{\sin(2n+1)\pi x}{\sin\pi x}\cdot\log\frac{n+1}{n}$$

The above formula suggests that integration may be fruitful employing the following representation

$$\frac{\sin(2n+1)\pi x}{\sin\pi x} = 1 + 2\sum_{k=1}^{n}\cos 2k\pi x$$

The integral $\int x^p\psi(x)$ may also produce interesting results.

We multiply (4.13) by $\log\Gamma(x)$ and integrate to obtain

$$\int_0^1\log\Gamma(x)\psi(x)\sin\pi x\,dx + \frac{\pi}{2}\int_0^1\log\Gamma(x)\cos\pi x\,dx + (\gamma + \log 2\pi)\int_0^1\log\Gamma(x)\sin\pi x\,dx$$

$$= -\sum_{n=1}^{\infty}\int_0^1\log\Gamma(x)\sin(2n+1)\pi x\cdot\log\frac{n+1}{n}\,dx$$

Integration by parts gives us

$$\int_0^1\log\Gamma(x)\psi(x)\sin\pi x\,dx = \frac{1}{2}\log^2\Gamma(x)\sin\pi x\Big|_0^1 - \frac{\pi}{2}\int_0^1\log^2\Gamma(x)\cos\pi x\,dx$$

and with regard to the integrated part we see that

$$\lim_{x\to 0}\log^2\Gamma(x)\sin\pi x = \lim_{x\to 0}x\log^2\Gamma(x)\frac{\sin\pi x}{x}$$

We have

$$\lim_{x\to 0}x\log^2\Gamma(x) = \lim_{x\to 0}\frac{\log^2\Gamma(x)}{1/x}$$

$$= \lim_{x\to 0}\frac{2\log\Gamma(x)\psi(x)}{-1/x^2}$$

$$= -2\lim_{x \to 0} x \log \Gamma(x) \lim_{x \to 0} x\psi(x) = 0$$

and we thereby obtain

$$\int_0^1 \log \Gamma(x)\psi(x)\sin \pi x\, dx = -\frac{\pi}{2}\int_0^1 \log^2 \Gamma(x)\cos \pi x\, dx$$

$$-\frac{\pi}{2}\int_0^1 \log^2 \Gamma(x)\cos \pi x\, dx + \frac{1}{\pi}\left[\log(2\pi) + \gamma + 2\sum_{n=1}^{\infty}\frac{\log n}{4n^2 - 1}\right] + \frac{\gamma + \log 2\pi}{\pi}\left[\log \frac{\pi}{2} + 1\right]$$

$$= -\sum_{n=1}^{\infty}\frac{1}{(2n+1)\pi}\left[\log\left(\frac{\pi}{2}\right) + \frac{1}{2n+1} + 2\sum_{j=0}^{n-1}\frac{1}{2j+1}\right]\log\frac{n+1}{n}$$

Hence we have

$$(4.17)\ \frac{\pi^2}{2}\int_0^1 \log^2 \Gamma(x)\cos \pi x\, dx = \log(2\pi) + \gamma + \left[\log\frac{\pi}{2} + 1\right]\sum_{n=1}^{\infty}\frac{1}{2n+1}\log\frac{n+1}{n}$$

$$+ (\gamma + \log 2\pi)\left[\log\frac{\pi}{2} + 1\right] + \sum_{n=1}^{\infty}\frac{1}{2n+1}\left[\frac{1}{2n+1} + 2\sum_{j=0}^{n-1}\frac{1}{2j+1}\right]\log\frac{n+1}{n}$$

where we have used (5.6) and (5.7).

$\square$

With reference to (4.6) we make the summation

$$\int_0^1 \log \Gamma(x)\sum_{k=0}^{\infty}\frac{\cos(2k+1)\pi x}{(2k+1)^2}\, dx = \frac{2[\log(2\pi) + \gamma]}{\pi^2}\sum_{k=0}^{\infty}\frac{1}{(2k+1)^4}$$

$$+ \frac{4}{\pi^2}\sum_{k=0}^{\infty}\frac{1}{(2k+1)^2}\sum_{n=1}^{\infty}\frac{\log n}{4n^2 - (2k+1)^2}$$

where we have assumed that the interchange of the order of integration and summation is valid.

We have the Fourier series [55, p.149] which is valid for $0 \le x \le 1$

$$\sum_{k=0}^{\infty}\frac{\cos(2k+1)\pi x}{(2k+1)^2} = \frac{\pi^2}{8}(1 - 2x)$$

so that

$$\int_0^1 \log \Gamma(x) \sum_{k=0}^{\infty} \frac{\cos(2k+1)\pi x}{(2k+1)^2} \, dx = \frac{\pi^2}{8} \int_0^1 (1-2x) \log \Gamma(x) \, dx$$

Using (C.2) this becomes

$$= \frac{\pi^2}{8} \left[ \frac{3}{4} \log 2\pi - \log A \right]$$

and hence we determine that

$$(4.18) \quad \frac{4}{\pi^2} \sum_{k=0}^{\infty} \frac{1}{(2k+1)^2} \sum_{n=1}^{\infty} \frac{\log n}{4n^2 - (2k+1)^2} = \frac{\pi^2}{8} \left[ \frac{3}{4} \log(2\pi) - \log A \right] - 48[\log(2\pi) + \gamma]\pi^2$$

Similarly, one could also consider the representation [55, p.149]

$$\sum_{k=0}^{\infty} \frac{\cos(2k+1)\pi x}{2k+1} = -\frac{1}{2} \log \tan \frac{\pi x}{2}$$

but it seems that this may give rise to convergence issues.

$\square$

Differentiating (4.5) gives us

$$(4.19) \quad -\pi \int_0^1 x \log \Gamma(x) \sin p\pi x \, dx$$

$$= \frac{1}{2} \log(2\pi) \frac{p\pi \cos p\pi - \sin p\pi}{p^2 \pi} - \frac{p \sin p\pi}{2\pi} \sum_{n=1}^{\infty} \frac{1}{n} \frac{2p}{(4n^2 - p^2)^2} - \frac{p\pi \cos p\pi + \sin p\pi}{2\pi} \sum_{n=1}^{\infty} \frac{1}{n} \frac{1}{4n^2 - p^2}$$

$$+ \frac{4p(1 - \cos p\pi)}{\pi^2} \sum_{n=1}^{\infty} \frac{\gamma + \log(2\pi n)}{(4n^2 - p^2)^2} + \frac{2 \sin p\pi}{\pi} \sum_{n=1}^{\infty} \frac{\gamma + \log(2\pi n)}{4n^2 - p^2}$$

and with $p = 1$ we have

$$-\pi \int_0^1 x \log \Gamma(x) \sin \pi x \, dx$$

$$= -\frac{1}{2} \log(2\pi) + \frac{1}{2} \sum_{n=1}^{\infty} \frac{1}{n} \frac{1}{4n^2 - 1} + \frac{8}{\pi^2} \sum_{n=1}^{\infty} \frac{\gamma + \log(2\pi n)}{(4n^2 - 1)^2}$$

$$= -\frac{1}{2}\log(2\pi) + \frac{1}{2}\sum_{n=1}^{\infty}\frac{1}{n}\frac{1}{4n^2-1} + \frac{8[\gamma+\log(2\pi)]}{\pi^2}\sum_{n=1}^{\infty}\frac{1}{(4n^2-1)^2} + \frac{8}{\pi^2}\sum_{n=1}^{\infty}\frac{\log n}{(4n^2-1)^2}$$

Differentiating (4.7) results in

$$\pi^2\operatorname{cosec}^2\alpha\pi = \frac{1}{\alpha^2} + 4\alpha^2\sum_{n=1}^{\infty}\frac{1}{(n^2-\alpha^2)^2} + 2\sum_{n=1}^{\infty}\frac{1}{n^2-\alpha^2}$$

and with $\alpha = 1/2$ we have

$$\pi^2 = 4 + 16\sum_{n=1}^{\infty}\frac{1}{(4n^2-1)^2} + 8\sum_{n=1}^{\infty}\frac{1}{4n^2-1}$$

$$= 4 + 16\sum_{n=1}^{\infty}\frac{1}{(4n^2-1)^2} + 4$$

so that

$$\sum_{n=1}^{\infty}\frac{1}{(4n^2-1)^2} = \frac{\pi^2-8}{16}$$

as reported in [33, p.9].

Therefore using (2.15) we obtain

(4.20)

$$\pi\int_0^1 x\log\Gamma(x)\sin\pi x\,dx = \frac{1}{2}\log(2\pi) - \frac{1}{2}(2\log 2 - 1) - \frac{[\gamma+\log(2\pi)]}{\pi^2}\frac{\pi^2-8}{2} - \frac{8}{\pi^2}\sum_{n=1}^{\infty}\frac{\log n}{(4n^2-1)^2}$$

$\square$

Applying Lerch's identity [10]

$$\varsigma'(0,x) = \log\Gamma(x) - \frac{1}{2}\log(2\pi)$$

we see that

$$\frac{\partial}{\partial s}\int_0^1 \varsigma(s,x)\varsigma(p,x)\cos q\pi x\,dx\bigg|_{s=0} = \int_0^1 \varsigma(p,x)\log\Gamma(x)\cos q\pi x\,dx - \frac{1}{2}\log(2\pi)\int_0^1 \varsigma(p,x)\cos q\pi x\,dx$$

With $p = -n$ where $n \geq 0$ is a positive integer we have [6, p.264] in terms of the Bernoulli polynomials

$$\varsigma(-n, x) = -\frac{B_{n+1}(x)}{n+1}$$

so that

$$(n+1)\frac{\partial}{\partial s}\int_0^1 \varsigma(s,x)\varsigma(-n,x)\cos q\pi x\,dx\bigg|_{s=0} = \frac{1}{2}\log(2\pi)\int_0^1 B_{n+1}(x)\cos q\pi x\,dx$$

$$-\int_0^1 B_{n+1}(x)\log\Gamma(x)\cos q\pi x\,dx$$

The integral on the left-hand side may then be evaluated by using Hurwitz's formula (3.1).

$\square$

We now let $p \to ip$ where $i = \sqrt{-1}$ and hence we have

$$\int_0^1 \log\Gamma(x)\cos ip\pi x\,dx = \int_0^1 \log\Gamma(x)\cosh p\pi x\,dx$$

$$\int_0^1 \log\Gamma(x)\sin ip\pi x\,dx = i\int_0^1 \log\Gamma(x)\sinh p\pi x\,dx$$

With $\alpha \to i\alpha$ in (4.7) we see that

$$\pi\coth\alpha\pi = \frac{1}{\alpha} + 2\alpha\sum_{n=1}^{\infty}\frac{1}{n^2+\alpha^2}$$

(4.21) $\quad \int_0^1 \log\Gamma(x)\cosh p\pi x\,dx$

$$= \frac{1}{2}\log(2\pi)\frac{\sinh p\pi}{p\pi} + \frac{p\sinh p\pi}{2\pi}\sum_{n=1}^{\infty}\frac{1}{n}\frac{1}{4n^2+p^2} + \frac{2(1-\cosh p\pi)}{\pi^2}\sum_{n=1}^{\infty}\frac{\gamma+\log(2\pi n)}{4n^2+p^2}$$

Espinosa and Moll [30] showed that

$$\int_0^1 e^{tx}\varsigma(s,x)\,dx = 2(e^t-1)\Gamma(1-s)(2\pi)^{s-2}\sum_{n=0}^{\infty}(-1)^n\left(\frac{t}{2\pi}\right)^n\varsigma(n+2-s)\cos\left(\frac{\pi}{2}[s-n]\right)$$

and differentiating this with respect to $s$ will give us integrals involving hyperbolic functions.

## 5. The $\int_0^1 \log \Gamma(x) \sin p\pi x\, dx$ integral

Provided $a \neq b$ we readily determine that

$$\int \sin ax \sin bx\, dx = \frac{1}{2}\left[\frac{\sin(a-b)x}{a-b} - \frac{\sin(a+b)x}{a+b}\right] + c$$

and hence we have

$$\int_0^1 \sin 2\pi nx \sin p\pi x\, dx = \frac{1}{2}\left[\frac{\sin(2n-p)\pi}{(2n-p)\pi} - \frac{\sin(2n+p)\pi}{(2n+p)\pi}\right]$$

$$(5.1) \qquad\qquad = -\frac{2n \sin p\pi}{\pi}\frac{1}{4n^2 - p^2}$$

Similarly we have

$$\int \cos ax \sin bx\, dx = \frac{1}{2}\left[\frac{\cos(a-b)x}{a-b} - \frac{\cos(a+b)x}{a+b}\right] + c$$

which gives us

$$\int_0^1 \cos 2\pi nx \sin p\pi x\, dx = \frac{1}{2}\left[\frac{\cos(2n-p)\pi}{(2n-p)\pi} - \frac{\cos(2n+p)\pi}{(2n+p)\pi}\right] - \frac{1}{2}\left[\frac{1}{(2n-p)\pi} - \frac{1}{(2n+p)\pi}\right]$$

$$(5.2) \qquad\qquad = -\frac{p}{\pi}\frac{1-\cos p\pi}{4n^2 - p^2}$$

We now multiply (3.1) by $\sin p\pi x$ and integrate to get

$$(5.3) \quad \int_0^1 \varsigma(s,x) \sin p\pi x\, dx$$

$$= -2\Gamma(1-s)\left[\sin\left(\frac{\pi s}{2}\right)\sum_{n=1}^{\infty}\frac{1}{(2\pi n)^{1-s}}\frac{p}{\pi}\frac{1-\cos p\pi}{4n^2 - p^2} + \cos\left(\frac{\pi s}{2}\right)\sum_{n=1}^{\infty}\frac{1}{(2\pi n)^{1-s}}\frac{2n \sin p\pi}{\pi}\frac{1}{4n^2 - p^2}\right]$$

Differentiation with respect to $s$ gives us

(5.4) $\quad \int_0^1 \varsigma'(s,x) \sin p\pi x \, dx$

$$= -2\Gamma(1-s)\left[\sin\left(\frac{\pi s}{2}\right)\sum_{n=1}^{\infty}\frac{\log(2\pi n)}{(2\pi n)^{1-s}}\frac{p}{\pi}\frac{1-\cos p\pi}{4n^2-p^2}+\frac{\pi}{2}\cos\left(\frac{\pi s}{2}\right)\sum_{n=1}^{\infty}\frac{1}{(2\pi n)^{1-s}}\frac{p}{\pi}\frac{1-\cos p\pi}{4n^2-p^2}\right]$$

$$-2\Gamma(1-s)\left[\cos\left(\frac{\pi s}{2}\right)\sum_{n=1}^{\infty}\frac{\log(2\pi n)}{(2\pi n)^{1-s}}\frac{2n\sin p\pi}{\pi}\frac{1}{4n^2-p^2}-\frac{\pi}{2}\sin\left(\frac{\pi s}{2}\right)\sum_{n=1}^{\infty}\frac{1}{(2\pi n)^{1-s}}\frac{2n\sin p\pi}{\pi}\frac{1}{4n^2-p^2}\right]$$

$$+2\Gamma'(1-s)\left[\sin\left(\frac{\pi s}{2}\right)\sum_{n=1}^{\infty}\frac{1}{(2\pi n)^{1-s}}\frac{p}{\pi}\frac{1-\cos p\pi}{4n^2-p^2}+\cos\left(\frac{\pi s}{2}\right)\sum_{n=1}^{\infty}\frac{1}{(2\pi n)^{1-s}}\frac{2n\sin p\pi}{\pi}\frac{1}{4n^2-p^2}\right]$$

and with $s=0$ we have

(5.5) $\quad \int_0^1 \log\Gamma(x)\sin p\pi x \, dx$

$$= \frac{1}{2}\log(2\pi)\frac{1-\cos p\pi}{p\pi}-\frac{p(1-\cos p\pi)}{2\pi}\sum_{n=1}^{\infty}\frac{1}{n}\frac{1}{4n^2-p^2}-\frac{2\sin p\pi}{\pi^2}\sum_{n=1}^{\infty}\frac{\gamma+\log(2\pi n)}{4n^2-p^2}$$

With a little algebra and using (2.15) and (4.7) this may be expressed as

(5.5.1)

$$\int_0^1 \log\Gamma(x)\sin p\pi x \, dx = \frac{(1-\cos p\pi)[\gamma+\log(2\pi)]}{2p\pi}+\frac{(1-\cos p\pi)}{4p\pi}\left[\psi\left(1+\frac{p}{2}\right)+\psi\left(1-\frac{p}{2}\right)\right]$$

$$-\frac{2\sin p\pi[\gamma+\log(2\pi)]}{\pi^2}\left[\frac{1}{2p^2}-\frac{\pi}{4p}\cot\left(\frac{p\pi}{2}\right)\right]-\frac{2\sin p\pi}{\pi^2}\sum_{n=1}^{\infty}\frac{\log n}{4n^2-p^2}$$

With $p=1$ we have

$$\int_0^1 \log\Gamma(x)\sin \pi x \, dx = \frac{[\gamma+\log(2\pi)]}{\pi}+\frac{1}{2\pi}\left[\psi\left(1+\frac{1}{2}\right)+\psi\left(\frac{1}{2}\right)\right]$$

or

$$\int_0^1 \log \Gamma(x) \sin \pi x \, dx = \frac{1}{\pi}\left[ \log \frac{\pi}{2} + 1 \right]$$

With $p = 1/2$ we have

$$\int_0^1 \log \Gamma(x) \sin(\pi x / 2) \, dx = \frac{\gamma + \log(2\pi)}{\pi} + \frac{1}{2\pi}\left[ \psi\left(1 + \frac{1}{4}\right) + \psi\left(\frac{3}{4}\right) \right]$$

$$- \frac{2[\gamma + \log(2\pi)]}{\pi^2}\left[ 2 - \frac{\pi}{2} \right] - \frac{8}{\pi^2}\sum_{n=1}^{\infty}\frac{\log n}{16n^2 - 1}$$

$$= \frac{\gamma + \log(2\pi)}{\pi} + \frac{1}{\pi}\left[ 2 - \gamma - 3\log 2 \right]$$

$$- \frac{2[\gamma + \log(2\pi)]}{\pi^2}\left[ 2 - \frac{\pi}{2} \right] - \frac{8}{\pi^2}\sum_{n=1}^{\infty}\frac{\log n}{16n^2 - 1}$$

(5.5.1)    $$\int_0^1 \log \Gamma(x) \sin(\pi x / 2) \, dx = \frac{1}{\pi}\left[ 2 - \gamma - 3\log 2 \right] - \frac{4[\gamma + \log(2\pi)]}{\pi^2} - \frac{8}{\pi^2}\sum_{n=1}^{\infty}\frac{\log n}{16n^2 - 1}$$

With $p = 2k + 1$ in (5.5) we obtain

$$\int_0^1 \log \Gamma(x) \sin(2k+1)\pi x \, dx = \frac{1}{(2k+1)\pi}\log(2\pi) - \frac{(2k+1)}{\pi}\sum_{n=1}^{\infty}\frac{1}{n}\frac{1}{4n^2 - (2k+1)^2}$$

and using (2.21) we have the well-known integral [33, p.650, 6.443.2]

(5.6)    $$\int_0^1 \log \Gamma(x) \sin(2k+1)\pi x \, dx = \frac{1}{(2k+1)\pi}\left[ \log\left(\frac{\pi}{2}\right) + \frac{1}{2k+1} + 2\sum_{j=0}^{k-1}\frac{1}{2j+1} \right]$$

As noted by Kölbig [40], this integral is also incorrectly reported in Nielsen's book [45, p.203]. The result was corrected in the 5th edition of Gradshteyn and Ryzhik [33, p.650, 6.443 2]). A different proof of (5.6) is contained in [23]

We obtain with $k = 0$

(5.7)    $$\int_0^1 \log \Gamma(x) \sin \pi x \, dx = \frac{1}{\pi}\left[ \log\left(\frac{\pi}{2}\right) + 1 \right]$$

Letting $p \to 2k$ in (5.5) and using (2.21.1) and (2.22.1) we easily obtain another derivation of (1.5)

$$\int_0^1 \log \Gamma(x) \sin 2k\pi x \, dx = \frac{\gamma + \log(2\pi k)}{2\pi k}$$

$\square$

We have

$$I = \int_0^1 \log \Gamma(x) \sin(2k+1)\pi x \, dx = \int_0^1 \log \Gamma(1-x) \sin(2k+1)\pi x \, dx$$

so that using Euler's reflection formula we have

$$\int_0^1 [\log \pi - \log \sin \pi x] \sin(2k+1)\pi x \, dx = 2I$$

Using (5.6) we then obtain

$$(5.8) \qquad \int_0^1 \log \sin \pi x \sin(2k+1)\pi x \, dx = \frac{2}{(2k+1)\pi} \left[ \log 2 - \frac{1}{2k+1} - 2\sum_{j=0}^{k-1} \frac{1}{2j+1} \right]$$

which is reported in [33, p.577]. A different derivation is given in Section 6.

$\square$

Using integration by parts we have

$$\int_a^1 \log \Gamma(x) \sin 2\pi x \, dx = -\frac{\log \Gamma(x) \cos 2\pi x}{2\pi} \bigg|_a^1 + \frac{1}{2\pi} \int_a^1 \psi(x) \cos 2\pi x \, dx$$

$$= \frac{\log \Gamma(a)[\cos 2\pi a - 1]}{2\pi} + \frac{1}{2\pi} \int_a^1 \psi(x)[\cos 2\pi x - 1] \, dx$$

We see that

$$\log \Gamma(a)[\cos 2\pi a - 1] = [\log \Gamma(a+1) - \log a][\cos 2\pi a - 1]$$

$$= \log \Gamma(a+1)[\cos 2\pi a - 1] - a \log a \left[ \frac{\cos 2\pi a - 1}{a} \right]$$

and taking the limit as $a \to 0$ we find that

$$\int_0^1 \log \Gamma(x) \sin 2\pi x \, dx = \frac{1}{2\pi} \int_0^1 \psi(x)[\cos 2\pi x - 1] \, dx$$

$$= -\frac{1}{\pi} \int_0^1 \psi(x) \sin^2 2\pi x \, dx$$

Since from (1.5)

$$\int_0^1 \log \Gamma(x) \sin 2\pi x \, dx = \frac{\gamma + \log 2\pi}{2\pi}$$

we obtain

$$(5.9) \qquad \int_0^1 \psi(x) \sin^2 2\pi x \, dx = -\frac{\gamma + \log 2\pi}{2}$$

as noted in [33, p.652, 6.468] and [45, p.204].

□

We make the following summation by reference to (5.6)

$$\int_0^1 \log \Gamma(x) \sum_{k=0}^{\infty} \frac{\sin(2k+1)\pi x}{2k+1} \, dx = \sum_{k=0}^{\infty} \frac{1}{(2k+1)^2 \pi} \left[ \log\left(\frac{\pi}{2}\right) + \frac{1}{2k+1} + 2\sum_{j=0}^{k-1} \frac{1}{2j+1} \right]$$

and using the Fourier series [55, p.149] which is valid for $0 < x < 1$

$$\sum_{k=0}^{\infty} \frac{\sin(2k+1)\pi x}{2k+1} = \frac{\pi}{4}$$

together with the well known formula for the Riemann zeta function

$$\sum_{k=0}^{\infty} \frac{1}{(2k+1)^s} = (1 - 2^{-s})\varsigma(s)$$

we obtain

$$\frac{\pi^2}{4} \int_0^1 \log \Gamma(x) \, dx = \frac{3}{4} \varsigma(2) \log\left(\frac{\pi}{2}\right) + \frac{7}{8} \varsigma(3) + 2\sum_{k=0}^{\infty} \frac{1}{(2k+1)^2} \sum_{j=0}^{k-1} \frac{1}{2j+1}$$

Using Raabe's integral (2.11) we have the Euler sum

$$\sum_{k=0}^{\infty} \frac{1}{(2k+1)^2} \sum_{j=0}^{k-1} \frac{1}{2j+1} = \frac{\pi^2}{8} \log 2 - \frac{7}{16} \varsigma(3)$$

which was previously reported in [23].

We note that [52, p.20]

$$\psi\left(n+\frac{1}{2}\right) = -\gamma - 2\log 2 + 2\sum_{j=0}^{n-1} \frac{1}{2j+1} = \psi\left(\frac{1}{2}\right) + 2\sum_{j=0}^{n-1} \frac{1}{2j+1}$$

and we therefore have

$$\sum_{n=0}^{\infty} \frac{\psi\left(n+\frac{1}{2}\right)}{(2n+1)^2} = -\frac{1}{8}\left[\gamma\pi^2 + 7\varsigma(3)\right]$$

$\square$

Using Euler's reflection formula $\Gamma(x)\Gamma(1-x) = \dfrac{\pi}{\sin \pi x}$ we obtain from (3.12)

(5.10) $\qquad \dfrac{1}{2} \log \dfrac{\Gamma(x)}{\Gamma(1-x)} = \sum_{n=1}^{\infty} \dfrac{(\gamma + \log 2\pi n) \sin 2\pi n x}{\pi n}$

We now multiply this by $\cos p\pi x$ and integrate to obtain

$$\int_0^1 \log \Gamma(x) \cos p\pi x \, dx - \int_0^1 \log \Gamma(1-x) \cos p\pi x \, dx = 2\sum_{n=1}^{\infty} \frac{(\gamma + \log 2\pi n)}{\pi n} \int_0^1 \sin 2\pi n x \cos p\pi x \, dx$$

and using the integral (4.1) we have

$$\int_0^1 \log \Gamma(x) \cos p\pi x \, dx - \int_0^1 \log \Gamma(1-x) \cos p\pi x \, dx = \frac{4(1-\cos p\pi)}{\pi^2} \sum_{n=1}^{\infty} \frac{(\gamma + \log 2\pi n)}{4n^2 - p^2}$$

and using (2.9)

$$\int_0^1 \log \Gamma(x) \cos p\pi x \, dx + \int_0^1 \log \Gamma(1-x) \cos p\pi x \, dx = \frac{\sin p\pi}{p\pi} \log(2\pi) - \frac{p\sin p\pi}{\pi} \sum_{n=1}^{\infty} \frac{1}{n} \frac{1}{4n^2 - p^2}$$

we obtain (4.5) again

$$2\int_0^1 \log\Gamma(x)\cos p\pi x\,dx = \frac{\sin p\pi}{p\pi}\log(2\pi) - \frac{p\sin p\pi}{\pi}\sum_{n=1}^{\infty}\frac{1}{n}\frac{1}{4n^2-p^2} + \frac{4(1-\cos p\pi)}{\pi^2}\sum_{n=1}^{\infty}\frac{(\gamma+\log 2\pi n)}{4n^2-p^2}$$

## 6. Some integrals involving $\log\sin\pi x$

We note from (4.2) that

$$\int_0^1 \cos 2\pi n x\cos p\pi x\,dx = -\frac{p}{\pi}\frac{\sin p\pi}{4n^2-p^2}$$

and make the summation

$$\sum_{n=1}^{\infty}\frac{1}{n}\int_0^1 \cos 2\pi n x\cos p\pi x\,dx = -\frac{p}{\pi}\sum_{n=1}^{\infty}\frac{1}{n}\frac{\sin p\pi}{4n^2-p^2}$$

Assuming that it is valid to interchange the order of the integration and summation

$$\sum_{n=1}^{\infty}\frac{1}{n}\int_0^1 \cos 2\pi n x\cos p\pi x\,dx = \int_0^1 \sum_{n=1}^{\infty}\frac{\cos 2\pi n x}{n}\cos p\pi x\,dx$$

and using the Fourier series

$$\sum_{n=1}^{\infty}\frac{\cos 2\pi n x}{n} = -\log(2\sin\pi x)\qquad 0<x<1$$

we obtain

(6.1) $$\int_0^1 \log(2\sin\pi x)\cos p\pi x\,dx = \frac{p}{\pi}\sum_{n=1}^{\infty}\frac{1}{n}\frac{\sin p\pi}{4n^2-p^2}$$

Letting $p\to p+1$ we see that

$$\int_0^1 \log(2\sin\pi x)\cos p\pi x\,dx = \frac{(p+1)}{\pi}\sum_{n=1}^{\infty}\frac{1}{n}\frac{\sin p\pi}{4n^2-(p+1)^2}$$

and hence we obtain

$$\sum_{n=1}^{\infty}\frac{1}{n}\frac{1}{4n^2-p^2} = \frac{p+1}{p}\sum_{n=1}^{\infty}\frac{1}{n}\frac{1}{4n^2-(p+1)^2}$$

With $p = 2k$ in (6.1) we have

$$\int_0^1 \log(2\sin \pi x)\cos 2k\pi x\, dx = \lim_{p \to 2k} \frac{p}{\pi}\frac{1}{k}\frac{\sin p\pi}{4k^2 - p^2}$$

and we see that

$$\lim_{p \to 2k} \frac{p}{\pi}\frac{1}{k}\frac{\sin p\pi}{4k^2 - p^2} = 2\lim_{p \to 2k}\frac{\sin p\pi}{4k^2 - p^2}$$

L'Hôpital's rule gives us

$$\lim_{p \to 2k}\frac{\sin p\pi}{4k^2 - p^2} = \lim_{p \to 2k}\frac{\pi\cos p\pi}{-2p} = -\frac{\pi}{4k}$$

and hence we obtain (2.24)

$$\int_0^1 \log\sin \pi x\cos 2k\pi x\, dx = -\frac{1}{2k}\quad k > 0$$

We now make a further summation

$$\sum_{k=1}^{\infty}\frac{1}{k}\int_0^1 \log\sin \pi x\cos 2k\pi x\, dx = -\frac{1}{2}\varsigma(2)$$

and assuming that

$$\sum_{k=1}^{\infty}\frac{1}{k}\int_0^1 \log\sin \pi x\cos 2k\pi x\, dx = \int_0^1 \log\sin \pi x\sum_{k=1}^{\infty}\frac{\cos 2k\pi x}{k}\, dx$$

we have

$$= -\int_0^1 \log\sin \pi x\log(2\sin \pi x)\, dx$$

$$= -\int_0^1 \log^2\sin \pi x\, dx - \log 2\int_0^1 \log\sin \pi x\, dx$$

$$= -\int_0^1 \log^2\sin \pi x\, dx + \log^2 2$$

Accordingly we obtain the well known integral

(6.3)
$$\int_0^1 \log^2 \sin \pi x \, dx = \frac{1}{2}\varsigma(2) + \log^2 2$$

which is a problem posed by Bremekamp [13] more than fifty years ago in 1957.

There are many ways to evaluate this integral; for example, we could also have employed Parseval's theorem [8, p.343]

$$\frac{1}{\pi}\int_{-\pi}^{\pi} f(x)g(x)dx = \frac{1}{2}a_0\alpha_0 + \sum_{n=1}^{\infty}(a_n\alpha_n + b_n\beta_n)$$

to evaluate it by utilising the known Fourier coefficients (2.29), (2.31) and (2.32).

With $p = 2k+1$ in (6.1) we immediately determine (2.30)

$$\int_0^1 \log \sin \pi x \cos(2k+1)\pi x \, dx = 0$$

This may also be derived as follows. We designate $I$ as

$$I = \int_0^1 \log \Gamma(x)\cos(2k+1)\pi x \, dx$$

and with the substitution $x \to 1-x$ we see that

$$I = -\int_0^1 \log \Gamma(1-x)\cos(2k+1)\pi x \, dx$$

Accordingly we have

$$\int_0^1 [\log \Gamma(x) + \Gamma(1-x)]\cos(2k+1)\pi x \, dx = 0$$

or equivalently

$$\int_0^1 [\log \pi - \log \sin \pi x]\cos(2k+1)\pi x \, dx = 0$$

We then deduce that

$$\int_0^1 \log \sin \pi x \cos(2k+1)\pi x\, dx = 0$$

We have the particular integral with $k = 0$ which may be directly evaluated as

$$\int_0^1 \log \sin \pi x \cdot \cos \pi x\, dx = \frac{1}{\pi} \sin \pi x [\log \sin \pi x - 1]\Big|_0^1 = 0$$

because using L'Hôpital's rule we see that

$$\lim_{x \to 1} \sin \pi x \log \sin \pi x = \lim_{x \to 1} \frac{\log \sin \pi x}{1/\sin \pi x} = -\lim_{x \to 1} \frac{\pi \cot \pi x}{\pi \cos \pi x / \sin^2 \pi x} = 0$$

$\square$

We see from (6.1) that with $p \to 2p$

$$\frac{2p\pi}{\sin 2p\pi} \int_0^1 \log(2\sin \pi x)\cos 2p\pi x\, dx = p^2 \sum_{n=1}^{\infty} \frac{1}{n} \frac{1}{n^2 - p^2}$$

and it is easily seen that

$$I = \frac{2p\pi}{\sin 2p\pi} \int_0^1 \log(\sin \pi x)\cos 2p\pi x\, dx = \frac{2p\pi}{\sin 2p\pi} \int_0^1 \log(2\sin \pi x)\cos 2p\pi x\, dx - \log 2$$

$$= p^2 \sum_{n=1}^{\infty} \frac{1}{n} \frac{1}{n^2 - p^2} - \log 2$$

From (A.3) we have

$$\psi(1+p) + \psi(1-p) + 2\gamma = -2p^2 \sum_{n=1}^{\infty} \frac{1}{n} \frac{1}{n^2 - p^2}$$

and therefore we have

$$I = -\frac{1}{2}[\psi(1+p) + \psi(1-p) + 2\gamma] - \log 2$$

$$= -\frac{1}{2}[\psi(1+p) + \psi(1-p) + 2\gamma + 2\log 2]$$

$$= -\frac{1}{2}[\psi(p) + \psi(-p) + 2\gamma + 2\log 2]$$

This corresponds with the formula given by Dwilewicz and Mináč [28]

(6.4) $$\frac{2p\pi}{\sin 2p\pi} \int_0^1 \log(\sin \pi x) \cos 2p\pi x \, dx = -\frac{1}{2}[\psi(p) + \psi(-p) + 2\gamma + 2\log 2]$$

$\square$

We note from (5.2) that

$$\int_0^1 \cos 2\pi n x \sin p\pi x \, dx = -\frac{p}{\pi} \frac{1 - \cos p\pi}{4n^2 - p^2}$$

and following the above procedure we obtain

(6.5) $$\int_0^1 \log(2\sin \pi x) \sin p\pi x \, dx = \frac{p}{\pi} \sum_{n=1}^{\infty} \frac{1}{n} \frac{1 - \cos p\pi}{4n^2 - p^2}$$

As above, we may also show that (6.5) is equivalent to

(6.6) $$\frac{2p\pi}{1 - \cos 2p\pi} \int_0^1 \log(\sin \pi x) \sin 2p\pi x \, dx = -\frac{1}{2}[\psi(1 + p) + \psi(1 - p) + 2\gamma + 2\log 2]$$

which was derived in a very different manner by Dwilewicz and Mináč [28].

With $p = 2k$ we have

$$\int_0^1 \log(2\sin \pi x) \sin 2k\pi x \, dx = 2 \lim_{p \to 2k} \frac{1 - \cos p\pi}{4k^2 - p^2}$$

$$= -2 \lim_{p \to 2k} \frac{\sin p\pi}{2p}$$

Hence we obtain (2.32)

$$\int_0^1 \log \sin \pi x \sin 2k\pi x \, dx = 0 \qquad k \geq 0$$

With $p = 2k + 1$ in (6.5) we immediately obtain

$$\int_0^1 \log(2\sin\pi x)\sin(2k+1)\pi x\,dx = -\frac{2(2k+1)}{\pi}\sum_{n=1}^{\infty}\frac{1}{n}\frac{1}{4n^2-(2k+1)^2}$$

and using (2.18) this becomes (2.33)

$$\int_0^1 \log\sin\pi x\sin(2k+1)\pi x\,dx = \frac{2}{(2k+1)\pi}\left[\log 2 - \frac{1}{2k+1} - 2\sum_{j=0}^{k-1}\frac{1}{2j+1}\right]$$

$\square$

We now consider the limit of (6.6) as $p\to 0$; we have

$$\lim_{p\to 0}\frac{h(p)}{1-\cos 2p\pi} = -\frac{1}{2}\lim_{p\to 0}[\psi(1+p)+\psi(1-p)+2\gamma+2\log 2] = -\log 2$$

where

$$h(p) = 2\pi\int_0^1 \log(\sin\pi x)\,p\sin 2p\pi x\,dx$$

Applying L'Hôpital's rule twice, we easily find that

$$\lim_{p\to 0}\frac{h(p)}{1-\cos 2p\pi} = 2\int_0^1 x\log(\sin\pi x)\,dx$$

so that we have

(6.7) $\qquad \displaystyle\int_0^1 x\log(\sin\pi x)\,dx = -\frac{1}{2}\log 2$

in agreement with [30].

$\square$

We note from (5.1) that

$$\int_0^1 \sin 2\pi nx\sin p\pi x\,dx = -\frac{2n\sin p\pi}{\pi}\frac{1}{4n^2-p^2}$$

and make the summation

$$\sum_{n=1}^{\infty}\frac{1}{n}\int_0^1 \sin 2\pi nx\sin p\pi x\,dx = -\frac{2\sin p\pi}{\pi}\sum_{n=1}^{\infty}\frac{1}{4n^2-p^2}$$

We again assume that it is valid to interchange the order of the integration and summation

$$\sum_{n=1}^{\infty} \frac{1}{n} \int_0^1 \sin 2\pi nx \sin p\pi x \, dx = \int_0^1 \sum_{n=1}^{\infty} \frac{\sin 2\pi nx}{n} \sin p\pi x \, dx$$

and use the Fourier series

$$\sum_{n=1}^{\infty} \frac{\sin 2\pi nx}{n} = \frac{\pi}{2}(1-2x) \qquad 0 < x < 1$$

to obtain

$$\sum_{n=1}^{\infty} \frac{1}{n} \int_0^1 \sin 2\pi nx \sin p\pi x \, dx = \frac{\pi}{2} \int_0^1 (1-2x) \sin p\pi x \, dx$$

$$= \frac{\pi}{2} \left[ \frac{1+\cos p\pi}{p\pi} - \frac{2\sin p\pi}{(p\pi)^2} \right]$$

Hence we obtain

$$\sum_{n=1}^{\infty} \frac{1}{4n^2 - p^2} = \frac{\pi^2}{4} \left[ \frac{2}{(p\pi)^2} - \frac{1+\cos p\pi}{p\pi \sin p\pi} \right]$$

or equivalently

$$\sum_{n=1}^{\infty} \frac{1}{4n^2 - p^2} = \frac{1}{4p} \left[ \frac{2}{p} - \pi \frac{1+\cos p\pi}{\sin p\pi} \right]$$

which corresponds with (4.7).

$\square$

Referring to (6.1)

$$\int_0^1 \log(2\sin \pi x) \cos p\pi x \, dx = \frac{p}{\pi} \sum_{n=1}^{\infty} \frac{1}{n} \frac{\sin p\pi}{4n^2 - p^2}$$

we denote $f(p)$ as

$$f(p) = \int_0^1 \log(2\sin \pi x) \cos p\pi x \, dx$$

and take the limit as $p \to 0$ to obtain

$$\lim_{p \to 0} \frac{f(p)}{p \sin p\pi} = \frac{\varsigma(3)}{4\pi}$$

Since $f(0) = 0$ we may apply L'Hôpital's rule to get

$$\lim_{p \to 0} \frac{f(p)}{p \sin p\pi} = \lim_{p \to 0} \frac{f'(p)}{p\pi \cos p\pi + \sin p\pi}$$

We have

$$f'(p) = -\pi \int_0^1 x \log(2 \sin \pi x) \sin p\pi x \, dx$$

and since $f'(0) = 0$ we may apply L'Hôpital's rule again to obtain

$$\lim_{p \to 0} \frac{f(p)}{p \sin p\pi} = \lim_{p \to 0} \frac{f''(p)}{-p\pi^2 \sin p\pi + 2\pi \cos p\pi} = -\frac{f''(0)}{2\pi}$$

We have

$$f''(p) = -\pi^2 \int_0^1 x^2 \log(2 \sin \pi x) \cos p\pi x \, dx$$

and thus

$$f''(0) = -\pi^2 \int_0^1 x^2 \log(2 \sin \pi x) \, dx$$

Hence we obtain

(6.8)  $$\int_0^1 x^2 \log(2 \sin \pi x) \, dx = -\frac{\varsigma(3)}{2\pi^2}$$

in accordance with [30].

Employing the same procedure with (6.2) produces the same result.

It is clear that higher derivatives will produce integrals of the form $\int_0^1 x^{2N} \log(2 \sin \pi x) \, dx$.

**APPENDIX A**

**Some series for** $\log \Gamma(x)$

We have the Weierstrass canonical form of the gamma function [52, p.1]

(A.1)
$$\frac{1}{\Gamma(x)} = xe^{\gamma x} \prod_{n=1}^{\infty} \left\{ \left(1 + \frac{x}{n}\right) e^{-\frac{x}{n}} \right\}$$

and taking logarithms results in

(A.2)
$$\log \Gamma(x) = -\log x - \gamma x - \sum_{n=1}^{\infty} \left[ \log\left(1 + \frac{x}{n}\right) - \frac{x}{n} \right]$$

The following identity is easily derived by differentiating (A.2)

(A.3)
$$\psi(x) = \frac{\Gamma'(x)}{\Gamma(x)} = -\frac{1}{x} - \gamma - \sum_{n=1}^{\infty} \left( \frac{1}{x+n} - \frac{1}{n} \right)$$

Since [52, p.14] $\psi(1+x) = \psi(x) + \frac{1}{x}$ we obtain

$$\psi(1+x) + \gamma = -\sum_{n=1}^{\infty} \left( \frac{1}{x+n} - \frac{1}{n} \right)$$

and we easily see that

$$\psi(1+x) + \psi(1-x) + 2\gamma = -\sum_{n=1}^{\infty} \left( \frac{1}{n+x} + \frac{1}{n-x} - \frac{2}{n} \right)$$

or equivalently

(A.3)
$$\psi(1+x) + \psi(1-x) + 2\gamma = -2x^2 \sum_{n=1}^{\infty} \frac{1}{n} \frac{1}{n^2 - x^2}$$

This identity was employed in (2.22) above.

$\square$

We see that

$$\frac{\psi(1+x)+\psi(1-x)+2\gamma}{x^2} = -2\sum_{n=1}^{\infty}\frac{1}{n}\frac{1}{n^2-x^2}$$

$$\lim_{x\to 0}\frac{\psi(1+x)+\psi(1-x)+2\gamma}{x^2} = -2\varsigma(3)$$

and applying L'Hôpital's rule twice we obtain

$$= \lim_{x\to 0}\frac{\psi'(1+x)-\psi'(1-x)}{2x}$$

$$= \lim_{x\to 0}\frac{\psi''(1+x)+\psi''(1-x)}{2} = \frac{1}{2}\psi''(1)$$

This simply results in the well known identity

$$\psi''(1) = -2\varsigma(3)$$

□

Differentiating (A.3) results in

(A.4)   $$\psi'(1+x)-\psi'(1-x) = -4x\sum_{n=1}^{\infty}\frac{n}{\left(n^2-x^2\right)^2}$$

Integrating (A.3) results in for $x<1$

$$\log\Gamma(1+x)-\log\Gamma(1-x)+2\gamma x = -2\sum_{n=1}^{\infty}\left[\tanh^{-1}\frac{x}{n}-\frac{x}{n}\right]$$

$$= -\sum_{n=1}^{\infty}\left[\log\left(1+\frac{x}{n}\right)-\log\left(1+\frac{x}{n}\right)-\frac{2x}{n}\right]$$

(A.5)   $$\log\Gamma(1+x)-\log\Gamma(1-x)+2\gamma x = -\sum_{n=1}^{\infty}\left[\log\left(n+x\right)-\log\left(n-x\right)-\frac{2x}{n}\right]$$

We also have Euler's formula for the gamma function [52, p.2]

$$\Gamma(x) = \frac{1}{x}\prod_{n=1}^{\infty}\left[\left(1+\frac{1}{n}\right)^x\left(1+\frac{x}{n}\right)^{-1}\right]$$

and multiplying by $x$ and taking the logarithm of both sides gives us

$$(A.6) \qquad \log \Gamma(1+x) = \sum_{n=1}^{\infty} \left[ x \log\left(1 + \frac{1}{n}\right) - \log\left(1 + \frac{x}{n}\right) \right]$$

We then see that

$$\log \Gamma(1+x) - \log \Gamma(1-x) = \sum_{n=1}^{\infty} \left[ 2x \log\left(1 + \frac{1}{n}\right) - \log\left(1 + \frac{x}{n}\right) + \log\left(1 - \frac{x}{n}\right) \right]$$

$$= \sum_{n=1}^{\infty} \left[ 2x \log\left(1 + \frac{1}{n}\right) - \log\left(n + x\right) + \log\left(n - x\right) \right]$$

We have

$$\sum_{n=1}^{N} \left[ 2x \log\left(1 + \frac{1}{n}\right) - \log\left(n + x\right) + \log\left(n - x\right) \right]$$

$$= \sum_{n=1}^{N} \left[ 2x \log\left(1 + \frac{1}{n}\right) - \frac{2x}{n} - \log\left(n + x\right) + \log\left(n - x\right) + \frac{2x}{n} \right]$$

$$= 2x \sum_{n=1}^{N} \left[ \log\left(1 + \frac{1}{n}\right) - \frac{1}{n} \right] + \sum_{n=1}^{N} \left[ -\log\left(n + x\right) + \log\left(n - x\right) + \frac{2x}{n} \right]$$

As a result of telescoping we have

$$\sum_{n=1}^{N} \left[ \log\left(1 + \frac{1}{n}\right) - \frac{1}{n} \right] = \log(N+1) - H_N$$

and hence as $N \to \infty$ we rediscover (A.5).

We may also express (A.6) as

$$(A.7) \qquad \log \Gamma(x) = -\log x - \gamma x + \sum_{n=1}^{\infty} \left[ \log n - \log\left(n + x\right) + \frac{x}{n} \right]$$

and we note that this is equivalent to (A.2). We see that with $x \to x + a$

$$(A.8) \qquad \log \Gamma(x+a) = -\log(x+a) - \gamma(x+a) + \sum_{n=1}^{\infty} \left[ \log n - \log\left(n + a + x\right) + \frac{x+a}{n} \right]$$

Equation (A.8) was employed in deriving (1.1).

In passing, we note a similarity with equation (E.43b) in [22]

$$\sum_{n=1}^{\infty} \frac{H_n}{n^s} + \frac{1}{2}\varsigma'(s) - \frac{1}{2}\sum_{n=1}^{\infty} \frac{\log(n+1)}{n^s} - \gamma\varsigma(s) = \sum_{n=1}^{\infty}\left[\tanh^{-1}\frac{1}{n} - \frac{1}{n}\right]\left(\varsigma(s) - H_n^{(s)}\right)$$

Since $\lim_{n\to\infty}\left[\varsigma(s) - H_n^{(s)}\right] = 0$, we see that the factor $\varsigma(s) - H_n^{(s)}$ is necessary to ensure convergence of the series.

## APPENDIX B

**Fourier coefficients for** $\log x$

We have the Fourier series for $\log x$

$$\log x = \frac{1}{2}a_0 + \sum_{n=1}^{\infty}(a_n\cos nx + b_n\sin nx)$$

where, as proved below, we have

(B.1) $\qquad a_n = \frac{1}{\pi}\int_0^{2\pi}\log x.\cos nx\,dx = -\frac{1}{n\pi}\left[si(2n\pi) + \frac{\pi}{2}\right] = -\frac{Si(2n\pi)}{n\pi}$

(B.2) $\qquad b_n = \frac{1}{\pi}\int_0^{2\pi}\log x.\sin nx\,dx = \frac{1}{n\pi}\left[Ci(2n\pi) - \gamma - \log(2n\pi)\right]$

With integration by parts we get

$$\int_\alpha^x \log t.\sin at\,dt = -\log t\frac{\cos at}{a}\Big|_\alpha^x + \int_\alpha^x \frac{\cos at}{at}\,dt$$

$$= -\log t\frac{\cos at}{a}\Big|_\alpha^x + \int_\alpha^x \frac{\cos at - 1}{at}\,dt + \int_\alpha^x \frac{1}{at}\,dt$$

$$= -\log x\frac{\cos ax}{a} + \int_\alpha^{ax} \frac{\cos u - 1}{u}\,du + \frac{\log\alpha}{a}(\cos a\alpha - 1) + \frac{\log x}{a}$$

Therefore, in the limit as $\alpha \to 0$, we get

$$\int_0^x \log t. \sin at\, dt = -\log x \frac{\cos ax}{a} + \int_0^{ax} \frac{\cos u - 1}{u}\, du + \frac{\log x}{a}$$

and, as reported in Nielsen's book [44, p.12], we have

(B.3) $$\int_0^x \log t. \sin at\, dt = \frac{1}{a}\Big[ Ci(ax) - \gamma - \cos(ax)\log x - \log a \Big]$$

We therefore obtain

$$b_n = \frac{1}{\pi} \int_0^{2\pi} \log x. \sin nx\, dx = \frac{1}{n\pi}\Big[ Ci(2n\pi) - \gamma - \log(2n\pi) \Big]$$

Let us now consider the integral

$$\int_0^x \log t. \cos nt\, dt = \log t \frac{\sin nt}{n}\Big|_0^x - \int_0^x \frac{\sin nt}{nt}\, dt$$

We have

$$\lim_{t\to 0} \log t \sin nt = \lim_{t\to 0} t \log t \frac{\sin nt}{t} = \lim_{t\to 0} t \log t \lim_{t\to 0} \frac{\sin nt}{t} = 0$$

Therefore we get

$$\int_0^x \log t. \cos nt\, dt = \frac{\sin nx \log x}{n} - \int_0^x \frac{\sin nt}{nt}\, dt$$

$$= \frac{\sin nx \log x}{n} - \frac{1}{n}\int_0^{nx} \frac{\sin u}{u}\, du$$

and hence we get

$$\int_0^x \log t. \cos nt\, dt = \frac{\sin nx \log x}{n} - \frac{Si(nx)}{n}$$

We also have

$$a_0 = \frac{1}{\pi} \int_0^{2\pi} \log x\, dx = \log(2\pi) - 1$$

and making use of the familiar Fourier series

$$\frac{1}{2}(\pi - x) = \sum_{n=1}^{\infty} \frac{\sin nx}{n} \quad (0 < x < 2\pi)$$

$$\log\left[2\sin(x/2)\right] = -\sum_{n=1}^{\infty} \frac{\cos nx}{n} \quad (0 < x < 2\pi)$$

we therefore obtain the Fourier series (where we have made the substitution $x \rightarrow 2\pi x$)

(B.4) $\qquad 1 + \log x - \frac{1}{2}\log(2\sin \pi x) + (\gamma + \log 2\pi)\left(\frac{1}{2} - x\right) =$

$$\sum_{n=1}^{\infty} \frac{\left[Ci(2n\pi) - \log n\right]\sin 2n\pi x - si(2n\pi)\cos 2n\pi x}{n\pi}$$

In (B.4) we let $x \rightarrow 1 - x$ and combine the resulting identity with (B.4) to obtain

(B.5) $\qquad \log x + \log(1 - x) - \log(2\sin \pi x) + 2 = -\frac{2}{\pi}\sum_{n=1}^{\infty} \frac{si(2n\pi)\cos 2n\pi x}{n}$

and subtraction results in

(B.6) $\qquad \log x - \log(1 - x) + (\gamma + \log 2\pi)(1 - 2x) = 2\sum_{n=1}^{\infty} \frac{\left[Ci(2n\pi) - \log n\right]\sin 2n\pi x}{n\pi}$

We refer to (3.12) for Kummer's Fourier series for $\log \Gamma(x)$ which is valid for $0 < x < 1$

$$\log \Gamma(x) = \frac{1}{2}\log \frac{\pi}{\sin \pi x} + \left(\frac{1}{2} - x\right)(\gamma + \log 2\pi) + \frac{1}{\pi}\sum_{n=1}^{\infty} \frac{\log n}{n}\sin 2\pi nx$$

and combining this with (B.4) we obtain Nielsen's formula [45, p.79] for $0 < x < 1$

(B.7) $\qquad \log \Gamma(x) = \frac{1}{2}\log(2\pi) - 1 - \log x + \sum_{n=1}^{\infty} \frac{\sin 2n\pi x \, Ci(2n\pi) - \cos 2n\pi x \, si(2n\pi)}{n\pi}$

Integration gives us

(B.8)

$$\int_0^t \log \Gamma(x)\,dx = \frac{t}{2}\log(2\pi) - t \log t + \frac{1}{2\pi^2}\sum_{n=1}^{\infty} \frac{[1 - \cos 2n\pi t]Ci(2n\pi) - \sin 2n\pi t \, si(2n\pi)}{n^2}$$

and for example we have

$$\int_0^{1/2} \log \Gamma(x)\, dx = \frac{1}{4}\log(2\pi) + \frac{1}{2}\log 2 + \frac{1}{2\pi^2}\sum_{n=1}^{\infty}\frac{[1-\cos n\pi]\, Ci(2n\pi)}{n^2}$$

Integration of (B.6) results in

(B.9)

$$u\log u + (1-u)\log(1-u) + (\gamma + \log 2\pi)\left(u - u^2\right) = -\frac{1}{\pi^2}\sum_{n=1}^{\infty}\frac{\left[Ci(2n\pi) - \log n\right][\cos 2n\pi u - 1]}{n^2}$$

With $u = 1/2$ we obtain

$$\frac{1}{4}(\gamma + \log 2\pi) - \log 2 = -\frac{1}{\pi^2}\sum_{n=1}^{\infty}\frac{\left[Ci(2n\pi) - \log n\right][(-1)^n - 1]}{n^2}$$

and we therefore have

(B.10) $$\frac{1}{4}(\gamma + \log 2\pi) - \log 2 = \frac{2}{\pi^2}\sum_{n=0}^{\infty}\frac{\left[Ci[2(2n+1)\pi] - \log(2n+1)\right]}{(2n+1)^2}$$

We have

$$(1 - 2^{-s})\varsigma(s) = \sum_{n=0}^{\infty}\frac{1}{(2n+1)^s}$$

and therefore

$$(1 - 2^{-s})\varsigma'(s) + 2^{-s}\varsigma(s)\log 2 = -\sum_{n=0}^{\infty}\frac{\log(2n+1)}{(2n+1)^s}$$

We then obtain

$$\frac{1}{4}(\gamma + \log 2\pi) - \log 2 = \frac{2}{\pi^2}\sum_{n=0}^{\infty}\frac{Ci[2(2n+1)\pi]}{(2n+1)^2} + \frac{2}{\pi^2}\left[\frac{3}{4}\varsigma'(2) + \frac{1}{4}\varsigma(2)\log 2\right]$$

$$= \frac{2}{\pi^2}\sum_{n=0}^{\infty}\frac{Ci[2(2n+1)\pi]}{(2n+1)^2} + \frac{3\varsigma'(2)}{2\pi^2} + \frac{1}{12}\log 2$$

and this becomes

(B.11) 
$$\frac{2}{\pi^2}\sum_{n=0}^{\infty}\frac{Ci[2(2n+1)\pi]}{(2n+1)^2}=\frac{1}{4}-3\varsigma'(-1)-\frac{13}{12}\log 2$$

**APPENDIX C**

**Some aspects of the multiple gamma functions**

We note from [52, p.209] that

$$(C.1)\quad 2\int_0^z t\log\Gamma(a+t)\,dt=\left(\frac{1}{4}-\frac{1}{2}a+\frac{1}{2}a^2-2\log A\right)z+\left(\frac{1}{2}\log(2\pi)-\frac{1}{2}a+\frac{1}{4}\right)z^2$$

$$-\frac{1}{2}z^2+z^2\log\Gamma(z+a)-(a-1)^2\left[\log\Gamma(z+a)-\log\Gamma(a)\right]+(2a-3)\left[\log G(z+a)-\log G(a)\right]$$

$$-2\left[\log\Gamma_3(z+a)-\log\Gamma_3(a)\right]$$

(i) where $G(z)$ is the Barnes double gamma function

$$G(z+1)=(2\pi)^{z/2}\exp\left[-\frac{1}{2}\left(z^2+\gamma z^2+z\right)\right]\prod_{n=1}^{\infty}\left(1+\frac{z}{n}\right)^n e^{-z+z^2/2n}$$

(ii) where $\Gamma_3(z)$ is the triple gamma function defined in [52, p.42]

$$\Gamma_3(1+z)=\exp(c_1 z+c_2 z^2+c_3 z^3)\prod_{k=1}^{\infty}\left\{\left(1+\frac{z}{k}\right)^{-\frac{1}{2}k(k+1)}\exp\left[\left(1+\frac{1}{k}\right)\left(\frac{1}{2}kz-\frac{1}{4}z^2+\frac{1}{6k}z^3\right)\right]\right\}$$

(iii) where the constants are

$$c_1=\frac{3}{8}-\frac{1}{4}\log(2\pi)-\log A\ ,\qquad c_2=\frac{1}{4}\left[\gamma+\log(2\pi)+\frac{1}{2}\right]$$

$$c_3=-\frac{1}{6}\left[\gamma+\varsigma(2)+\frac{3}{2}\right]\qquad\log A=\frac{1}{12}-\varsigma'(-1)$$

(iv) where $A\cong 1.28242\cdots$ is the Glaisher-Kinkelin constant

$$\log A=\lim_{n\to\infty}\left[\sum_{k=1}^{n}k\log k-\left(\frac{n^2}{2}+\frac{n}{2}+\frac{1}{12}\right)\log n+\frac{n^2}{4}\right]$$

and

$$\log A = \frac{1}{12} - \varsigma'(-1)$$

With $a = z = 1$ we obtain the much simplified result

(C.2)     $$\int_0^1 x \log \Gamma(x)\, dx = \frac{1}{4}\log(2\pi) - \log A$$

which we used in (4.18). This integral may also be obtained by using (C.4).

Choi and Srivastava (see [17] and [52, p.209]) have evaluated a number of integrals of the form $\int_0^z x^n \psi(x+a)\, dx$ and, in the particular case, $a = z = 1$ we have

(C.3)     $$\int_0^1 x\, \psi(x+1)\, dx = 1 - \frac{1}{2}\log(2\pi)$$

$$\int_0^1 x^2 \psi(x+1)\, dx = \frac{1}{2} - \frac{1}{2}\log(2\pi) + 2\log A$$

$$\int_0^1 x^3 \psi(x+1)\, dx = \frac{1}{3} - \frac{7}{12}\log(2\pi) + 3\log A - 2\int_0^1 \log \Gamma_3(x+1)\, dx$$

where $\Gamma_3(x)$ is the triple gamma function [52, p.40].

Since for $n \geq 1$

$$\int_0^1 x^n \psi(x)\, dx = -\frac{1}{n} + \int_0^1 x^n \psi(x+1)\, dx$$

we immediately see that

(C.4)     $$\int_0^1 x\, \psi(x)\, dx = -\frac{1}{2}\log(2\pi)$$

$$\int_0^1 x^2 \psi(x)\, dx = -\frac{1}{2}\log(2\pi) + 2\log A$$

$$\int_0^1 x^3 \psi(x) dx = -\frac{7}{12}\log(2\pi) + 3\log A - 2\int_0^1 \log\Gamma_3(x+1) dx$$

Integration by parts gives us

$$\int_0^1 x^3 \psi(x) dx = -3\int_0^1 x^2 \log\Gamma(x) dx$$

and we obtain

$$\int_0^1 x^2 \log\Gamma(x) dx = \frac{7}{36}\log(2\pi) - \log A + \frac{2}{3}\int_0^1 \log\Gamma_3(x+1) dx$$

The following integral was given by Choi and Srivastava [17]

(C.5)     $$\int_0^1 \log\Gamma_3(x+1) dx = -\frac{1}{24}\log(2\pi) + \frac{3\varsigma(3)}{8\pi^2}$$

and this was also derived in [20] using an observation made by Glasser [32]. Therefore we obtain

(C.6)     $$\int_0^1 x^2 \log\Gamma(x) dx = \frac{1}{6}\log(2\pi) - \log A + \frac{\varsigma(3)}{4\pi^2}$$

as derived in a different manner in (2.39) above.

## APPENDIX D

## Table of relevant integrals and series

### Cosine function

(1.27)     $$\int_0^1 \log\Gamma(x+a)\cos 2k\pi x \, dx = -\frac{1}{2k\pi}\big[-\sin(2k\pi a)Ci(2k\pi a) + \cos(2k\pi a)si(2k\pi a)\big]$$

(1.90)     $$\int_0^1 \log\Gamma(x)\cos p\pi x \, dx = \left[\frac{1}{2}\log(2\pi) - 1\right]\frac{\sin p\pi}{p\pi} + \frac{Si(p\pi)}{p\pi}$$

$$+ \frac{2(1-\cos p\pi)}{\pi^2}\sum_{n=1}^{\infty}\frac{Ci(2n\pi)}{4n^2 - p^2} + \frac{p\sin p\pi}{\pi^2}\sum_{n=1}^{\infty}\frac{1}{n}\frac{si(2n\pi)}{4n^2 - p^2}$$

(1.92) $\displaystyle\int_0^1 \log\Gamma(x)\cos\pi x\,dx = \frac{Si(\pi)}{\pi} + \frac{4}{\pi^2}\sum_{n=1}^{\infty}\frac{Ci(2n\pi)}{4n^2-1}$

(1.93) $\displaystyle\int_0^1 \log\Gamma(x)\cos(2k+1)\pi x\,dx = \frac{Si[(2k+1)\pi]}{(2k+1)\pi} + \frac{4}{\pi^2}\sum_{n=1}^{\infty}\frac{Ci(2n\pi)}{4n^2-(2k+1)^2}$

(4.5) $\displaystyle\int_0^1 \log\Gamma(x)\cos p\pi x\,dx$

$$= \frac{1}{2}\log(2\pi)\frac{\sin p\pi}{p\pi} - \frac{p\sin p\pi}{2\pi}\sum_{n=1}^{\infty}\frac{1}{n}\frac{1}{4n^2-p^2} + \frac{2(1-\cos p\pi)}{\pi^2}\sum_{n=1}^{\infty}\frac{\gamma+\log(2\pi n)}{4n^2-p^2}$$

(4.5.2) $\displaystyle\int_0^1 \log\Gamma(x)\cos\pi x\,dx = \frac{4}{\pi^2}\sum_{n=1}^{\infty}\frac{\gamma+\log(2\pi n)}{4n^2-1}$

(4.9) $\displaystyle\int_0^1 \log\Gamma(x)\cos(2k+1)\pi x\,dx = \frac{2}{\pi^2}\left[\frac{\log(2\pi)+\gamma}{(2k+1)^2} + 2\sum_{n=1}^{\infty}\frac{\log n}{4n^2-(2k+1)^2}\right]$

(2.28) $\displaystyle\int_0^1 \log\Gamma(x)\cos 2k\pi x\,dx = \frac{1}{4k}$

(4.10) $\displaystyle\int_0^1 \log\Gamma(x)\cos(\pi x/2)\,dx = \frac{1}{\pi}\log(2\pi) - \frac{3\log 2 - 2}{\pi} + \frac{[\gamma+\log(2\pi)](4-\pi)}{\pi^2}$

$$+ \frac{8}{\pi^2}\sum_{n=1}^{\infty}\frac{\log n}{16n^2-1}$$

(1.28.1) $\displaystyle\int_0^1 \log\Gamma\left(x+\frac{1}{2}\right)\cos 2k\pi x\,dx = \frac{(-1)^{k+1}}{2k\pi}si(k\pi)$

Differentiating (5.5) results in

$$\pi\int_0^1 x\log\Gamma(x)\cos p\pi x\,dx$$

$$= \frac{1}{2}\log(2\pi)\frac{p\pi\sin p\pi - 1 + \cos p\pi}{p^2\pi} - \frac{p^2(1-\cos p\pi)}{\pi}\sum_{n=1}^{\infty}\frac{1}{n}\frac{1}{(4n^2-p^2)^2}$$

$$- \frac{1-\cos p\pi + p\pi\sin p\pi}{2\pi}\sum_{n=1}^{\infty}\frac{1}{n}\frac{1}{4n^2-p^2} - \frac{4p\sin p\pi}{\pi^2}\sum_{n=1}^{\infty}\frac{\gamma+\log(2\pi n)}{(4n^2-p^2)^2}$$

$$- \frac{2\cos p\pi}{\pi}\sum_{n=1}^{\infty}\frac{\gamma+\log(2\pi n)}{4n^2-p^2}$$

$$\pi\int_0^1 x\log\Gamma(x)\cos\pi x\,dx$$

$$= -\frac{1}{\pi}\log(2\pi) - \frac{2}{\pi}\sum_{n=1}^{\infty}\frac{1}{n}\frac{1}{(4n^2-1)^2} - \frac{1}{\pi}\sum_{n=1}^{\infty}\frac{1}{n}\frac{1}{4n^2-1} + \frac{2}{\pi}\sum_{n=1}^{\infty}\frac{\gamma+\log(2\pi n)}{4n^2-1}$$

(4.17) $$\frac{\pi^2}{2}\int_0^1\log^2\Gamma(x)\cos\pi x\,dx = \log(2\pi) + \gamma + \left[\log\frac{\pi}{2}+1\right]\sum_{n=1}^{\infty}\frac{1}{2n+1}\log\frac{n+1}{n}$$

$$+ (\gamma+\log 2\pi)\left[\log\frac{\pi}{2}+1\right] + \sum_{n=1}^{\infty}\frac{1}{2n+1}\left[\frac{1}{2n+1}+2\sum_{j=0}^{n-1}\frac{1}{2j+1}\right]\log\frac{n+1}{n}$$

(1.78) $$\int_0^1\psi(x+1)\cos 2n\pi x\,dx = Ci(2n\pi)$$

## Sine function

(1.1) $$\int_0^1\log\Gamma(x+a)\sin 2k\pi x\,dx = -\frac{1}{2k\pi}\left[\log a - \cos(2k\pi a)\,Ci(2k\pi a) - \sin(2k\pi a)\,si(2k\pi a)\right]$$

(1.97) $$\int_0^1\log\Gamma(x)\sin p\pi x\,dx = \left[\frac{1}{2}\log(2\pi)-1\right]\frac{1-\cos p\pi}{p\pi} - \frac{Ci(p\pi)}{p\pi} + \frac{\gamma+\log(p\pi)}{p\pi}$$

$$- \frac{2\sin p\pi}{\pi^2}\sum_{n=1}^{\infty}\frac{Ci(2n\pi)}{4n^2-p^2} + \frac{p(1-\cos p\pi)}{\pi^2}\sum_{n=1}^{\infty}\frac{1}{n}\frac{si(2n\pi)}{4n^2-p^2}$$

(1.98) $$\int_0^1\log\Gamma(x)\sin\pi x\,dx = \left[\frac{1}{2}\log(2\pi)-1\right]\frac{2}{\pi} - \frac{Ci(\pi)}{\pi} + \frac{\gamma+\log\pi}{\pi} + \frac{2}{\pi^2}\sum_{n=1}^{\infty}\frac{1}{n}\frac{si(2n\pi)}{4n^2-1}$$

(1.100)

$$\int_0^1 \log \Gamma(x) \sin(2k+1)\pi x \, dx = \left[\frac{1}{2}\log(2\pi) - 1\right]\frac{2}{(2k+1)\pi} - \frac{Ci[(2k+1)\pi]}{(2k+1)\pi} + \frac{\gamma + \log[(2k+1)\pi]}{(2k+1)\pi}$$

$$+ \frac{2(2k+1)}{\pi^2}\sum_{n=1}^{\infty}\frac{1}{n}\frac{si(2n\pi)}{4n^2 - (2k+1)^2}$$

(5.5) $\displaystyle\int_0^1 \log \Gamma(x)\sin p\pi x \, dx$

$$= \frac{1}{2}\log(2\pi)\frac{1 - \cos p\pi}{p\pi} - \frac{p(1 - \cos p\pi)}{2\pi}\sum_{n=1}^{\infty}\frac{1}{n}\frac{1}{4n^2 - p^2} - \frac{2\sin p\pi}{\pi^2}\sum_{n=1}^{\infty}\frac{\gamma + \log(2\pi n)}{4n^2 - p^2}$$

(5.6) $\displaystyle\int_0^1 \log \Gamma(x)\sin(2k+1)\pi x \, dx = \frac{1}{(2k+1)\pi}\left[\log\left(\frac{\pi}{2}\right) + \frac{1}{2k+1} + 2\sum_{j=0}^{k-1}\frac{1}{2j+1}\right]$

(5.7) $\displaystyle\int_0^1 \log \Gamma(x)\sin \pi x \, dx = \frac{1}{\pi}\left[\log\frac{\pi}{2} + 1\right]$

(1.5) $\displaystyle\int_0^1 \log \Gamma(x)\sin 2k\pi x \, dx = \frac{\gamma + \log(2\pi k)}{2\pi k}$

(1.10) $\displaystyle\int_0^1 \log \Gamma\left(x + \frac{1}{2}\right)\sin 2k\pi x \, dx = \frac{1}{2k\pi}\left[\log 2 + (-1)^k Ci(k\pi)\right]$

(5.5.1) $\displaystyle\int_0^1 \log \Gamma(x)\sin(\pi x / 2) \, dx = \frac{1}{\pi}\log(2\pi) - \frac{1}{\pi}\sum_{n=1}^{\infty}\frac{1}{n}\frac{1}{16n^2 - 1} - \frac{8}{\pi^2}\sum_{n=1}^{\infty}\frac{\gamma + \log(2\pi n)}{16n^2 - 1}$

(4.19) $\displaystyle -\pi\int_0^1 x\log \Gamma(x)\sin p\pi x \, dx$

$$= \frac{1}{2}\log(2\pi)\frac{p\pi\cos p\pi - \sin p\pi}{p^2\pi} - \frac{p^2\sin p\pi}{\pi}\sum_{n=1}^{\infty}\frac{1}{n}\frac{1}{(4n^2 - p^2)^2}$$

$$- \frac{p\pi\cos p\pi + \sin p\pi}{2\pi}\sum_{n=1}^{\infty}\frac{1}{n}\frac{1}{4n^2 - p^2} + \frac{4p(1 - \cos p\pi)}{\pi^2}\sum_{n=1}^{\infty}\frac{\gamma + \log(2\pi n)}{(4n^2 - p^2)^2}$$

$$+\frac{2\sin p\pi}{\pi}\sum_{n=1}^{\infty}\frac{\gamma+\log(2\pi n)}{4n^2-p^2}$$

(4.20)

$$\pi\int_0^1 x\log\Gamma(x)\sin\pi x\,dx=\frac{1}{2}\log(2\pi)-\frac{1}{2}(2\log 2-1)-\frac{[\gamma+\log(2\pi)]}{\pi^2}\frac{\pi^2-8}{2}-\frac{8}{\pi^2}\sum_{n=1}^{\infty}\frac{\log n}{(4n^2-1)^2}$$

$$-\pi^2\int_0^1 x^2\log\Gamma(x)\sin p\pi x\,dx$$

$$=\frac{1}{2}\log(2\pi)\frac{p^3\pi^3\cos p\pi-(p\pi\sin p\pi-1+\cos p\pi)2p\pi}{p^4\pi^2}$$

$$-\frac{4p^3(1-\cos p\pi)}{\pi}\sum_{n=1}^{\infty}\frac{1}{n}\frac{1}{(4n^2-p^2)^3}-\frac{-2p+p^2\pi\sin p\pi-2p\cos p\pi}{\pi}\sum_{n=1}^{\infty}\frac{1}{n}\frac{1}{(4n^2-p^2)^2}$$

$$-\frac{p(1-\cos p\pi+p\pi\sin p\pi)}{\pi}\sum_{n=1}^{\infty}\frac{1}{n}\frac{1}{(4n^2-p^2)^2}-\frac{p\pi^2\cos p\pi+2p\sin p\pi}{2\pi}\sum_{n=1}^{\infty}\frac{1}{n}\frac{1}{4n^2-p^2}$$

$$-\frac{16p^2\sin p\pi}{\pi^2}\sum_{n=1}^{\infty}\frac{\gamma+\log(2\pi n)}{(4n^2-p^2)^3}-\frac{4p\pi\cos p\pi+4\sin p\pi}{\pi^2}\sum_{n=1}^{\infty}\frac{\gamma+\log(2\pi n)}{(4n^2-p^2)^2}$$

$$-\frac{4p\cos p\pi}{\pi}\sum_{n=1}^{\infty}\frac{\gamma+\log(2\pi n)}{(4n^2-p^2)^2}+\frac{2\pi\sin p\pi}{\pi}\sum_{n=1}^{\infty}\frac{\gamma+\log(2\pi n)}{4n^2-p^2}$$

$$-\pi^2\int_0^1 x^2\log\Gamma(x)\sin\pi x\,dx$$

$$=\frac{1}{2}\log(2\pi)\frac{4-\pi^2}{\pi}-\frac{8}{\pi}\sum_{n=1}^{\infty}\frac{1}{n}\frac{1}{(4n^2-1)^3}-\frac{2}{\pi}\sum_{n=1}^{\infty}\frac{1}{n}\frac{1}{(4n^2-1)^2}+\frac{\pi}{2}\sum_{n=1}^{\infty}\frac{1}{n}\frac{1}{4n^2-1}$$

$$+\frac{4}{\pi}\sum_{n=1}^{\infty}\frac{\gamma+\log(2\pi n)}{(4n^2-1)^2}+\frac{4}{\pi}\sum_{n=1}^{\infty}\frac{\gamma+\log(2\pi n)}{(4n^2-1)^2}$$

(4.12) $\quad\displaystyle\int_0^1\psi(x)\sin\pi x\,dx=-\frac{2}{\pi}\left[\log(2\pi)+\gamma+2\sum_{n=1}^{\infty}\frac{\log n}{4n^2-1}\right]$

$$(5.9) \qquad \int_0^1 \psi(x)\sin^2 2\pi x \, dx = -\frac{\gamma + \log 2\pi}{2}$$

$$(1.77) \qquad \int_0^1 \psi(x+1)\sin 2n\pi x \, dx = si(2n\pi)$$

## ACKNOWLEDGEMENTS

I thank R.J. Hughes for drawing my attention to the identity (A.3).

Princeton, NJ, 2003.

Donal F. Connon
Elmhurst
Dundle Road
Matfield
Kent TN12 7HD
dconnon@btopenworld.com